\documentclass[12pt]{amsart}

\usepackage{amsmath,amssymb,amsfonts, amsthm, amsrefs}
\usepackage[margin=2cm]{geometry}
\usepackage{verbatim}
\usepackage{tikz-cd}
\usepackage{freetikz}
\usetikzlibrary{calc}

\theoremstyle{plain}

\newtheorem{prop}{Proposition}[section]

\newtheorem{thm}[prop]{Theorem}

\newtheorem{cor}[prop]{Corollary}
\newtheorem*{cor*}{Corollary}

\newtheorem*{eg}{Example}

\newtheorem{defn}[prop]{Definition}         
\newtheorem{notation}[prop]{Notation}         

\numberwithin{equation}{section}

\usepackage{tikz}
\usetikzlibrary{shapes}
\usetikzlibrary{decorations,decorations.pathreplacing}

\newcommand\uchair{}
 \def\uchair[#1](#2)(#3:#4)
 { \draw[#1] ($(#2)$) --  ($(#2)+(-#3,#3)$) --  ($(#2)+(-#3+1,#3+1)$) -- ($(#2)+(-#3+1-#4,#3+1+#4)$); 
   \draw[#1]  ($(#2)+(-#3+0.5,#3+0.5)$) -- ($(#2)+(-#3+0.5-#4,#3+0.5+#4)$);
   \draw[#1] ($(#2)+(-#3,#3)$) -- ($(#2)+(-#3-#4,#3+#4)$);
 }

\newcommand\chair{}
 \def\chair[#1](#2)(#3:#4)
 { \draw[#1] ($(#2)$) --  ($(#2)+(-#3,-#3)$) --  ($(#2)+(-#3+1,-#3-1)$) -- ($(#2)+(-#3+1-#4,-#3-1-#4)$); 
   \draw[#1]  ($(#2)+(-#3+0.5,-#3-0.5)$) -- ($(#2)+(-#3+0.5-#4,-#3-0.5-#4)$);
   \draw[#1] ($(#2)+(-#3,-#3)$) -- ($(#2)+(-#3-#4,-#3-#4)$);
 }

\begin{document}
\title{Conway Rational Tangles and the Thompson Group}

\author{Ariana Grymski}
\address{Department of Mathematics and Statistics, Loyola University Chicago}
\email{agrymski@luc.edu}

\author{Emily Peters}
\address{Department of Mathematics and Statistics, Loyola University Chicago}
\email{epeters3@luc.edu}

\begin{abstract}
There is a map, defined and studied by Jones, from Thompson's group $F$ to knots.  Jones proved that every knot is in the image of this map -- that is, that every knot can be seen as the ``knot closure" of a Thompson group element.
We approach the question of methodologically finding Thompson group elements to generate a particular knot or link through the lens of Conway's rational tangles.  We are able to give methods to construct any product or concatenation of simple tangles, and we hope these are seeds for a more skein-theoretic approach to the construction question.
\end{abstract}

\maketitle

\section{Introduction}

Thompson's groups $F \subset T \subset V$ are infinite but approachable groups which exhibit some particularly interesting characteristics \cite{MR1426438}.  Most famously is the open question of whether $F$ is amenable.  They have applications in a variety of fields including logic, von Neumann algebras, category and homotopy theory, and shape theory \cite{MR2856142}.

Recently V.F.R. Jones has considered connections between the Thompson groups and constructions of conformal field theories \cite{MR3589908}.  The idea of this work is to replace the very large group $\text{Diff}(\mathbb{S}^1)$ with the more approachable Thompson group $T$, which acts on the circle via piecewise linear functions which are differentiable except at finitely many dyadic rational numbers.  Sadly, one cannot achieve a CFT as a semicontinuous limit of finite spin chains with an action of $T$, as shown in \cite{MR3764571}.  Nevertheless, this work has been the catalyst for recent interesting work on the Thompson group, including study of it representations \cites{MR4156130, MR4000616, MR3989149, MR4294116} and a connection to knots, which is the subject of this paper.

One presentation of the Thompson group is by pairs of rooted binary trees \cite{MR1426438}.  Jones observed that these pairs of diagrams can be easily turned into knot/link diagrams \cite{MR3986040}, and proved that all knots and links can be thus generated.  However, this method of producing knots is somewhat piecemeal, and the question of whether there is a better way to do this was asked by Jones himself.  

In this article, we take a tangle-theoretic approach to the creation of knots and links from Thompson diagrams, using Conway's rational tangles \cite{MR0258014} as an organizational principle.  We explain how to construct Thompson group elements that correspond to tangle sums or concatenations of (positive) $n$-crossing tangles.   

After this work was completed, we were made aware of the article \cite{MR4403999} in which Aiello and Baader construct Thompson group elements corresponding to a similar class of knots (this is their Corollary 2, concerning arborescent links associated to certain trees with positive weights).  Their algorithm, and the Thompson group elements they construct, are different from ours.  From the point of view of the Thompson group, their results are more interesting, as they answer the question of which links come from the positive monoid $F_+$ of $F$.  But our results may be more interesting from a skein-theoretic point of view, and perhaps generalize in different directions.  Regardless, there is now a large class of knots for which three distinct presentations are known, which could be fruitful fodder in a search for a Markov-type theorem for the group $F$.  Section \ref{future} includes an explanation of this question.

Both authors were partially supported by NSF DMS grant 1501116. 
Emily Peters was introduced to this problem by Vaughan Jones and would like to thank him for interesting conversations on the topic.  We also want to thank Valeriano Aiello for bringing his work to our attention, and a thoughtful conversation around the similarities and differences of our results.  Thanks also to Matt Brin and Arnaud Brothier for feedback on the preprint.

\section{Background}

\subsection{Thompson's group $F$}  
An element of the Thompson group $F$ is a piecewise linear function, with slopes which are powers of 2 and only change at dyadic rational numbers.  
An example element is shown in graph form below.
The function is completely described by the domain and range partitions created by the non-differentiable points, shown in the second diagram below.
 The binary trees, also in the second picture, are constructed from the partition \cite{MR1426438}.
 If we join the two trees along their leaves, we get a Thompson diagram (third diagram below).  
 We finally turn this Thompson diagram into a knot \cite{MR3986040}, as shown in the fourth diagram.

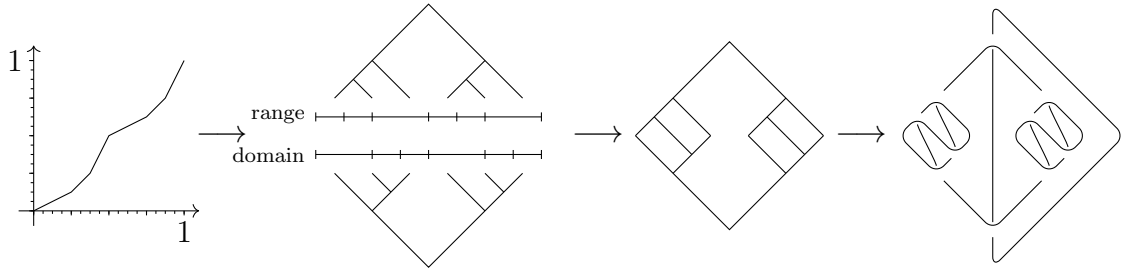
\begin{center}
\begin{tikzpicture}[scale=.5]
\begin{scope}[scale=4]
  \draw[->] (-0.1,0) -- (1.1,0);
  \draw[->] (0,-.1) -- (0,1.1);
 \foreach \x in {0,.0625,...,1} \draw (\x,-.02)--(\x,0); 
  \foreach \x in {0,.25,...,1} \draw (\x,-.03)--(\x,0); 
  \node at (1,0) [below] {$1$};
 \foreach \y in {0,.0625,...,1} \draw (-.02,\y)--(0,\y); 
  \foreach \y in {0,.25,...,1} \draw (-.03,\y)--(0,\y); 
  \node at (0,1) [left] {$1$};
  \draw (0,0) -- (1/4,1/8)--(3/8,1/4)--(1/2,1/2)--(6/8,5/8)--(7/8,3/4)--(1,1);
  \end{scope}
 \node at (5,2) {$\longrightarrow$};
  \begin{scope}[shift={(8,-1)}]
  \begin{scope}[shift={(-.5,2.5)}, scale=6]
     \node at (0,0) [left] {\tiny{domain}};
    \draw (0,0) -- (1,0);
 \foreach \x in {0,.25,.375,.5,.75,.875,1} \draw (\x,-.02)--(\x,.02); 
 \end{scope}
  \begin{scope}[shift={(-.5,3.5)},scale=6]
       \node at (0,0) [left] {\tiny{range}};
    \draw (0,0) -- (1,0);
 \foreach \x in {0,.125, .25,.5,.625,.75,1} \draw (\x,-.02)--(\x,.02); 
  \end{scope}
\begin{scope}[shift={(0,2)},yscale=-1]
	\draw(0,0)--(2.5,2.5)--(5,0);
	\draw(1,0)--(1.5,.5);
	\draw(2,0)--(1,1);
	\draw(3,0)--(4,1);
	\draw(4,0)--(4.5,.5);
\end{scope}
\begin{scope}[shift={(0,4)}]
	\draw(0,0)--(2.5,2.5)--(5,0);
	\draw(1,0)--(.5,.5);
	\draw(2,0)--(1,1);
	\draw(3,0)--(4,1);
	\draw(4,0)--(3.5,.5);
\end{scope}
\end{scope}
 \node at (15,2) {$\longrightarrow$};
  \begin{scope}[shift={(16,-1)}]
\begin{scope}[shift={(0,3)},yscale=-1]
	\draw(0,0)--(2.5,2.5)--(5,0);
	\draw(1,0)--(1.5,.5);
	\draw(2,0)--(1,1);
	\draw(3,0)--(4,1);
	\draw(4,0)--(4.5,.5);
\end{scope}
\begin{scope}[shift={(0,3)}]
	\draw(0,0)--(2.5,2.5)--(5,0);
	\draw(1,0)--(.5,.5);
	\draw(2,0)--(1,1);
	\draw(3,0)--(4,1);
	\draw(4,0)--(3.5,.5);
\end{scope}
\end{scope}
 \node at (22,2) {$\longrightarrow$};
\begin{scope}[shift={(23,2)}]
	\draw[rounded corners] (.6,.6)--(1,1)--(2,0)--(1.5,-.5)--(.5,.5)--(0,0)--(1,-1)--(1.4,-.6);
	\draw[rounded corners] (1.5,-.3) .. controls +(90:1mm) and +(-90:1mm) .. (1,.8);
	\draw[rounded corners] (1,-.8) .. controls +(90:1mm) and +(-90:1mm) .. (.5,.3);
	\draw[rounded corners] (3.6,.6)--(4,1)--(5,0)--(4.5,-.5)--(3.5,.5)--(3,0)--(4,-1)--(4.4,-.6);
	\draw[rounded corners] (4.5,-.3) .. controls +(90:1mm) and +(-90:1mm) .. (4,.8);
	\draw[rounded corners] (4,-.8) .. controls +(90:1mm) and +(-90:1mm) .. (3.5,.3);
	\draw[rounded corners] (1.2,1.2)--(2.5,2.5)--(3.8,1.2);
	\draw[rounded corners] (1.2,-1.2)--(2.5,-2.5)--(3.8,-1.2);
	\draw (2.5,2.3)--(2.5,-2.3);
	\draw[rounded corners]  (2.5, 2.7) --(2.5,3.5)--(6,0)--(2.5,-3.5)--(2.5,-2.7);
\end{scope}
\end{tikzpicture}
\end{center}

\begin{defn}
$R$ denotes reflection over the horizontal axis, and acts on Thompson diagrams.
\end{defn}

\begin{defn}
The building blocks of our Thompson diagrams will be chairs.  \emph{Chair} always refers to the left diagram below, with the given orientation.  A \emph{reflected chair} is shown on the right.
$$
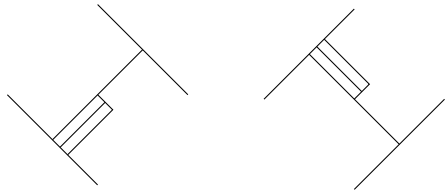
\hspace{12pt}
\begin{tikzpicture}[scale=0.6,baseline=-70pt]
\begin{scope}[xscale=-1,yscale=1]
\draw(4,0)--(2,-2);
\draw(6,-2)--(4,-4);
\draw(3,-1)--(5,-3);
\draw(4,-2)--(3.666,-2.333)--(4.666,-3.333);
\draw(3.833,-2.166)--(4.833,-3.166);
\end{scope}
\end{tikzpicture}
\hspace*{1cm}
\begin{tikzpicture}[rotate=180,scale=0.6,baseline]
\draw(4,0)--(2,-2);
\draw(6,-2)--(4,-4);
\draw(3,-1)--(5,-3);
\draw(4,-2)--(3.666,-2.333)--(4.666,-3.333);
\draw(3.833,-2.166)--(4.833,-3.166);
\end{tikzpicture}
$$
\end{defn}

\begin{notation}\label{mchairs}
We abbreviate the Thompson diagram
consisting of a rectangle with $n$ chairs inside with a single chair labelled with $n$.

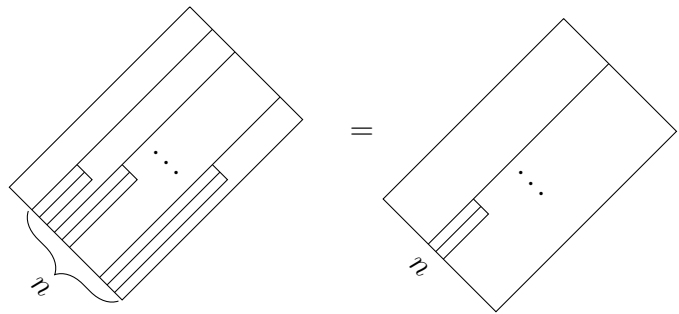
$$
\begin{tikzpicture}[scale=0.6,rotate=-90,baseline=-50]
\draw (0,0)--(2.5,2.5)--(6.5,-1.5)--(4,-4)--(0,0);

\draw(1,1)--(5,-3);
\draw(3.5,-1.5)--(3.833,-1.166)--(5.333,-2.666);
\draw(3.666,-1.333)--(5.166,-2.833);

\begin{scope}[shift={(1,1)}]
\draw(1,1)--(5,-3);
\draw(2.5,-0.5)--(2.833,-0.166)--(5.333,-2.666);
\draw(2.666,-0.333)--(5.166,-2.833);
\end{scope}

\begin{scope}[shift={(-0.5,-0.5)}]
\draw(1,1)--(5,-3);
\draw(4,-2)--(4.333,-1.666)--(5.333,-2.666);
\draw(4.166,-1.833)--(5.166,-2.833);
\end{scope}
\draw [decorate,decoration={brace,amplitude=10pt},xshift=-5pt,yshift=-8pt] (6.7,-1.3)--(4.7,-3.3);
\node at (6.2,-3.3)[rotate=-45]{$n$};
\node at (3.5,-0.5)[rotate=-45]{$\cdots$};

\end{tikzpicture}
\hspace{0.5cm}
=
\begin{tikzpicture}[rotate=-90,scale=0.6, baseline=-1.6cm]
\draw (0,0)--(2.5,2.5)--(6.5,-1.5)--(4,-4)--(0,0);

\draw(1,1)--(5,-3);
\draw(4,-2)--(4.333,-1.666)--(5.333,-2.666);
\draw(4.166,-1.833)--(5.166,-2.833);

\node at (5.5,-3.2)[rotate=-45]{$n$};
\node at (3.7,-0.7)[rotate=-45]{$\cdots$};

\end{tikzpicture}
$$
\end{notation}

Note that above, $n$ refers to the number of chairs, not the number of legs on chairs.  (All our chairs are three-legged.)

\begin{defn}
A \emph{snipped Thompson diagram} $\tilde{T}$ is a Thompson diagram $T$ with a small neighborhood of one (or more) of its corners deleted.  Subscripts right (r), left(l), top (t), or bottom (b) will be added to $\tilde{T}$ to indicate which corners are removed.
\end{defn}

Here is an example of $\tilde{T}_{ud}$:

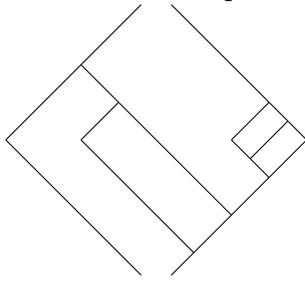
\begin{tikzpicture}
\draw(0,0)--(1.8,1.8);
\draw(2.2,1.8)--(4,0);
\draw(1,1)--(2,0);
\draw(1.5,0.5)--(1,0);
\draw(3.5,0.5)--(3,0);
\draw(3.75,0.25)--(3.25,-0.25);
\begin{scope}[xscale=1,yscale=-1]

\draw(0,0)--(1.8,1.8);
\draw(2.2,1.8)--(4,0);
\draw(2.5,1.5)--(1,0);
\draw(3,1)--(2,0);
\draw(3.5,0.5)--(3,0);
\end{scope}

\end{tikzpicture}

\subsection{Knots and tangles}

A mathematical knot is a smooth embedding $\mathbb{S}^1 \rightarrow \mathbb{R}^3$, considered up to equivalence by ambient isotopy of  $\mathbb{R}^3$.  Pragmatically speaking, knots are almost two dimensional -- one may choose a plane to project into so that the projection map is one-to-one except for a finite set of crossings where two strands come together.  Thus we illustrate knots in the plane by drawing crossings with a continuous overstand and drawing the understand with a break in it.  

To analyze the knots generated by Thompson diagrams, we are going to think about tangles as the building blocks for knots.

\begin{defn}

A \emph{tangle} is a collection of potentially intertwined strings inside a disk, with the ends of the strings anchored at the boundary of the disk.  To specify that a tangle has $n$ fixed strings at the boundary, we call it an $n$-tangle.  

\end{defn}

\begin{defn}
\emph{R} denotes reflection through a horizontal mirror, and can act on knots and tangles.
\end{defn}

%
%
%
%
%
%
%
%
%
%
%
%
%
%
%
%
%
%
%

\begin{defn}\label{simple}

The simplest tangles consist of two strings that we twist together $n$ times.  The tangle anchors are points N, S, E and W and our strings originate from N and W.  If $n$ is positive, illustrated below, we twist the strings $n$ times with the horizontal strand passing over the vertical strand.  If $n$ is negative, we twist the strings $-n$ times with the vertical strand passing over the horizontal.  This is called an \emph{$n$-crossing-tangle} and denoted by a circle labeled by $n$.

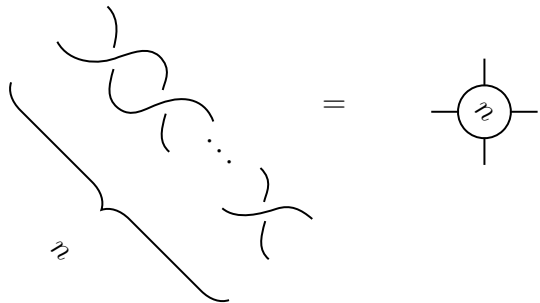
\tikzset{every picture/.style={line width=0.75pt}} 
$$
\begin{tikzpicture}[x=0.75pt,y=0.75pt,yscale=-1,xscale=1, baseline=-150]

\draw    (236.88,167.65) .. controls (237.96,163.13) and (240.51,153.31) .. (233.68,146.8) ;
\draw    (261.56,188.96) .. controls (262.64,184.45) and (266.58,178.34) .. (259.76,171.83) .. controls (253.91,166.82) and (246.31,170.15) .. (239.31,172.33) .. controls (229.4,176.6) and (214.4,175.6) .. (208.4,164.6) ;
\draw    (287.09,204.84) .. controls (283.52,195.02) and (274.25,190.55) .. (260.24,195.73) .. controls (252.81,198.16) and (245.86,203.95) .. (238.53,197.64) .. controls (233.34,193.06) and (234.15,187.21) .. (236.76,178.38) ;
\draw    (264.86,220.17) .. controls (259.68,215.59) and (260.48,209.74) .. (263.09,200.91) ;
\draw    (312.66,245.46) .. controls (313.73,240.95) and (317.68,234.84) .. (310.85,228.33) ;
\draw    (336.99,254.06) .. controls (331.14,249.05) and (326.04,246.71) .. (319.04,248.9) .. controls (311.61,251.32) and (298.93,254.96) .. (291.6,248.65) ;
\draw    (315.05,274.39) .. controls (309.86,269.81) and (310.67,263.96) .. (313.28,255.12) ;

\draw (285.99,208.62) node [anchor=north west][inner sep=0.75pt]  [rotate=-46.71]  {$\cdots $};

\draw[decorate,decoration= {brace,amplitude=10pt},rotate=180,shift={(-290,-290)}] (-5,-5)--(105,105);
\node at (210,270)[rotate=-45]{$n$};

\end{tikzpicture}
=
\hspace{1cm}
\begin{tikzpicture}[rotate=-45,baseline,scale=0.5]
\draw(-1,-1)--(1,1);
\draw(1,-1)--(-1,1);
\fill[color=white] (0,0) circle (20pt);
\draw (0,0) circle (20pt);
\node at (0,0) [rotate=-45]{$n$};
\end{tikzpicture}
$$

\end{defn}

\begin{defn}
\emph{Multiplication} is a tangle operation that attaches the upper and left strings of a reflected tangle to the lower and right strings of a second tangle, so that the product $F \cdot G$ is the tangle
$$
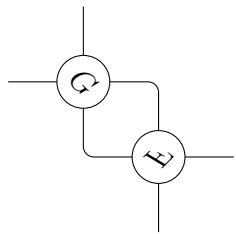
\begin{tikzpicture}

\draw[rounded corners](0,0)--(0,1)--(-1,1);
\draw[rounded corners](0,0)--(-1,0)--(-1,1);
\draw(1,0)--(0,0)--(0,-1);
\draw(-2,1)--(-1,1)--(-1,2);

\filldraw[fill=white](0,0) circle(10pt);
\node at (0,0)[rotate=-135,xscale=-1,yscale=1]{$F$};

\filldraw[fill=white](-1,1) circle(10pt);
\node at (-1,1)[rotate=-45]{$G$};

\end{tikzpicture}
$$

\end{defn}

Notice that tangle multiplication is not associative nor is it commutative.  Also, our multiplication is a rotation of the usual definition of tangle multiplication \cite{MR0258014}.

%
%
%
%
%
%
%
%
%
%
%
%
%
%
%

\begin{defn}
A \emph{rational tangle} is any tangle that can be formed by multiplication of $n$-crossing-tangles.

\end{defn}

\begin{defn}
\emph{Addition} is a tangle operation that attaches the right most strands of a tangle to the left-most strands of another tangle, so that the sum $F+G$ is the tangle

\end{defn}

$$
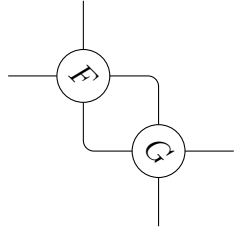
\begin{tikzpicture}

\draw[rounded corners](0,0)--(0,1)--(-1,1);
\draw[rounded corners](0,0)--(-1,0)--(-1,1);
\draw(1,0)--(0,0)--(0,-1);
\draw(-2,1)--(-1,1)--(-1,2);

\filldraw[fill=white](0,0) circle(10pt);
\node at (0,0)[rotate=-45]{$G$};

\filldraw[fill=white](-1,1) circle(10pt);
\node at (-1,1)[rotate=-45]{$F$};

\end{tikzpicture}
$$
\begin{defn}
\emph{Concatenation} is an operation which takes two reflected tangles and connects the bottom and right ends of the first tangle to the top and left ends of the second tangle, so that the concatenation $F,G$ is
$$
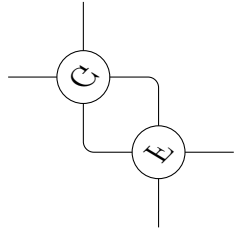
\begin{tikzpicture}

\draw[rounded corners](0,0)--(0,1)--(-1,1);
\draw[rounded corners](0,0)--(-1,0)--(-1,1);
\draw(1,0)--(0,0)--(0,-1);
\draw(-2,1)--(-1,1)--(-1,2);

\filldraw[fill=white](0,0) circle(10pt);
\node at (0,0)[rotate=-135,xscale=-1,yscale=1]{$F$};

\filldraw[fill=white](-1,1) circle(10pt);
\node at (-1,1)[rotate=-135,xscale=-1,yscale=1]{$G$};

\end{tikzpicture}
$$
\end{defn}

Note that concatenation is the result of multiplying each individual tangle by zero, and adding together the result.  Thus, it is a kind of algebraic tangle:
\begin{defn}
Any tangle that can be formed by iterated addition and multiplication of $n$-crossing-tangles is an \emph{algebraic tangle}.
\end{defn}

\begin{defn}

The \emph{closure} of a 4-tangle is the knot that results from
connecting the north and east strands, and the south and west strands,
outside of the tangle.

\end{defn}

$$
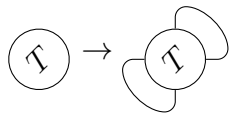
\begin{tikzpicture}[baseline]
\node (T)  at (0,0)  [shape=circle, draw, rotate=45]{$T$};
\end{tikzpicture}
\rightarrow
\begin{tikzpicture}[baseline]
\node (T)  at (0,0)  [shape=circle, draw, rotate=45]{$T$};
\draw[rounded corners] (T.45)--(90:7mm) arc (90:0:7mm) -- (T.-45);
\draw[rounded corners] (T.135)--(180:7mm) arc (180:270:7mm) -- (T.-135);
\end{tikzpicture}
$$

\begin{defn}

A \emph{knot diagram} (\emph{tangle diagram}) is a projection of a knot (tangle) into two dimensions, arranged so that the only intersections are between two non-parallel strands.  At these intersections, the overstrand is drawn continuously, and the understrand is drawn with a break.

\end{defn}

\begin{defn}[Conway Notation]
Recall John Conway's notation for knots consisting of rational tangles and algebraic tangles.  Where multiplications occur, the number of twists in each $n$-crossing-tangle are written side-by-side.  Parentheses are necessary as multiplication is not associative.  Where concatenations occur, a comma separates the tangles being concatenated.  The entire thing is surrounded by square brackets to show the closure of the tangle into a knot.  
\end{defn}

This is the tangle with Conway notation  $3 \, 2,4$.   :
\tikzset{every picture/.style={line width=0.75pt}} 

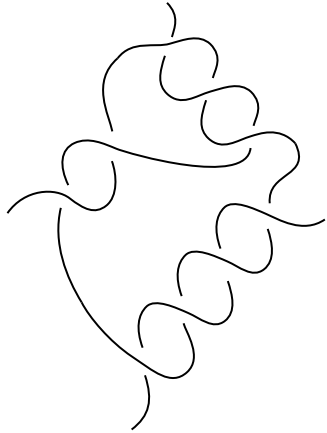
$$
\tikzset{every picture/.style={line width=0.75pt}} 
\begin{tikzpicture}[x=0.75pt,y=0.75pt,yscale=-1,xscale=1]

\draw    (220.58,126.13) .. controls (218.81,117.71) and (209.62,101.11) .. (223.27,89.21) ;
\draw    (198.37,153.25) .. controls (193.71,143.54) and (194.35,135.26) .. (201.51,132.01) .. controls (208.66,128.75) and (216.39,132.32) .. (222.55,134.86) .. controls (228.72,137.39) and (288.3,154.42) .. (290.29,134.52) ;
\draw    (250.63,78.52) .. controls (251.84,74.04) and (254.79,67.96) .. (248.16,61.24) ;
\draw    (312.84,132.16) .. controls (304.68,123.3) and (294.35,127.47) .. (290.82,128.46) .. controls (283.32,130.66) and (276.2,136.24) .. (269.07,129.71) .. controls (264.02,124.97) and (265,119.15) .. (267.88,110.4) ;
\draw    (271.25,99.19) .. controls (272.46,94.72) and (276.59,88.73) .. (269.96,82.02) .. controls (264.27,76.83) and (256.57,79.93) .. (249.51,81.9) .. controls (242.01,84.1) and (230.66,79.7) .. (223.27,89.21) ;
\draw    (291.79,123.32) .. controls (293.01,118.84) and (297.13,112.85) .. (290.51,106.14) .. controls (284.81,100.96) and (277.11,104.05) .. (270.05,106.02) .. controls (262.55,108.22) and (255.43,113.8) .. (248.3,107.27) .. controls (243.25,102.54) and (244.24,96.71) .. (247.11,87.97) ;
\draw    (167.68,167.42) .. controls (176.36,155.03) and (191.9,154.37) .. (199.05,159.9) .. controls (206.2,165.42) and (211.76,168.25) .. (217.59,163.56) .. controls (223.41,158.87) and (222.88,148.11) .. (220.36,141.14) ;
\draw    (205.1,212.42) .. controls (197.5,200) and (190.4,182.6) .. (194.6,164.73) ;
\draw    (255.56,209.27) .. controls (252.13,199.77) and (253.15,192.69) .. (257.57,188.47) .. controls (261.99,184.25) and (272.44,189.1) .. (278.55,191.76) .. controls (284.66,194.43) and (291.08,200.38) .. (297,195.81) .. controls (302.92,191.24) and (301.31,182.37) .. (298.94,175.35) ;
\draw    (235.85,235.3) .. controls (232.42,225.8) and (233.43,218.73) .. (237.85,214.5) .. controls (242.28,210.28) and (252.73,215.13) .. (258.84,217.8) .. controls (264.95,220.46) and (271.37,226.41) .. (277.29,221.85) .. controls (283.21,217.28) and (281.19,208.69) .. (278.82,201.67) ;
\draw    (205.1,212.42) .. controls (208.96,220.04) and (219.39,232.34) .. (230.31,239.67) .. controls (241.24,246.99) and (250.07,255.11) .. (258.25,247.69) .. controls (266.43,240.27) and (257.76,227.54) .. (256.5,222.99) ;
\draw    (230.31,276.61) .. controls (239.36,269.72) and (241.15,261.86) .. (237.05,249.25) ;
\draw    (275.13,185.29) .. controls (271.7,175.79) and (272.71,168.71) .. (277.13,164.49) .. controls (281.55,160.26) and (292,165.12) .. (298.11,167.78) .. controls (304.22,170.44) and (316.76,178.13) .. (327.76,170.59) ;
\draw    (299.96,162.42) .. controls (298.57,146.2) and (321.35,149.49) .. (312.84,132.16) ;

\end{tikzpicture}
$$

\subsection{Producing a knot from a Thompson tree diagram}

Here is Jones' method for creating knots from Thompson tree diagrams:

\begin{defn}\label{psi}
We define a function $\psi$ that maps Thompson diagrams to knot diagrams.  It is defined skein-theoretically:  
$\psi$ replaces branches with crossings by having an overstrand connect the two branches, and an understrand beginning along the trunk and continuing vertically through the crossing:
$$
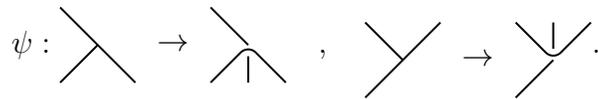
\psi:
\begin{tikzpicture}[baseline=.4cm, scale=0.5]
\draw(0,0)--(1,1);
\draw(0,2)--(2,0);
\node at (3,1){$\rightarrow$};
\draw[rounded corners](4,0)--(5,1)--(6,0);
\draw(4,2)--(5,1);
\draw(5,0.7)--(5,0);
\end{tikzpicture}
\quad , \quad
\begin{tikzpicture}[xscale=1,yscale=-1,scale=0.5, baseline=-.4cm]
\draw(0,0)--(1,1);
\draw(0,2)--(2,0);
\node at (3,1){$\rightarrow$};
\draw[rounded corners](4,0)--(5,1)--(6,0);
\draw(4,2)--(5,1);
\draw(5,0.7)--(5,0);
\end{tikzpicture}.
$$
The top and bottom corners of the diagram also become crossings, by making the two sides into an overstrand and introducing a vertical understrand:
$$
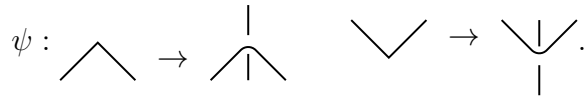
\psi: \begin{tikzpicture}[scale=0.5, baseline=.4cm]
\draw(0,0)--(1,1)--(2,0);
\node at (3,0.5){$\rightarrow$};
\draw[rounded corners](4,0)--(5,1)--(6,0);
\draw(5,0.7)--(5,0);
\draw(5,1.2)--(5,2);
\end{tikzpicture}
\quad \quad
\begin{tikzpicture}[xscale=1,yscale=-1,scale=0.5, baseline=-.4cm]
\draw(0,0)--(1,1)--(2,0);
\node at (3,0.5){$\rightarrow$};
\draw[rounded corners](4,0)--(5,1)--(6,0);
\draw(5,0.7)--(5,0);
\draw(5,1.2)--(5,2);
\end{tikzpicture}. 
$$
Then, each connected region of the modified diagram (including the exterior) has two newly created loose ends; connect these without introducing any further crossings \cite{MR3986040}.
\end{defn}

Here is an example:

$$
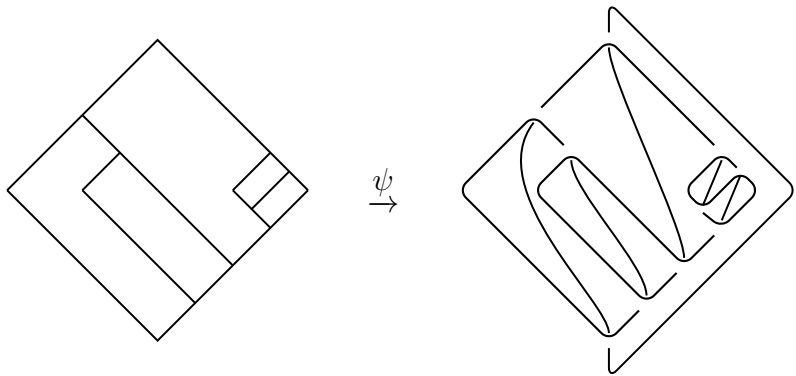
\begin{tikzpicture}

\draw(0,0)--(2,2)--(4,0);
\draw(1,1)--(2,0);
\draw(1.5,0.5)--(1,0);
\draw(3.5,0.5)--(3,0);
\draw(3.75,0.25)--(3.25,-0.25);

\begin{scope}[xscale=1,yscale=-1]
\draw(0,0)--(2,2)--(4,0);
\draw(2.5,1.5)--(1,0);
\draw(3,1)--(2,0);
\draw(3.5,0.5)--(3,0);

\end{scope}
\node at (5,0){$\underrightarrow{\psi}$};
\begin{scope}[xshift=6cm]
\draw[thick,rounded corners](2.4,-1.6)--(2,-2)--(0,0)--(1,1)--(1.4,0.6);
\draw[thick,rounded corners](1.1,1.1)--(2,2)--(3.4,0.6);
\draw[thick,rounded corners](3.25,-0.3)--(3.5,-0.5)--(4,0)--(3.75,0.25)--(3.25,-0.25)--(3,0)--(3.5,0.5)--(3.7,0.3);

\draw[thick,rounded corners](2.9,-1.1)--(2.5,-1.5)--(1,0)--(1.5,0.5)--(3,-1)--(3.4,-0.6);

\draw(1,0.9) .. controls (0.3,0) and (2,-1.5) .. (2,-1.9);
\draw(1.5,0.4) .. controls (1.5,0) and (2.5,-1)..(2.5,-1.4);
\draw(3.5,0.4)--(3.25,-0.2);
\draw(3.5,-0.4)--(3.75,0.2);
\draw(2,1.9) .. controls (2,1.5) and (3,-0.5)..(3,-0.9);
\draw[rounded corners](2,2.1)--(2,2.5)--(4.5,0)--(2,-2.5)--(2,-2.1);

\end{scope}

\end{tikzpicture}
$$

Note that $\psi$ can produce a knot or a link, and that $\psi$ turns a snipped Thompson diagram into a tangle.

\begin{thm}\label{thm:chairs}
Consider a chair, as shown below, which is part of a Thompson diagram.  The chair appears exactly as shown, with no additional edges or connections in the white region.  Outside of the white region, the rest of the diagram is unknown.  
After applying $\psi$ and performing ambient isotopies, this diagram becomes a tangle consisting of a 1-crossing-tangle above a line.

$$
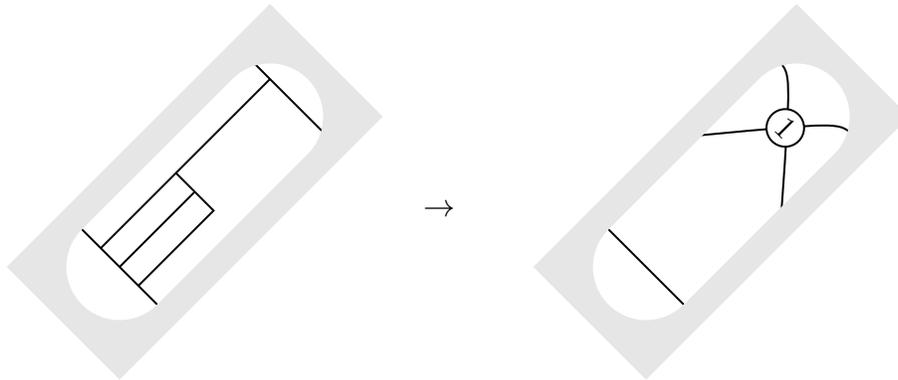
\begin{tikzpicture}[scale=0.5,baseline=-10pt]
	\draw (0.5,-0.5)--(2.5,-2.5);
	\draw (5,4)--(7,2);
	\draw (5.5,3.5)--(1,-1);
	\draw (3,1)--(4,0)--(2,-2);
	\draw ( 3.5,0.5)--(1.5,-1.5);

			\fill [fill=black!10!white] (0.5,-0.5)--(4.5,3.5) arc (135:-45:1.414) -- (2.5,-2.5) -- (1.5,-4.5) -- (8.5,2.5) -- (5.5,5.5) -- (-1.5,-1.5) -- (1.5,-4.5) -- (2.5,-2.5) arc (-45:-225:1.414);

	\node at (10,0) {$\rightarrow$};

\begin{scope}[xshift=14cm]
		\draw (0.5,-0.5)--(2.5,-2.5);

		\draw (5,4) .. controls ++(-45:0.5cm) .. (5.1,0);
		\draw (7,2) .. controls ++(135:0.5cm) .. (2.9,2);
		\filldraw[fill=white] (5.2,2.2) circle (0.5cm) node [rotate=-45]{$1$};

		\fill [fill=black!10!white] (0.5,-0.5)--(4.5,3.5) arc (135:-45:1.414) -- (2.5,-2.5) -- (1.5,-4.5) -- (8.5,2.5) -- (5.5,5.5) -- (-1.5,-1.5) -- (1.5,-4.5) -- (2.5,-2.5) arc (-45:-225:1.414);
\end{scope}
\end{tikzpicture}
$$

\end{thm}

\begin{proof}

Beginning with the chair in the Thompson diagram, first we apply $\psi$, then ambient isotopies:

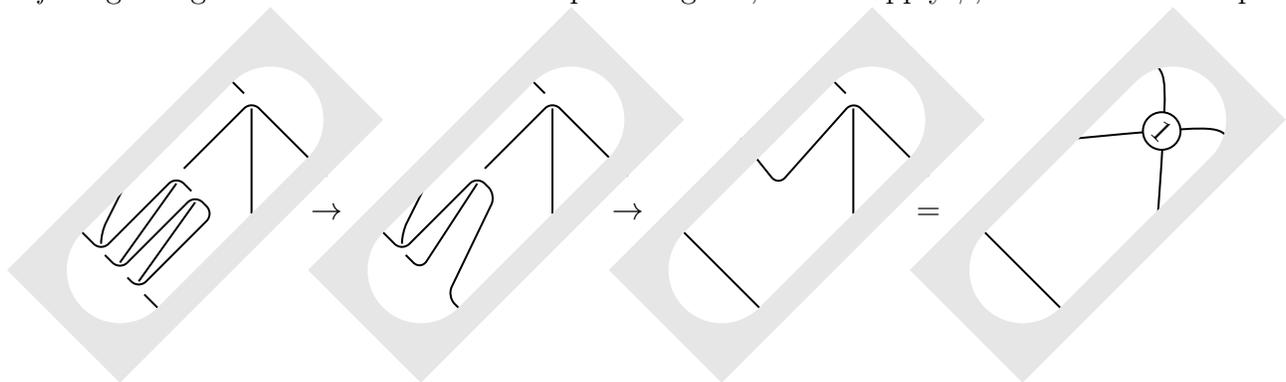
\begin{tikzpicture}[scale=.5,baseline]
	\draw[rounded corners](1.1,-1.1)--(1.5,-1.5)--(3.5,0.5)--(4,0)--(2,-2)--(1.7,-1.7);
	\draw[rounded corners](0.5,-0.5)--(1,-1)--(3,1)--(3.4,0.6);
	\draw(4.5,3.5)--(4.8,3.2);
	\draw[rounded corners](7,1)--(5,3)--(3.2,1.2);
	\draw(3,0.8)--(1.5,-1.3);
	\draw(3.5,0.3)--(2,-1.8);
	\draw[rounded corners](1,-0.8)--(1,-0.5)--(2,1.5);
	\draw(5,2.8)--(5,0);
	\draw(2.15,-2.15)--(2.5,-2.5);

	\fill [fill=black!10!white] (0.5,-0.5)--(4.5,3.5) arc (135:-45:1.414) -- (2.5,-2.5) -- (1.5,-4.5) -- (8.5,2.5) -- (5.5,5.5) -- (-1.5,-1.5) -- (1.5,-4.5) -- (2.5,-2.5) arc (-45:-225:1.414);

	\node at (7,0) {$\rightarrow$};

\begin{scope}[xshift=8cm]
	\draw [color=white,fill=white] (5,3) circle (20pt);
	\draw[color=white,fill=white] (1.45,-1.45) circle (40pt);
	\draw[rounded corners](0.5,-0.5)--(1,-1)--(3,1)--(3.5,0.5)--(2.2,-2.2)--(2.5,-2.5);
	\draw[rounded corners](1.1,-1.1)--(1.5,-1.5)--(3,0.8);
	\draw[rounded corners](7,1)--(5,3)--(3.2,1.2);
	\draw[rounded corners](1,-0.8)--(1,-0.5)--(2,1.5);
	\draw(5,2.8)--(5,0);
	\draw (4.5,3.5)--(4.8,3.2);
	\fill [fill=black!10!white] (0.5,-0.5)--(4.5,3.5) arc (135:-45:1.414) -- (2.5,-2.5) -- (1.5,-4.5) -- (8.5,2.5) -- (5.5,5.5) -- (-1.5,-1.5) -- (1.5,-4.5) -- (2.5,-2.5) arc (-45:-225:1.414);

\end{scope}
	\node at (15,0) {$\rightarrow$};

	\begin{scope}[xshift=16cm]
		\draw[rounded corners](7,1)--(5,3)--(3,0.75)--(2,2);
		\draw (4.5,3.5)--(4.8,3.2);
		\draw(5,2.8)--(5,0);
		\draw(0.5,-0.5)--(2.5,-2.5);
	\fill [fill=black!10!white] (0.5,-0.5)--(4.5,3.5) arc (135:-45:1.414) -- (2.5,-2.5) -- (1.5,-4.5) -- (8.5,2.5) -- (5.5,5.5) -- (-1.5,-1.5) -- (1.5,-4.5) -- (2.5,-2.5) arc (-45:-225:1.414);
	\end{scope}
	
	\node at (23,0) {$=$};

	\begin{scope}[xshift=24cm]
		\draw (0.5,-0.5)--(2.5,-2.5);

		\draw (5,4) .. controls ++(-45:0.5cm) .. (5.1,0);
		\draw (7,2) .. controls ++(135:0.5cm) .. (2.9,2);
		\filldraw[fill=white] (5.2,2.2) circle (0.5cm) node [rotate=-45]{$1$};

	\fill [fill=black!10!white] (0.5,-0.5)--(4.5,3.5) arc (135:-45:1.414) -- (2.5,-2.5) -- (1.5,-4.5) -- (8.5,2.5) -- (5.5,5.5) -- (-1.5,-1.5) -- (1.5,-4.5) -- (2.5,-2.5) arc (-45:-225:1.414);
	\end{scope}
\end{tikzpicture}
%
%

\end{proof}

\begin{thm}\label{nchairs}
 Applying $\psi$ to $n$ chairs in a row produces an $n$-crossing tangle above a line, as shown below.

$$
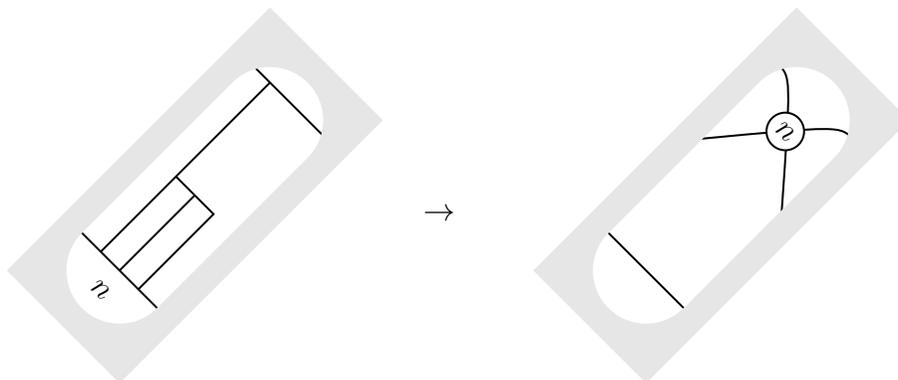
\begin{tikzpicture}[scale=.5,baseline]

	\draw (0.5,-0.5)--(2.5,-2.5);
	\draw (5,4)--(7,2);
	\draw (5.5,3.5)--(1,-1);
	\draw (3,1)--(4,0)--(2,-2);
	\draw ( 3.5,0.5)--(1.5,-1.5);

\node at (1,-2)[rotate=-45]{$n$};

	\fill [fill=black!10!white] (0.5,-0.5)--(4.5,3.5) arc (135:-45:1.414) -- (2.5,-2.5) -- (1.5,-4.5) -- (8.5,2.5) -- (5.5,5.5) -- (-1.5,-1.5) -- (1.5,-4.5) -- (2.5,-2.5) arc (-45:-225:1.414);

	\node at (10,0) {$\rightarrow$};

\begin{scope}[xshift=14cm]
		\draw (0.5,-0.5)--(2.5,-2.5);

		\draw (5,4) .. controls ++(-45:0.5cm) .. (5.1,0);
		\draw (7,2) .. controls ++(135:0.5cm) .. (2.9,2);
		\filldraw[fill=white] (5.2,2.2) circle (0.5cm) node [rotate=-45]{$n$};

		\fill [fill=black!10!white] (0.5,-0.5)--(4.5,3.5) arc (135:-45:1.414) -- (2.5,-2.5) -- (1.5,-4.5) -- (8.5,2.5) -- (5.5,5.5) -- (-1.5,-1.5) -- (1.5,-4.5) -- (2.5,-2.5) arc (-45:-225:1.414);
\end{scope}
\end{tikzpicture}
$$

\end{thm}

\begin{proof}

Each chair will become a $1$-crossing-tangle connected in the following manner:

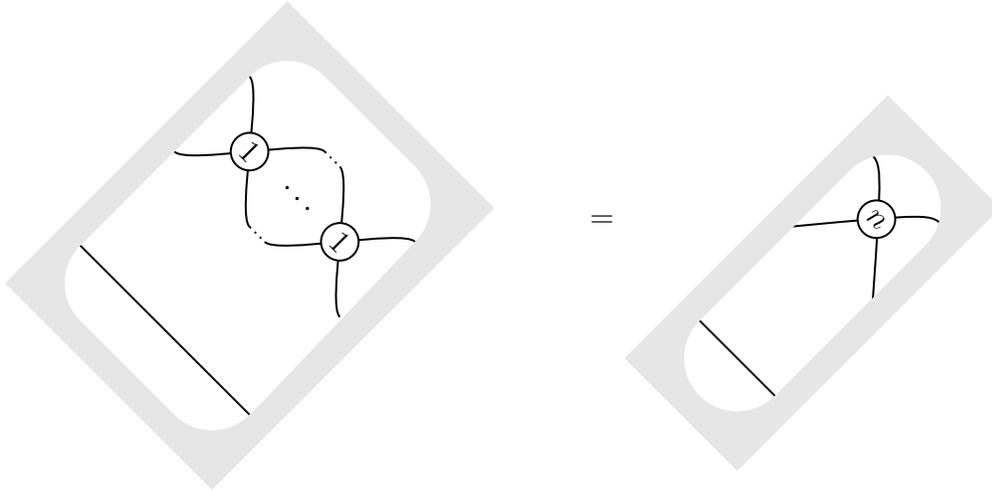
\begin{tikzpicture}[scale=.5, baseline]
		\draw (0.5,-0.5)--(5,-5);

		\draw (5,4) .. controls ++(-45:0.5cm)  and ++(135:0.5cm).. (5,0);
		\draw[dotted] (5,0) -- (5.4,-0.4);
		\draw (7,2) .. controls   ++(135:0.5cm) and ++(-45:0.5cm) .. (3,2);
		\draw[dotted] (7,2) -- (7.4,1.6);
		\filldraw[fill=white] (5,2) circle (0.5cm) node [rotate=-45]{$1$};
		
		\node[rotate=-45] at (6.3,0.7) {$\cdots$};
	\begin{scope}[shift={(2.4,-2.4)}]
]]		\draw (5,4) .. controls ++(-45:0.5cm)  and ++(135:0.5cm).. (5,0);
		\draw (7,2) .. controls   ++(135:0.5cm) and ++(-45:0.5cm) .. (3,2);
		\filldraw[fill=white] (5,2) circle (0.5cm) node [rotate=-45]{$1$};
	\end{scope}

		\fill [fill=black!10!white] (0.5,-0.5)--(5,4) arc (135:45:1.414) -- (9.4,1.6) arc (45:-45:1.414) -- (5.4,-4.6) -- (6,-5) -- (11.5,0.5) -- (6,6) -- (-1.5,-1.5) -- (4,-7) -- (6,-5) -- (5.4,-4.6) -- (5,-5) arc (-45:-135:1.414) -- (0.5,-2.5)  arc (-135:-225:1.414);
\end{tikzpicture}
\hspace{1cm} =
%
%
%
\begin{tikzpicture}[scale=0.5,baseline=1cm]
		\draw (0.5,-0.5)--(2.5,-2.5);

		\draw (5,4) .. controls ++(-45:0.5cm) .. (5.1,0);
		\draw (7,2) .. controls ++(135:0.5cm) .. (2.9,2);
		\filldraw[fill=white] (5.2,2.2) circle (0.5cm) node [rotate=-45]{$n$};

		\fill [fill=black!10!white] (0.5,-0.5)--(4.5,3.5) arc (135:-45:1.414) -- (2.5,-2.5) -- (1.5,-4.5) -- (8.5,2.5) -- (5.5,5.5) -- (-1.5,-1.5) -- (1.5,-4.5) -- (2.5,-2.5) arc (-45:-225:1.414);
\end{tikzpicture}

\end{proof}

\begin{cor}
$n$ reflected chairs gives a reflected $n$-crossing-tangle.
\end{cor}

\begin{defn}\label{psiprime}
We define a variation on $\psi$ which sends Thompson diagrams to knot diagrams (and snipped Thompson diagrams to tangle diagrams) which are a few steps more simplified.  Let $\psi'$ work by directly replacing $n$ parallel chairs with an $n$-crossing-tangle,
and corners which are not part of chairs are mapped according to the same rules as $\psi$.  
\end{defn}

Here is an example with a Thompson diagram and knot diagram:

$$
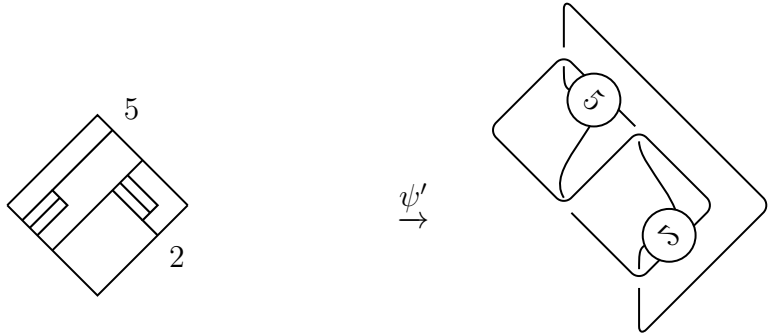
\begin{tikzpicture}[baseline]
\begin{scope}[scale=0.2] 
    \uchair[thick](-2,-2)(1:2);
	
    \chair[thick](-5,5)(4:2);

    \draw (0,0)--(-6,6)--(-12,0);
    \draw (0,0)--(-6,-6)--(-12,0);
    \draw (-3,3)--(-9,-3);

    \node at (-2,-2) [below right] {$2$};
    \node at (-5,5) [above right] {$5$};
\end{scope}

\node at (3,0){$\underrightarrow{\psi'}$};

\begin{scope}[shift={(5,0)}]
\draw[rounded corners](0.95,1.05)--(0.5,1.5)--(0,2)--(-1,1)--(0,0)--(0.5,0.5)--(1,1)--(2,0)--(1,-1)--(0.1,-0.1);
\draw (1.5,-0.5) .. controls (1,-0.5) .. (1,-0.85);
\draw(1.5,-0.5) .. controls (1.7,0) and (1,0.5) .. (1,0.85);
\draw(0.5,1.5) .. controls (0.5,1) and (-0.25,0.5) .. (0,0.1);
\draw(0.5,1.5) .. controls (0,1.5) .. (0,1.85);
\draw[rounded corners](0,2.1)--(0,2.75)--(2.75,0)--(1,-1.75)--(1,-1.1);
\filldraw[fill=white](0.4,1.4) circle(10pt);
\filldraw[fill=white](1.4,-0.4) circle (10pt);
\node at (0.4,1.4)[rotate=-45]{$5$};
\node at (1.4,-0.4)[rotate=45,xscale=1,yscale=-1]{$2$};
\end{scope}
\end{tikzpicture}
$$

\begin{thm}
Let $\mathcal{T}$ represent the set of tangles, $\mathcal{KD}$ the set of knot diagrams, and $\mathcal{K}$ the set of knots.  $\phi$ is the usual interpretation of knot diagrams as representatives of knots.  The following diagram commutes:
$$\begin{tikzcd}
\mathcal{T} \arrow[r, "\psi'"] \arrow [d, "\psi"] & \mathcal{KD} \arrow[d, "\phi"] \\
\mathcal{KD} \arrow[r,"\phi"] & \mathcal{K}
\end{tikzcd}
$$

\end{thm}

\begin{proof}
As in the proof of Theorem \ref{thm:chairs}, we go from the bottom left $\mathcal{KD}$ to the top right $\mathcal{KD}$ by ambient isotopy.  As knots are defined up to ambient isotopy, the end result in both cases is the same.  
\end{proof}

\subsection{Jones' algorithm}

Vaughan Jones provided  an algorithm which, given a knot $K$, produces a Thompson diagram $T$ whose associated knot is again $K$:  $\psi(T)=K$ \cite{MR3589908}.   As told to us by Jones himself, this algorithm has a few unsatisfying aspects.  Essentially, this algorithm is very far from functorial (even though we don't know exactly what functorial should mean here), and not in any sense skein-theoretic or tangle-theoretic.

Jones's algorithm is based on the following observation:  If one begins with a Thompson diagram $T$ and produces a knot $\psi(T)$, then the diagram for $\psi(T)$ will have a checkerboard shading on its regions (as do all knot diagrams), in which every region intersects the horizontal midline exactly once.  

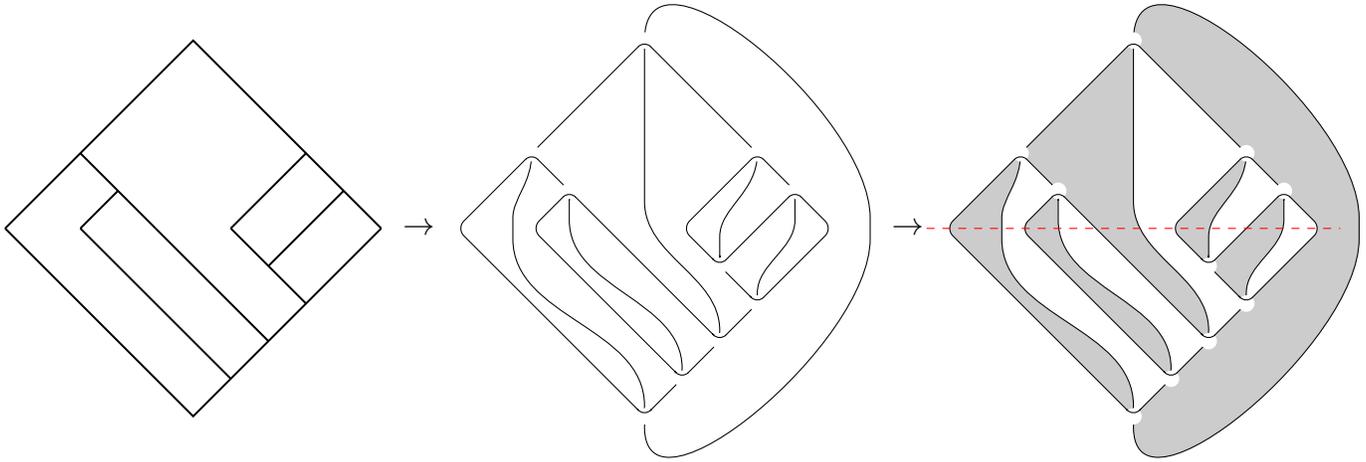
\hspace*{-1cm}
\begin{tikzpicture}

\draw(0,0)--(2.5,2.5)--(5,0);
\draw(1,1)--(2,0);
\draw(1.5,0.5)--(1,0);
\draw(4,1)--(3,0);
\draw(4.5,0.5)--(3.5,-0.5);

\begin{scope}[xscale=1,yscale=-1]
\draw(0,0)--(2.5,2.5)--(5,0);
\draw(3,2)--(1,0);
\draw (3.5,1.5)--(2,0);
\draw(4,1)--(3,0);
\end{scope}

\node at (5.5,0){$\rightarrow$};

\begin{scope}[xshift=6cm, yshift=-2.5cm, scale=.5]
  \node[cross] (d0) at (2, 7) {};
  \node[cross] (d1) at (3, 6) {};
  \node[cross] (d2) at (5,10) {};
  \node[cross] (d3) at (8,7) {};
  \node[cross] (d4) at (9,6) {};
  \node[cross] (d5) at (8,3) {};
  \node[cross] (d6) at (7,2) {};
  \node[cross] (d7) at (6,1) {};
  \node[cross] (d8) at (5,0) {};
  \node[cross] (d9) at (7,4) {};
  
\begin{pgfonlayer}{foreground}  
  \draw (d0.-45) to[out=135, in = 45] (d0.-135);
  \draw (d1.-45) to[out=135, in = 45] (d1.-135);
  \draw (d2.-45) to[out=135, in = 45] (d2.-135);
  \draw (d3.-45) to[out=135, in = 45] (d3.-135);
  \draw (d4.-45) to[out=135, in = 45] (d4.-135);
  \draw (d5.45) to[out=-135, in = -45] (d5.135);
  \draw (d6.45) to[out=-135, in = -45] (d6.135);
  \draw (d7.45) to[out=-135, in = -45] (d7.135);
  \draw (d8.45) to[out=-135, in = -45] (d8.135);
  \draw (d9.45) to[out=-135, in = -45] (d9.135);
\end{pgfonlayer}

  \draw (d7.center) to (d6.center);
  \draw (d2.center) to (d0.center);
  \draw (d0.center) to (d1.center);
  \draw[rounded corners] (d8.center) to (0,5) to (d0.center);
  \draw[rounded corners] (d1.center) to (2,5) to (d7.center);
  \draw (d8.center) to (d7.center);
  \draw (d1.center) to (d6.center);
  \draw (d3.center) to (d4.center);
  \draw[rounded corners] (d4.center) to (10,5) to (d5.center);
  \draw (d6.center) to (d5.center);
  \draw (d5.center) to (d9.center);
  \draw[rounded corners] (d9.center) to (6,5) to (d3.center);
  \draw (d3.center) to (d2.center);
  \draw (d4.center) to (d9.center);
  
  \draw[rounded corners] (d0.center) to[out=-90,in=90] (1.5,5) to[out=-90, in=90] (d8.center) to (d8.135);
  \draw[rounded corners] (d1.center) to[out=-90,in=90] (3,5.5) to[out=-90,in=90] (d7.center) to (d7.135);
  \draw[rounded corners] (d2.center) to[out=-90,in=90] (5,6) to[out=-90,in=90] (d6.center) to (d6.135);
  \draw[rounded corners] (d3.center) to[out=-90,in=90] (7,4.5) to[out=-90,in=90] (d9.center) to (d9.135);
  \draw[rounded corners] (d4.center) to[out=-90,in=90] (9,5.5) to[out=-90,in=90] (d5.center) to (d5.135);

  \draw[rounded corners] (d2.center) to[out=90, in=90] (11,5) to[out=-90, in=-90] (d8.center) to (d8.45);
\end{scope}

\node at (12,0){$\rightarrow$};

\begin{scope}[xshift=12.5cm, yshift=-2.5cm, scale=.5]
  \node[cross] (d0) at (2, 7) {};
  \node[cross] (d1) at (3, 6) {};
  \node[cross] (d2) at (5,10) {};
  \node[cross] (d3) at (8,7) {};
  \node[cross] (d4) at (9,6) {};
  \node[cross] (d5) at (8,3) {};
  \node[cross] (d6) at (7,2) {};
  \node[cross] (d7) at (6,1) {};
  \node[cross] (d8) at (5,0) {};
  \node[cross] (d9) at (7,4) {};
  
\begin{pgfonlayer}{foreground}  
  \draw (d0.-45) to[out=135, in = 45] (d0.-135);
  \draw (d1.-45) to[out=135, in = 45] (d1.-135);
  \draw (d2.-45) to[out=135, in = 45] (d2.-135);
  \draw (d3.-45) to[out=135, in = 45] (d3.-135);
  \draw (d4.-45) to[out=135, in = 45] (d4.-135);
  \draw (d5.45) to[out=-135, in = -45] (d5.135);
  \draw (d6.45) to[out=-135, in = -45] (d6.135);
  \draw (d7.45) to[out=-135, in = -45] (d7.135);
  \draw (d8.45) to[out=-135, in = -45] (d8.135);
  \draw (d9.45) to[out=-135, in = -45] (d9.135);
\end{pgfonlayer}

  \draw (d7.center) to (d6.center);
  \draw (d2.center) to (d0.center);
  \draw (d0.center) to (d1.center);
  \draw[rounded corners] (d8.center) to (0,5) to (d0.center);
  \draw[rounded corners] (d1.center) to (2,5) to (d7.center);
  \draw (d8.center) to (d7.center);
  \draw (d1.center) to (d6.center);
  \draw (d3.center) to (d4.center);
  \draw[rounded corners] (d4.center) to (10,5) to (d5.center);
  \draw (d6.center) to (d5.center);
  \draw (d5.center) to (d9.center);
  \draw[rounded corners] (d9.center) to (6,5) to (d3.center);
  \draw (d3.center) to (d2.center);
  \draw (d4.center) to (d9.center);
  
  \draw[rounded corners] (d0.center) to[out=-90,in=90] (1.5,5) to[out=-90, in=90] (d8.center) to (d8.135);
  \draw[rounded corners] (d1.center) to[out=-90,in=90] (3,5.5) to[out=-90,in=90] (d7.center) to (d7.135);
  \draw[rounded corners] (d2.center) to[out=-90,in=90] (5,6) to[out=-90,in=90] (d6.center) to (d6.135);
  \draw[rounded corners] (d3.center) to[out=-90,in=90] (7,4.5) to[out=-90,in=90] (d9.center) to (d9.135);
  \draw[rounded corners] (d4.center) to[out=-90,in=90] (9,5.5) to[out=-90,in=90] (d5.center) to (d5.135);
  
  \begin{pgfonlayer}{background}
  \fill[ black!20!white, rounded corners] (d0.center) to[out=-90,in=90] (1.5,5) to[out=-90, in=90] (d8.center) to (0,5) to (d0.center);
  \fill[black!20!white, rounded corners]  (d1.center) to[out=-90,in=90] (3,5.5) to[out=-90, in=90] (d7.center) to (2,5) to (d1.center);
  \fill[black!20!white, rounded corners] (d2.center) to[out=-90,in=90] (5,6) to[out=-90, in=90] (d6.center) to (2,7) to (d2.center);
  \fill[black!20!white, rounded corners] (d3.center) to[out=-90,in=90] (7,4.5) to[out=-90, in=90] (d9.center) to (6,5) to (d3.center);
  \fill[black!20!white, rounded corners] (d4.center) to[out=-90,in=90] (9,5.5) to[out=-90, in=90] (d5.center) to (7,4) to (d4.center);
  \fill[black!20!white, rounded corners] (d2.center) to[out=90, in=90] (11,5) to[out=-90, in=-90] (d8.center) to (10,5) to (d2.center);
\end{pgfonlayer}

\draw[rounded corners] (d2.center) to[out=90, in=90] (11,5) to[out=-90, in=-90] (d8.center) to (d8.45);

\draw[red, dashed] (-0.5,5)--(10.5,5);

\end{scope}

\end{tikzpicture}

We observe some further patterns:  There are two types of shaded crossing.  A positive crossing has the form
\begin{tikzpicture}[baseline=-0.1cm]
   \fill[ black!20!white] (-45:1cm) to (0,0) to (-135:1cm);
   \fill[ black!20!white] (45:1cm) to (0,0) to (135:1cm);
   \draw (-45:1cm) to (-45:2mm);
   \draw (135:2mm) -- (135:1cm);
   \draw (45:1cm)--(-135:1cm);
\end{tikzpicture}
and a negative crossing has the form
\begin{tikzpicture}[baseline=-0.1cm]
   \fill[ black!20!white] (-45:1cm) to (0,0) to (-135:1cm);
   \fill[ black!20!white] (45:1cm) to (0,0) to (135:1cm);
   \draw (45:1cm) to (45:2mm);
   \draw (-135:2mm) -- (-135:1cm);
   \draw (-45:1cm)--(135:1cm);
\end{tikzpicture}.
Now, looking at the crossing types in our shaded Thompson knot, one notices only positive crossings above the midline and negative crossings below.  

Next, we focus on the connection data between the unshaded regions.  Place a vertex in each unshaded region (and one representing the exterior region, to the left of the diagram), and if two unshaded regions meet at a vertex, connect them by an edge.  Label the edge by the type (positive or negative) of crossing.  Notice that each vertex except the  leftmost has exactly one positive edge coming in from a further left vertex, and exactly one negative edge coming from further left.  (This is because all unshaded regions have a `medial string' forming their left boundary.)  A consequence of this is that our planar graph is in fact the union of two trees -- one made of the positive edges, one of the negative.  

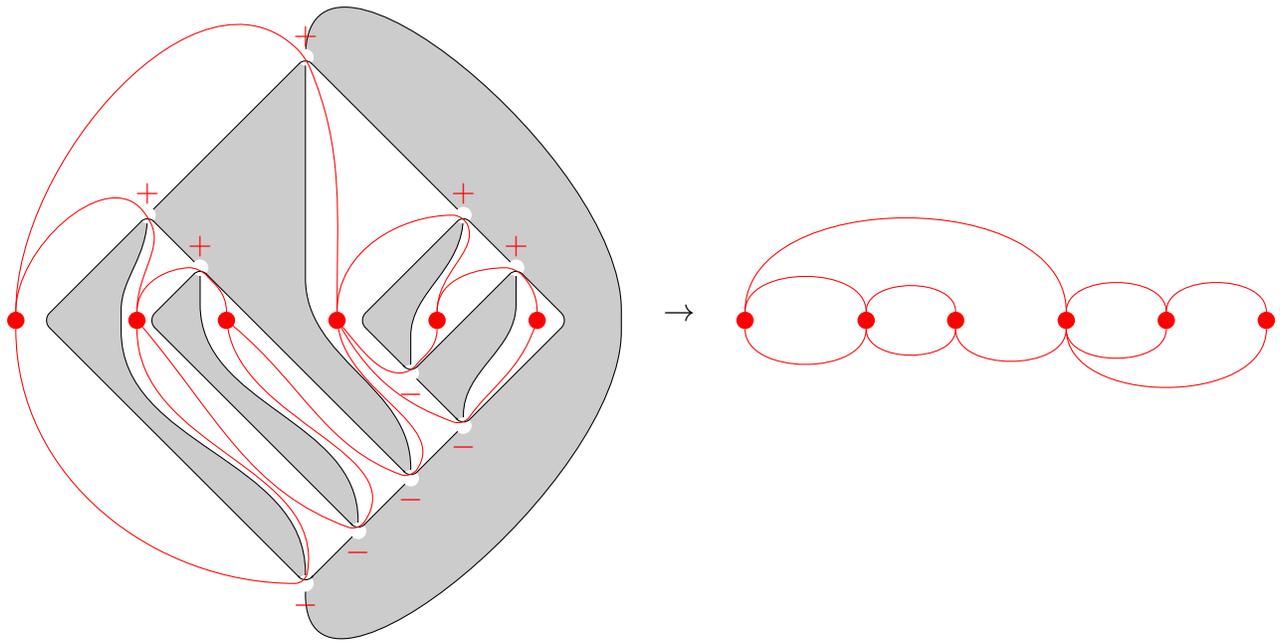
\begin{tikzpicture}[baseline=3.5cm, scale=.7]
  \node[cross] (d0) at (2, 7) {};
  \node[cross] (d1) at (3, 6) {};
  \node[cross] (d2) at (5,10) {};
  \node[cross] (d3) at (8,7) {};
  \node[cross] (d4) at (9,6) {};
  \node[cross] (d5) at (8,3) {};
  \node[cross] (d6) at (7,2) {};
  \node[cross] (d7) at (6,1) {};
  \node[cross] (d8) at (5,0) {};
  \node[cross] (d9) at (7,4) {};
  
\begin{pgfonlayer}{foreground}  
  \draw (d0.-45) to[out=135, in = 45] (d0.-135);
  \draw (d1.-45) to[out=135, in = 45] (d1.-135);
  \draw (d2.-45) to[out=135, in = 45] (d2.-135);
  \draw (d3.-45) to[out=135, in = 45] (d3.-135);
  \draw (d4.-45) to[out=135, in = 45] (d4.-135);
  \draw (d5.45) to[out=-135, in = -45] (d5.135);
  \draw (d6.45) to[out=-135, in = -45] (d6.135);
  \draw (d7.45) to[out=-135, in = -45] (d7.135);
  \draw (d8.45) to[out=-135, in = -45] (d8.135);
  \draw (d9.45) to[out=-135, in = -45] (d9.135);
\end{pgfonlayer}

  \draw (d7.center) to (d6.center);
  \draw (d2.center) to (d0.center);
  \draw (d0.center) to (d1.center);
  \draw[rounded corners] (d8.center) to (0,5) to (d0.center);
  \draw[rounded corners] (d1.center) to (2,5) to (d7.center);
  \draw (d8.center) to (d7.center);
  \draw (d1.center) to (d6.center);
  \draw (d3.center) to (d4.center);
  \draw[rounded corners] (d4.center) to (10,5) to (d5.center);
  \draw (d6.center) to (d5.center);
  \draw (d5.center) to (d9.center);
  \draw[rounded corners] (d9.center) to (6,5) to (d3.center);
  \draw (d3.center) to (d2.center);
  \draw (d4.center) to (d9.center);
  
  \draw[rounded corners] (d0.center) to[out=-90,in=90] (1.5,5) to[out=-90, in=90] (d8.center) to (d8.135);
  \draw[rounded corners] (d1.center) to[out=-90,in=90] (3,5.5) to[out=-90,in=90] (d7.center) to (d7.135);
  \draw[rounded corners] (d2.center) to[out=-90,in=90] (5,6) to[out=-90,in=90] (d6.center) to (d6.135);
  \draw[rounded corners] (d3.center) to[out=-90,in=90] (7,4.5) to[out=-90,in=90] (d9.center) to (d9.135);
  \draw[rounded corners] (d4.center) to[out=-90,in=90] (9,5.5) to[out=-90,in=90] (d5.center) to (d5.135);
  
  \node[dot] (a0) at (-0.5, 5) {};  
  \node[dot] (a1) at (1.8, 5) {};
  \node[dot] (a2) at (3.5,5) {};
  \node[dot] (a3) at (5.6,5) {};
  \node[dot] (a4) at (7.5,5) {};
  \node[dot] (a5) at (9.4,5) {};
  
\begin{pgfonlayer}{foreground}    
  \draw[red, rounded corners] (a0) to[out=90, in=135] (d0.center) node[above] {$+$} to[out=-67, in=90] (a1);
  \draw[red, rounded corners] (a0) to[out=90, in=135] (d2.center) node[above] {$+$} to[out=-67, in=90] (a3);
  \draw[red, rounded corners] (a1) to[out=90, in=180] (d1.center) node[above] {$+$} to[out=-50, in=90] (a2);
  \draw[red, rounded corners] (a3) to[out=90, in=180] (d3.center) node[above] {$+$} to[out=-67, in=90] (a4);
  \draw[red, rounded corners] (a4) to[out=90, in=180] (d4.center) node[above] {$+$} to[out=-50, in=90] (a5);
  
    \draw[red, rounded corners] (a3) to[out=-60, in=160] (d5.center) node[below] {$-$} to[out=55, in=-110] (a5);
    \draw[red, rounded corners] (a2) to[out=-40, in=160] (d6.center) node[below] {$-$} to[out=55, in=-80] (a3);
    \draw[red, rounded corners] (a1) to[out=-50, in=160] (d7.center) node[below] {$-$} to[out=55, in=-80] (a2);
    \draw[red, rounded corners] (a0) to[out=-90, in=180] (d8.center) node[below] {$-$} to[out=80, in=-90] (a1);
    \draw[red, rounded corners] (a3) to[out=-50, in=180] (d9.center) node[below] {$-$} to[out=55, in=-90] (a4);

\end{pgfonlayer}
  
  \begin{pgfonlayer}{background}
  \fill[black!20!white, rounded corners] (d0.center) to[out=-90,in=90] (1.5,5) to[out=-90, in=90] (d8.center) to (0,5) to (d0.center);
  \fill[black!20!white, rounded corners]  (d1.center) to[out=-90,in=90] (3,5.5) to[out=-90, in=90] (d7.center) to (2,5) to (d1.center);
  \fill[black!20!white, rounded corners] (d2.center) to[out=-90,in=90] (5,6) to[out=-90, in=90] (d6.center) to (2,7) to (d2.center);
  \fill[black!20!white, rounded corners] (d3.center) to[out=-90,in=90] (7,4.5) to[out=-90, in=90] (d9.center) to (6,5) to (d3.center);
  \fill[black!20!white, rounded corners] (d4.center) to[out=-90,in=90] (9,5.5) to[out=-90, in=90] (d5.center) to (7,4) to (d4.center);
  \fill[black!20!white, rounded corners] (d2.center) to[out=90, in=90] (11,5) to[out=-90, in=-90] (d8.center) to (10,5) to (d2.center);
\end{pgfonlayer}

\draw[rounded corners] (d2.center) to[out=90, in=90] (11,5) to[out=-90, in=-90] (d8.center) to (d8.45);
\end{tikzpicture}
$\quad \rightarrow \quad$
\begin{tikzpicture}[baseline=3.5cm, scale=.7]
  \node[dot] (a0) at (-0.5, 5) {};  
  \node[dot] (a1) at (1.8, 5) {};
  \node[dot] (a2) at (3.5,5) {};
  \node[dot] (a3) at (5.6,5) {};
  \node[dot] (a4) at (7.5,5) {};
  \node[dot] (a5) at (9.4,5) {};
  
\begin{pgfonlayer}{foreground}    
  \draw[red, rounded corners] (a0) to[out=90, in=90] (a1);
  \draw[red, rounded corners] (a0) to[out=90, in=90] (a3);
  \draw[red, rounded corners] (a1) to[out=90, in=90] (a2);
  \draw[red, rounded corners] (a3) to[out=90, in=90] (a4);
  \draw[red, rounded corners] (a4) to[out=90, in=90] (a5);
  
    \draw[red, rounded corners] (a3) to[out=-90, in=-90]  (a5);
    \draw[red, rounded corners] (a2) to[out=-90, in=-90] (a3);
    \draw[red, rounded corners] (a1) to[out=-90, in=-90] (a2);
    \draw[red, rounded corners] (a0) to[out=-90, in=-90] (a1);
    \draw[red, rounded corners] (a3) to[out=-90, in=-90]  (a4);

\end{pgfonlayer}
\end{tikzpicture}

If we wish to run this procedure backwards -- that is, start from a knot and render it as the image of a Thompson tree diagram -- we start by checkerboard shading the regions of the knot, putting a vertex in each unshaded region, and then isotoping the diagram so that the new vertices sit along a midline.

Our running example for constructing a Thompson diagram associated to a knot will be the following:

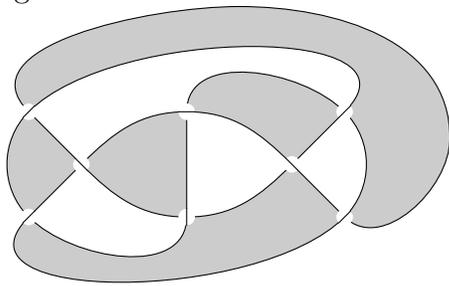
\begin{tikzpicture}[scale=0.7]
  \node[cross] (d0) at (0,0) {};				  \node[cross] (d1) at (0,2) {};
  \node[cross] (d2) at (1,1) {};				  \node[cross] (d3) at (3,0) {};
  \node[cross] (d4) at (3,2) {};				  \node[cross] (d5) at (5,1) {};
  \node[cross] (d6) at (6,0) {};				  \node[cross] (d7) at (6,2) {};

\begin{pgfonlayer}{foreground}  
  \draw (d0.45) -- (d0.-135);				  \draw (d1.45) -- (d1.-135);
  \draw (d2.-45) -- (d2.135);				  \draw (d3.90) -- (d3.-90);
  \draw (d4.0) -- (d4.180);					  \draw (d5.-45) -- (d5.135);
  \draw (d6.45) -- (d6.-135);				  \draw (d7.45) -- (d7.-135);
\end{pgfonlayer}

\draw (d0.center) to[out=135, in=-135] (d1.center);			\draw (d0.center) to[out= 45, in=-135] (d2.center);
\draw (d0.center) to[out=-45 , in=-90] (d3.center);			\draw (d0.center) to[out=-135 , in=-135] (d6.center);
\draw (d1.center) to[out=-45 , in=135] (d2.center);			\draw (d1.center) to[out=45 , in=45] (d7.center);
\draw (d1.center) to[out=135 , in=180] (4,4) to [out=0,in=90] (8,1.5) to [out=-90, in=-45](d6.center);
\draw (d2.center) to[out=-45 , in=180] (d3.center);			\draw (d2.center) to[out=45 , in=180] (d4.center);
\draw (d3.center) to[out=90 , in=-90] (d4.center);			\draw (d3.center) to[out=0 , in=-135] (d5.center);
\draw (d4.center) to[out=0 , in=135] (d5.center);				\draw (d4.center) to[out=90 , in=135] (d7.center);
\draw (d5.center) to[out=-45 , in=135] (d6.center);			\draw (d5.center) to[out=45 , in=-135] (d7.center);
\draw (d6.center) to[out=45 , in=-45] (d7.center);

\begin{pgfonlayer}{background}
    \fill[black!20!white] (d0.center) to[out=135,in=-135]  (d1.center) to (d2.center) to (d0.center);     
    \fill[black!20!white] (d2.center) to[out=45,in=180]  (d4.center) to (d3.center) to [out=180, in=-45] (d2.center);
    \fill[black!20!white] (d4.center) to[out=0,in=135]  (d5.center) to (d7.center) to[out=135, in=90] (d4.center);
    \fill[black!20!white] (d3.center) to[out=0,in=-135] (d5.center) to[out=-45, in=135] (d6.center) to [out=-135, in=-135] (d0.center) to [out=-45, in=-90] (d3.center);
    \fill[black!20!white] (d1.center) to[out=135 , in=180] (4,4) to [out=0,in=90] (8,1.5) to [out=-90, in=-45] (d6.center) to [out=45, in=-45] (d7.center) to [out=45, in=45] (d1.center);    
\end{pgfonlayer}

\end{tikzpicture}

\noindent 
By adding a vertex to each unshaded region, and isotoping the knot so the new vertices are arranged horizontally, we get:

\noindent
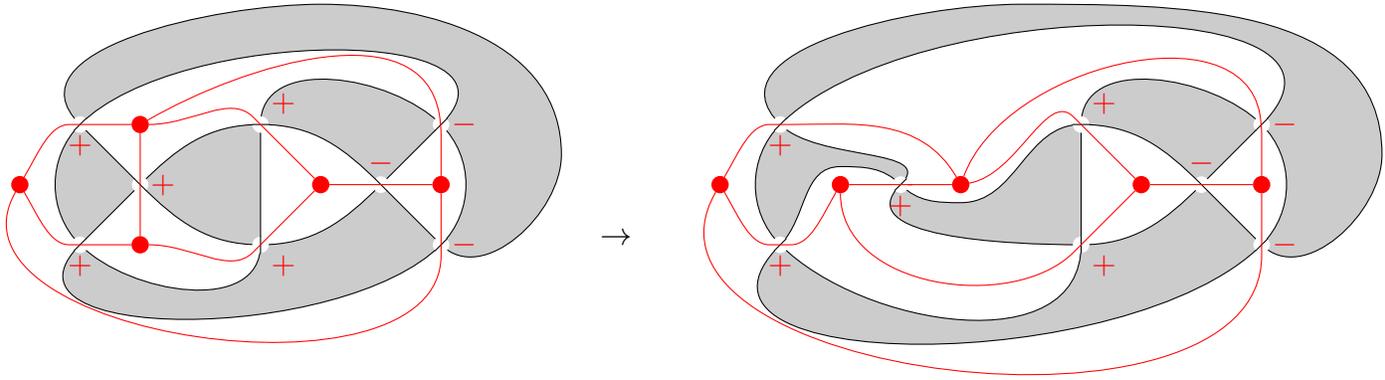
\begin{tikzpicture}[scale=0.8, baseline=0]
  \node[cross] (d0) at (0,0) {};				  \node[cross] (d1) at (0,2) {};
  \node[cross] (d2) at (1,1) {};				  \node[cross] (d3) at (3,0) {};
  \node[cross] (d4) at (3,2) {};				  \node[cross] (d5) at (5,1) {};
  \node[cross] (d6) at (6,0) {};				  \node[cross] (d7) at (6,2) {};
  
\begin{pgfonlayer}{foreground}  
  \draw (d0.45) -- (d0.-135);				  \draw (d1.45) -- (d1.-135);
  \draw (d2.-45) -- (d2.135);				  \draw (d3.90) -- (d3.-90);
  \draw (d4.0) -- (d4.180);					  \draw (d5.-45) -- (d5.135);
  \draw (d6.45) -- (d6.-135);				  \draw (d7.45) -- (d7.-135);
\end{pgfonlayer}

\draw (d0.center) to[out=135, in=-135] (d1.center);			\draw (d0.center) to[out= 45, in=-135] (d2.center);
\draw (d0.center) to[out=-45 , in=-90] (d3.center);			\draw (d0.center) to[out=-135 , in=-135] (d6.center);
\draw (d1.center) to[out=-45 , in=135] (d2.center);			\draw (d1.center) to[out=45 , in=45] (d7.center);
\draw (d1.center) to[out=135 , in=180] (4,4) to [out=0,in=90] (8,1.5) to [out=-90, in=-45](d6.center);
\draw (d2.center) to[out=-45 , in=180] (d3.center);			\draw (d2.center) to[out=45 , in=180] (d4.center);
\draw (d3.center) to[out=90 , in=-90] (d4.center);			\draw (d3.center) to[out=0 , in=-135] (d5.center);
\draw (d4.center) to[out=0 , in=135] (d5.center);				\draw (d4.center) to[out=90 , in=135] (d7.center);
\draw (d5.center) to[out=-45 , in=135] (d6.center);			\draw (d5.center) to[out=45 , in=-135] (d7.center);
\draw (d6.center) to[out=45 , in=-45] (d7.center);

  \node[dot] (a0) at (-1, 1) {};  				  \node[dot] (a1) at (1, 0) {};
  \node[dot] (a2) at (1,2) {};				  \node[dot] (a3) at (4,1) {};
  \node[dot] (a4) at (6,1) {};
  
\begin{pgfonlayer}{foreground}    
  \draw[red, rounded corners] (a0) to [out=-60, in=180] (d0.center) node [below] {$+$} to (a1);  
  \draw[red, rounded corners] (a0) to [out=60, in=180] (d1.center) node [below] {$+$} to [out=0, in=180] (a2);
  \draw[red, rounded corners] (a0) to [out=-120, in=-90] (d6.center)node [right] {$-$} to [out=90, in=-90] (a4);
  \draw[red, rounded corners] (a1) to [out=90, in=-90] (d2.center)node [right] {$+$} to [out=90, in=-90] (a2);
  \draw[red, rounded corners] (a1) to [out=0,in=-135] (d3.center)node [below right] {$+$} to [out=45, in=-135] (a3);
  \draw[red, rounded corners] (a2) to [out=0, in=135] (d4.center)node [above right] {$+$} to [out=-45, in=135] (a3);
  \draw[red, rounded corners] (a2) to [out=30, in=100] (d7.center)node [right] {$-$} to [out=-90, in=90] (a4); 
  \draw[red, rounded corners] (a3) to [out=0,in=180] (d5.center)node [above] {$-$} to [out=0, in=180] (a4); 
\end{pgfonlayer}

\begin{pgfonlayer}{background}
    \fill[black!20!white] (d0.center) to[out=135,in=-135]  (d1.center) to (d2.center) to (d0.center);     
    \fill[black!20!white] (d2.center) to[out=45,in=180]  (d4.center) to (d3.center) to [out=180, in=-45] (d2.center);
    \fill[black!20!white] (d4.center) to[out=0,in=135]  (d5.center) to (d7.center) to[out=135, in=90] (d4.center);
    \fill[black!20!white] (d3.center) to[out=0,in=-135] (d5.center) to[out=-45, in=135] (d6.center) to [out=-135, in=-135] (d0.center) to [out=-45, in=-90] (d3.center);
    \fill[black!20!white] (d1.center) to[out=135 , in=180] (4,4) to [out=0,in=90] (8,1.5) to [out=-90, in=-45] (d6.center) to [out=45, in=-45] (d7.center) to [out=45, in=45] (d1.center);    
\end{pgfonlayer}

\end{tikzpicture}
$\quad \rightarrow$ \hspace*{-5mm}
%
%
\begin{tikzpicture}[scale=0.8, baseline=0]
  \node[cross] (d0) at (0,0) {};				  \node[cross] (d1) at (0,2) {};
  \node[cross] (d2) at (2,1) {};				  \node[cross] (d3) at (5,0) {};
  \node[cross] (d4) at (5,2) {};				  \node[cross] (d5) at (7,1) {};
  \node[cross] (d6) at (8,0) {};				  \node[cross] (d7) at (8,2) {};
  
\begin{pgfonlayer}{foreground}  
  \draw (d0.45) -- (d0.-135);				  \draw (d1.45) -- (d1.-135);
  \draw (d2.-135) -- (d2.45);				  \draw (d3.90) -- (d3.-90);
  \draw (d4.0) -- (d4.180);					  \draw (d5.-45) -- (d5.135);
  \draw (d6.45) -- (d6.-135);				  \draw (d7.45) -- (d7.-135);
\end{pgfonlayer}

\draw (d0.center) to[out=135, in=-135] (d1.center);			
\draw (d0.center) to[out= 45, in=180] (1,1.3) to [out=0, in=135] (d2.center);
\draw (d0.center) to[out=-45 , in=-90] (d3.center);			\draw (d0.center) to[out=-135 , in=-135] (d6.center);
\draw (d1.center) to[out=-45 , in=45] (d2.center);			\draw (d1.center) to[out=45 , in=45] (d7.center);
\draw (d1.center) to[out=135 , in=180] (4,4) to [out=0,in=90] (10,1.5) to [out=-90, in=-45](d6.center);
\draw (d2.center) to[out=-135 , in=180] (d3.center);			
\draw (d2.center) to[out=-45 , in=180] (3,0.7) to [out=0,in=180] (d4.center);
\draw (d3.center) to[out=90 , in=-90] (d4.center);			\draw (d3.center) to[out=0 , in=-135] (d5.center);
\draw (d4.center) to[out=0 , in=135] (d5.center);				\draw (d4.center) to[out=90 , in=135] (d7.center);
\draw (d5.center) to[out=-45 , in=135] (d6.center);			\draw (d5.center) to[out=45 , in=-135] (d7.center);
\draw (d6.center) to[out=45 , in=-45] (d7.center);

  \node[dot] (a0) at (-1, 1) {};  				  \node[dot] (a1) at (1,1) {};
  \node[dot] (a2) at (3,1) {};				  \node[dot] (a3) at (6,1) {};
  \node[dot] (a4) at (8,1) {};
  
\begin{pgfonlayer}{foreground}    
  \draw[red, rounded corners] (a0) to [out=-60, in=180] (d0.center) node [below] {$+$} to [out=0,in=-120] (a1);  
  \draw[red, rounded corners] (a0) to [out=60, in=180] (d1.center) node [below] {$+$} to [out=0, in=120] (a2);
  \draw[red, rounded corners] (a0) to [out=-120, in=-90] (d6.center)node [right] {$-$} to [out=90, in=-90] (a4);
  \draw[red, rounded corners] (a1) to [out=0, in=0] (d2.center)node [below] {$+$} to [out=0, in=180] (a2);
  \draw[red, rounded corners] (a1) to [out=-90,in=-135] (d3.center)node [below right] {$+$} to [out=45, in=-135] (a3);
  \draw[red, rounded corners] (a2) to [out=10, in=135] (d4.center)node [above right] {$+$} to [out=-45, in=135] (a3);
  \draw[red, rounded corners] (a2) to [out=70, in=100] (d7.center)node [right] {$-$} to [out=-90, in=90] (a4); 
  \draw[red, rounded corners] (a3) to [out=0,in=180] (d5.center)node [above] {$-$} to [out=0, in=180] (a4); 
\end{pgfonlayer}

\begin{pgfonlayer}{background}
    \fill[black!20!white] (d0.center) to[out=135,in=-135]  (d1.center) to[out=-45 , in=45] (d2.center) to [out=135, in=0] (1,1.3) to [out=180,in=45] (d0.center);     
    \fill[black!20!white] (d2.center) to[out=-45 , in=180] (3,0.7) to [out=0,in=180] (d4.center) to (d3.center) to [out=180, in=-135] (d2.center);
    \fill[black!20!white] (d4.center) to[out=0,in=135]  (d5.center) to (d7.center) to[out=135, in=90] (d4.center);
    \fill[black!20!white] (d3.center) to[out=0,in=-135] (d5.center) to[out=-45, in=135] (d6.center) to [out=-135, in=-135] (d0.center) to [out=-45, in=-90] (d3.center);
    \fill[black!20!white] (d1.center) to[out=135 , in=180] (4,4) to [out=0,in=90] (10,1.5) to [out=-90, in=-45] (d6.center) to [out=45, in=-45] (d7.center) to [out=45, in=45] (d1.center);    
\end{pgfonlayer}

\end{tikzpicture}

We then extract the graph, and draw every edge as either having all positive, or all negative, $y$-coordinates (except at endpoints, where $y=0$).  Then we may omit the plus labels  on edges with all positive $y$-coordinates, and omit the minus on edges with negative $y$-coordinates.

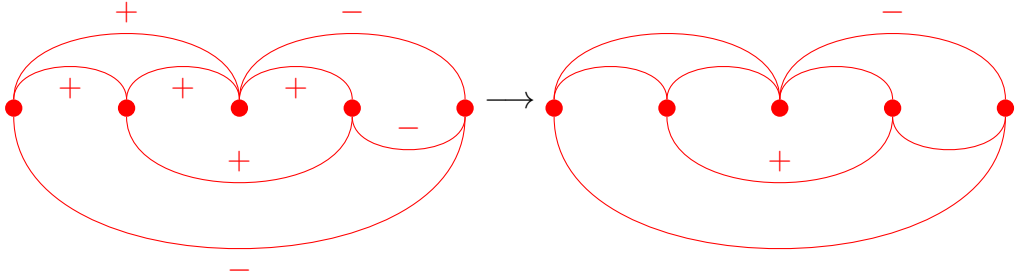
\begin{tikzpicture}[scale=1, baseline=1cm]
  \node[dot] (a0) at (0, 1) {};  				  \node[dot] (a1) at (1.5,1) {};
  \node[dot] (a2) at (3,1) {};				  \node[dot] (a3) at (4.5,1) {};
  \node[dot] (a4) at (6,1) {};
  
\begin{pgfonlayer}{foreground}    
  \draw[red, rounded corners] (a0) to [out=90, in=90] node [below, midway] {$+$} (a1) ;  
  \draw[red, rounded corners] (a0) to [out=90, in=90]  node [midway, above] {$+$} (a2);
  \draw[red, rounded corners] (a0)  to [out=-90, in=-90] node [midway, below] {$-$}(a4) ;
  \draw[red, rounded corners] (a1)  to  [out=90, in=90] node [midway, below] {$+$} (a2);
  \draw[red, rounded corners] (a1) to [out=-90, in=-90] node [midway, above] {$+$}  (a3);
  \draw[red, rounded corners] (a2) to  [out=90, in=90] node [midway, below ] {$+$} (a3);
  \draw[red, rounded corners] (a2) to  [out=90, in=90] node [midway, above] {$-$} (a4); 
  \draw[red, rounded corners] (a3) to  [out=-90, in=-90]node [midway, above] {$-$} (a4); 
\end{pgfonlayer}
\end{tikzpicture}
$\longrightarrow$
\begin{tikzpicture}[scale=1, baseline=1cm]
  \node[dot] (a0) at (0, 1) {};  				  \node[dot] (a1) at (1.5,1) {};
  \node[dot] (a2) at (3,1) {};				  \node[dot] (a3) at (4.5,1) {};
  \node[dot] (a4) at (6,1) {};
  
\begin{pgfonlayer}{foreground}    
  \draw[red, rounded corners] (a0) to [out=90, in=90] node [below, midway] {} (a1) ;  
  \draw[red, rounded corners] (a0) to [out=90, in=90]  node [midway, above] {} (a2);
  \draw[red, rounded corners] (a0)  to [out=-90, in=-90] node [midway, below] {}(a4) ;
  \draw[red, rounded corners] (a1)  to  [out=90, in=90] node [midway, below] {} (a2);
  \draw[red, rounded corners] (a1) to [out=-90, in=-90] node [midway, above] {$+$}  (a3);
  \draw[red, rounded corners] (a2) to  [out=90, in=90] node [midway, below ] {} (a3);
  \draw[red, rounded corners] (a2) to  [out=90, in=90] node [midway, above] {$-$} (a4); 
  \draw[red, rounded corners] (a3) to  [out=-90, in=-90]node [midway, above] {} (a4); 
\end{pgfonlayer}

\end{tikzpicture}

The resulting graph does not yet have the form we observed, above, as coming from a Thompson tree diagram knot.  First, there is one negative edge with positive $y$-coordinates, and one positive edge with negative $y$-coordinates.  Second, all of the vertices fail to have exactly one positive and one negative edge incoming.

To manipulate this graph into Thompson form without changing the knot it represents, we consider the effect of Reidemeister moves I and II on the signed planar graph obtained from unshaded regions.  Reidemeister I moves allow us to delete self-loops and degree-one vertices.  Reidemeister II moves allow us to contract a sequential positive and negative edge, or cancel a pair of parallel positive and negative edges.  

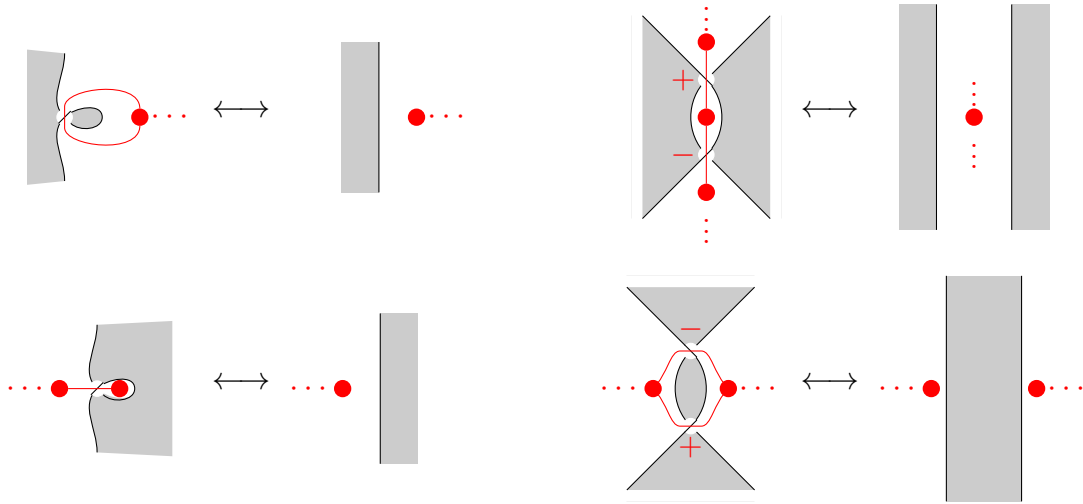
\begin{align*}
\begin{tikzpicture}[baseline=1cm]
      \node[cross] (d0) at (0,1) {};				 
      \node (s) at (0,0) {}; 					
      \node (n) at (0,2) {};
    \begin{pgfonlayer}{foreground}  
      \draw (d0.45) -- (d0.-135);				 
    \end{pgfonlayer}
    \draw (d0.center) to[out=45, in=90] (0.5,1) to [out=-90, in=-45] (d0.-45);			
    \draw (d0.center) to [out=135,in=-90] (n);
    \draw (d0.center) to [out=-135,in=90] (s);
         \node[dot] (a0) at (1,1) {};  				  
    \begin{pgfonlayer}{foreground}    
      \draw[red, rounded corners] (a0) to [out=90, in=90] (d0.center) to [out=-90, in=-90] (a0) ;  
      \node at (a0) [red, right] {$\cdots$};
    \end{pgfonlayer}
    \begin{pgfonlayer}{background}
        \fill[black!20!white] (-0.5,0.1) to (s.90) to [out=90,in=-135] (d0.center) to [out=135,in=-90] (n.-90) to (-0.5,1.9);
      \fill[black!20!white] (d0.center) to[out=45, in=90] (0.5,1) to [out=-90, in=-45] (d0.-45);			
    \end{pgfonlayer}
\end{tikzpicture}
&
 \longleftrightarrow \qquad
\begin{tikzpicture}[baseline=1cm]
      \begin{pgfonlayer}{foreground}  
        \draw (0,0) -- (0,2);
      \end{pgfonlayer}
      \node[dot] (a0) at (0.5,1) {};  		
      \node at (a0) [red, right] {$\cdots$};		  
      \begin{pgfonlayer}{background}
         \fill[black!20!white] (-0.5,0) rectangle (0,2);
      \end{pgfonlayer}
\end{tikzpicture}
&
\begin{tikzpicture}[baseline=1.5cm]
      \node[cross] (d0) at (1,1) {};				  \node[cross] (d1) at (1,2) {};  
      \node (sw) at (0,0) {}; 					  \node (se) at (2,0) {};
      \node (nw) at (0,3) {};					  \node (ne) at (2,3) {};
    \begin{pgfonlayer}{foreground}  
      \draw (d0.45) -- (d0.-135);				  \draw (d1.-45) -- (d1.135);
    \end{pgfonlayer}
    \draw (d0.center) to[out=45, in=-45] (d1.center);			\draw (d0.center) to[out=135, in=-135] (d1.center);
    \draw (d0.center) to (sw);								\draw (d0.center) to (se);
    \draw (d1.center) to (nw);								\draw (d1.center) to (ne);
     \node[dot] (a0) at (1,0.5) {};  				   \node[dot] (a1) at (1,1.5) {};  				
      \node[dot] (a2) at (1,2.5) {};  	  
    \begin{pgfonlayer}{foreground}    
      \draw[red, rounded corners] (a0) -- node [left] {$-$} (a1) -- node [left] {$+$} (a2) ;  
      \node at (a0) [red, below] {$\vdots$};
      \node at (a2) [red, above] {$\vdots$};
    \end{pgfonlayer}
    \begin{pgfonlayer}{background}
        \fill[black!20!white] (nw) -- (sw)--(d0.center) to [out=135, in =-135] (d1.center) -- (nw);
        \fill[black!20!white] (ne) -- (se)--(d0.center) to [out=45, in =-45] (d1.center) -- (ne);
    \end{pgfonlayer}
\end{tikzpicture}
&
 \longleftrightarrow \quad
\begin{tikzpicture}[baseline=1.5cm]
      \begin{pgfonlayer}{foreground}  
        \draw (0,0) -- (0,3);
        \draw (1,0) -- (1,3);
      \end{pgfonlayer}
        \node[dot] (a1) at (0.5,1.5) {};  	
       \begin{pgfonlayer}{foreground}    
      \node at (a1) [red, below] {$\vdots$};
      \node at (a1) [red, above] {$\vdots$};
    \end{pgfonlayer}
      \begin{pgfonlayer}{background}
         \fill[black!20!white] (-0.5,0) rectangle (0,3);
         \fill[black!20!white] (1.5,0) rectangle (1,3);
      \end{pgfonlayer}
\end{tikzpicture}
\\
\begin{tikzpicture}[baseline=1cm]
      \node[cross] (d0) at (0,1) {};				 
      \node (s) at (0,0) {}; 					
      \node (n) at (0,2) {};
    \begin{pgfonlayer}{foreground}  
      \draw (d0.45) -- (d0.-135);				 
    \end{pgfonlayer}
    \draw (d0.center) to[out=45, in=90] (0.5,1) to [out=-90, in=-45] (d0.-45);			
    \draw (d0.center) to [out=135,in=-90] (n);
    \draw (d0.center) to [out=-135,in=90] (s);
     \node[dot] (a0) at (-0.5,1) {};  				  
     \node[dot] (a1) at (0.3,1) {};  				  
    \begin{pgfonlayer}{foreground}    
      \draw[red, rounded corners] (a0) -- (a1) ;  
      \node at (a0) [red, left] {$\cdots$};
    \end{pgfonlayer}
    \begin{pgfonlayer}{background}
        \fill[black!20!white] (1,0.1) -- (s.90) to [out=90,in=-135] (d0.center) to[out=-45, in=-90] (0.5,1) to [out=90, in=45] (d0.center) to [out=135,in=-90] (n.-90) -- (1,1.9) -- (1,0);
    \end{pgfonlayer}
\end{tikzpicture}
\quad
&
 \longleftrightarrow 
\begin{tikzpicture}[baseline=1cm]
      \begin{pgfonlayer}{foreground}  
        \draw (1,0) -- (1,2);
      \end{pgfonlayer}
       \node[dot] (a0) at (0.5,1) {};  
       \begin{pgfonlayer}{foreground}    
      \node at (a0) [red, left] {$\cdots$};
    \end{pgfonlayer}	
      \begin{pgfonlayer}{background}
         \fill[black!20!white] (1.5,0) rectangle (1,2);
      \end{pgfonlayer}
\end{tikzpicture}
&
\begin{tikzpicture}[baseline=1.5cm]
      \node[cross] (d0) at (1,1) {};				  \node[cross] (d1) at (1,2) {};  
      \node (sw) at (0,0) {}; 					  \node (se) at (2,0) {};
      \node (nw) at (0,3) {};					  \node (ne) at (2,3) {};
    \begin{pgfonlayer}{foreground}  
      \draw (d0.45) -- (d0.-135);				  \draw (d1.-45) -- (d1.135);
    \end{pgfonlayer}
    \draw (d0.center) to[out=45, in=-45] (d1.center);			\draw (d0.center) to[out=135, in=-135] (d1.center);
    \draw (d0.center) to (sw);								\draw (d0.center) to (se);
    \draw (d1.center) to (nw);								\draw (d1.center) to (ne);
     \node[dot] (a0) at (0.5,1.5) {};  				   \node[dot] (a1) at (1.5,1.5) {};  	
     \begin{pgfonlayer}{foreground}    
      \draw[red, rounded corners] (a0.center) to [out=-45, in=180] (d0.center) node [below] {$+$} to [out=0, in=-135] (a1.center);  
      \draw[red, rounded corners] (a0.center) to [out=45, in=180] (d1.center) node [above] {$-$} to [out=0, in=135] (a1.center);  
      \node at (a0) [red, left] {$\cdots$};
      \node at (a1) [red, right] {$\cdots$};
    \end{pgfonlayer}
    \begin{pgfonlayer}{background}
        \fill[black!20!white] (se) -- (sw)--(d0.center) -- (se);
        \fill[black!20!white] (ne) -- (nw)-- (d1.center) -- (ne);
        \fill[black!20!white] (d0.center) to [out=135, in=-135] (d1.center) to [out=-45, in=45] (d0.center);
    \end{pgfonlayer}
\end{tikzpicture}
&
\longleftrightarrow  
\begin{tikzpicture}[baseline=1.5cm]
      \begin{pgfonlayer}{foreground}  
        \draw (0,0) -- (0,3);
        \draw (1,0) -- (1,3);
      \end{pgfonlayer}
      \node[dot] (a0) at (-0.2, 1.5) {};  				   \node[dot] (a1) at (1.2,1.5) {};  	
     \begin{pgfonlayer}{foreground}    
      \node at (a0) [red, left] {$\cdots$};
      \node at (a1) [red, right] {$\cdots$};
    \end{pgfonlayer}
      \begin{pgfonlayer}{background}
         \fill[black!20!white] (1,0) rectangle (0,3);
      \end{pgfonlayer}
\end{tikzpicture}
\end{align*}

A note on conventions:  a vertex with adjacent ellipses indicates that the vertex may have additional edges coming out of it; the additional edges are unaffected by the depicted move.

Remember, our goal is to make a planar signed graph into a planar signed graph coming from a Thompson tree.  We can force positive edges up and negative edges down, using a Reidemeister II move. In the examples that follow, vertices that are created by a Reidemeister I or II move are highlighted in pink during that step. 

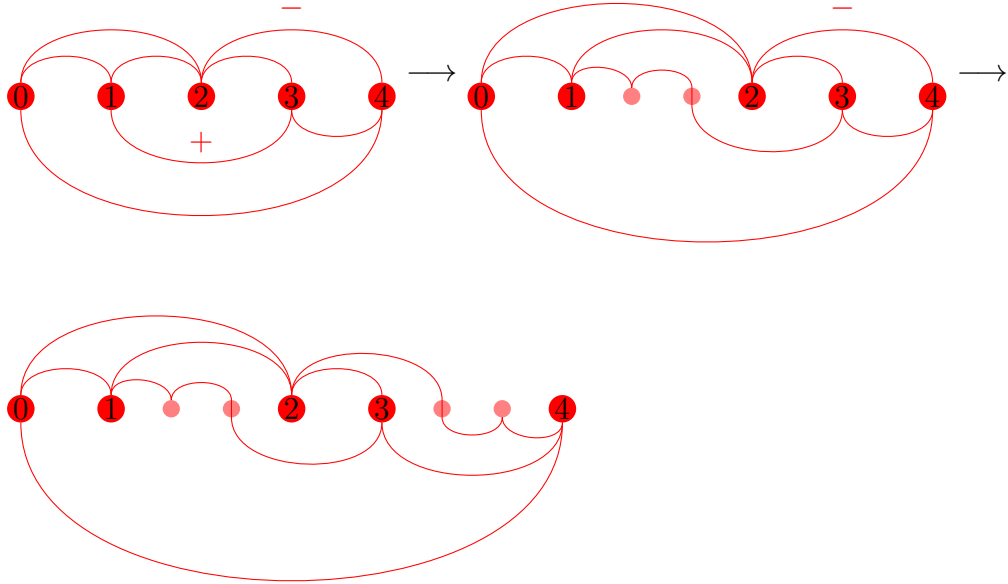
\begin{tikzpicture}[xscale=0.8, yscale=0.8, baseline=1cm]
  \node[dot, inner sep=0] (a0) at (0, 1) {$0$};  				  \node[dot, inner sep=0] (a1) at (1.5,1) {$1$};
  \node[dot, inner sep=0] (a2) at (3,1) {$2$};				  \node[dot, inner sep=0] (a3) at (4.5,1) {$3$};
  \node[dot, inner sep=0] (a4) at (6,1) {$4$};
  
\begin{pgfonlayer}{foreground}    
  \draw[red, rounded corners] (a0) to [out=90, in=90] node [below, midway] {} (a1) ;  
  \draw[red, rounded corners] (a0) to [out=90, in=90]  node [midway, above] {} (a2);
  \draw[red, rounded corners] (a0)  to [out=-90, in=-90] node [midway, below] {}(a4) ;
  \draw[red, rounded corners] (a1)  to  [out=90, in=90] node [midway, below] {} (a2);
  \draw[red, rounded corners] (a1) to [out=-90, in=-90] node [midway, above] {$+$}  (a3);
  \draw[red, rounded corners] (a2) to  [out=90, in=90] node [midway, below ] {} (a3);
  \draw[red, rounded corners] (a2) to  [out=90, in=90] node [midway, above] {$-$} (a4); 
  \draw[red, rounded corners] (a3) to  [out=-90, in=-90]node [midway, above] {} (a4); 
\end{pgfonlayer}

\end{tikzpicture}
$\longrightarrow$
\begin{tikzpicture}[xscale=0.8, yscale=0.8, baseline=1cm]
  \node[dot, inner sep=0] (a0) at (0, 1) {$0$};  				  \node[dot, inner sep=0] (a1) at (1.5,1) {$1$};
  \node[pinkdot] (b1) at (2.5,1) {};					  
  \node[pinkdot] (b2) at (3.5,1) {};	
  \node[dot, inner sep=0] (a2) at (4.5,1) {$2$};			 	 \node[dot, inner sep=0] (a3) at (6,1) {$3$};
  \node[dot, inner sep=0] (a4) at (7.5,1) {$4$};
  
\begin{pgfonlayer}{foreground}    
  \draw[red, rounded corners] (a0) to [out=90, in=90] node [below, midway] {} (a1) ;  
  \draw[red, rounded corners] (a0) to [out=90, in=90]  node [midway, above] {} (a2);
  \draw[red, rounded corners] (a0)  to [out=-90, in=-90] node [midway, below] {}(a4) ;
  \draw[red, rounded corners] (a1)  to  [out=90, in=90] node [midway, below] {} (a2);
  \draw[red, rounded corners] (a2) to  [out=90, in=90] node [midway, below ] {} (a3);
  \draw[red, rounded corners] (a2) to  [out=90, in=90] node [midway, above] {$-$} (a4); 
  \draw[red, rounded corners] (a3) to  [out=-90, in=-90]node [midway, above] {} (a4); 
  
  \draw [red] (a1.90) to [out=90, in=90] (b1.90) to [out=90, in=90] (b2.90) to (b2.-90) to [out=-90, in=-90] (a3.-90);
\end{pgfonlayer}
\end{tikzpicture}
$\longrightarrow$

\begin{tikzpicture}[xscale=0.8, yscale=0.8, baseline=1cm]
  \node[dot, inner sep=0] (a0) at (0, 1) {$0$};  				  \node[dot, inner sep=0] (a1) at (1.5,1) {$1$};
  \node[pinkdot] (b1) at (2.5,1) {};					  
  \node[pinkdot] (b2) at (3.5,1) {};	
  \node[dot, inner sep=0] (a2) at (4.5,1) {$2$};			 	 \node[dot, inner sep=0] (a3) at (6,1) {$3$};
    \node[pinkdot] (c1) at (7,1) {};					  
    \node[pinkdot] (c2) at (8,1) {};	
  \node[dot, inner sep=0] (a4) at (9,1) {$4$};
  
\begin{pgfonlayer}{foreground}    
  \draw[red, rounded corners] (a0) to [out=90, in=90] node [below, midway] {} (a1) ;  
  \draw[red, rounded corners] (a0) to [out=90, in=90]  node [midway, above] {} (a2);
  \draw[red, rounded corners] (a0)  to [out=-90, in=-90] node [midway, below] {}(a4) ;
  \draw[red, rounded corners] (a1)  to  [out=90, in=90] node [midway, below] {} (a2);
  \draw[red, rounded corners] (a2) to  [out=90, in=90] node [midway, below ] {} (a3);
  \draw[red, rounded corners] (a3) to  [out=-90, in=-90] (a4); 
  
  \draw [red] (a1.90) to [out=90, in=90] (b1.90) to [out=90, in=90] (b2.90) to (b2.-90) to [out=-90, in=-90] (a3.-90);  \draw[red] (a2) to  [out=90, in=90] (c1.90) to (c1.-90) to  [out=-90, in=-90] (c2.-90) to  [out=-90, in=-90]   (a4); 
\end{pgfonlayer}
\end{tikzpicture}

Now we have a signed planar graph where the signs of all edges are accurately reflected in their $y$-coordinates.  It still remains to manipulate our graph so that each vertex (other than the initial one) has one incoming positive edge and one incoming negative edge.  First we deal with vertices that have too many incoming edges, by splitting them up using a Reidemeister II move.  In our running example, this looks like

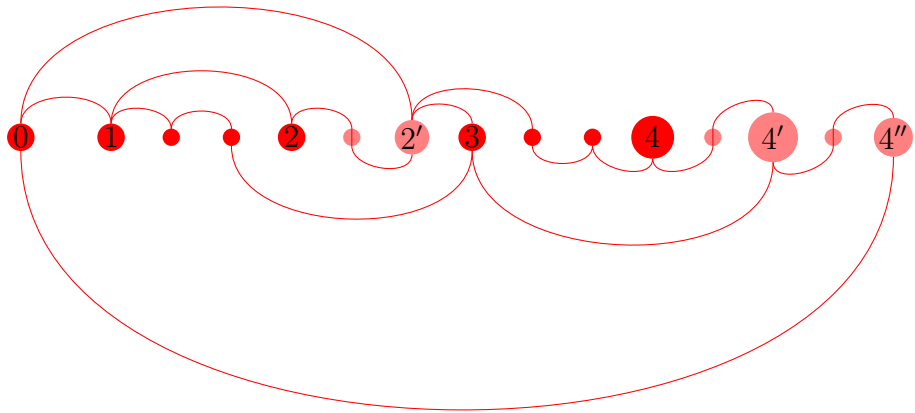
\begin{tikzpicture}[xscale=0.8, yscale=0.8, baseline=1cm]
  \node[dot, inner sep=0] (a0) at (0, 1) {$0$};  				  \node[dot, inner sep=0] (a1) at (1.5,1) {$1$};
  \node[dot] (b1) at (2.5,1) {};							  \node[dot] (b2) at (3.5,1) {};	
  \node[dot, inner sep=0] (a2) at (4.5,1) {$2$};		
    \node[pinkdot] (d1) at (5.5,1) {};	
  \node[pinkdot, inner sep =0] (d2) at (6.5,1) {$2'$};	
  \begin{scope}[xshift=1.5cm]
  \node[dot, inner sep=0] (a3) at (6,1) {$3$};
    \node[dot] (c1) at (7,1) {};							\node[dot] (c2) at (8,1) {};	
    \node[dot, inner sep] (e1) at (9,1) {$4$};					\node[pinkdot] (e2) at (10,1) {};
     \node[pinkdot, inner sep] (e3) at (11,1) {$4'$};				\node[pinkdot] (e4) at (12,1) {};
      \begin{scope}[xshift=4cm]
      \node[pinkdot, inner sep=0] (a4) at (9,1) {$4''$};
      \end{scope}
  \end{scope}

\begin{pgfonlayer}{foreground}    
  \draw[red, rounded corners] (a0) to [out=90, in=90] node [below, midway] {} (a1) ;  
  \draw[red, rounded corners] (a0) to [out=90, in=90]  node [midway, above] {} (d2);
  \draw[red, rounded corners] (a0)  to [out=-90, in=-90] node [midway, below] {}(a4) ;
  \draw[red, rounded corners] (a1)  to  [out=90, in=90] node [midway, below] {} (a2);
  \draw[red, rounded corners] (d2) to  [out=90, in=90] node [midway, below ] {} (a3);
  \draw[red, rounded corners] (a3) to  [out=-90, in=-90] (e3); 
  
  \draw [red] (a1.90) to [out=90, in=90] (b1.90) to [out=90, in=90] (b2.90) to (b2.-90) to [out=-90, in=-90] (a3.-90);
  \draw[red] (d2) to  [out=90, in=90] (c1.90) to (c1.-90) to  [out=-90, in=-90] (c2.-90) to  [out=-90, in=-90]   (e1); 
  \draw [red] (a2) to  [out=90, in=90] (d1.90) to (d1.-90) to [out=-90, in=-90] (d2);
    \draw[red] (e1) to  [out=-90, in=-90] (e2.-90) to (e2.90) to  [out=90, in=90]   (e3); 
  \draw[red] (e3) to  [out=-90, in=-90] (e4.-90) to (e4.90)  to  [out=90, in=90]   (a4); 

\end{pgfonlayer}
\end{tikzpicture}

\noindent where the multiple edges coming into a vertex $k$ in the earlier diagram have been split between vertices $k$ and $k'$ in the new diagram.

Our final step is to give out extra edges as needed.  If we use a Reidemeister I move to create a new vertex and edge, and then a Reidemeister II move to create a double bond with the vertex to its left, then our new vertex also has incoming edges of both signs. 

\noindent
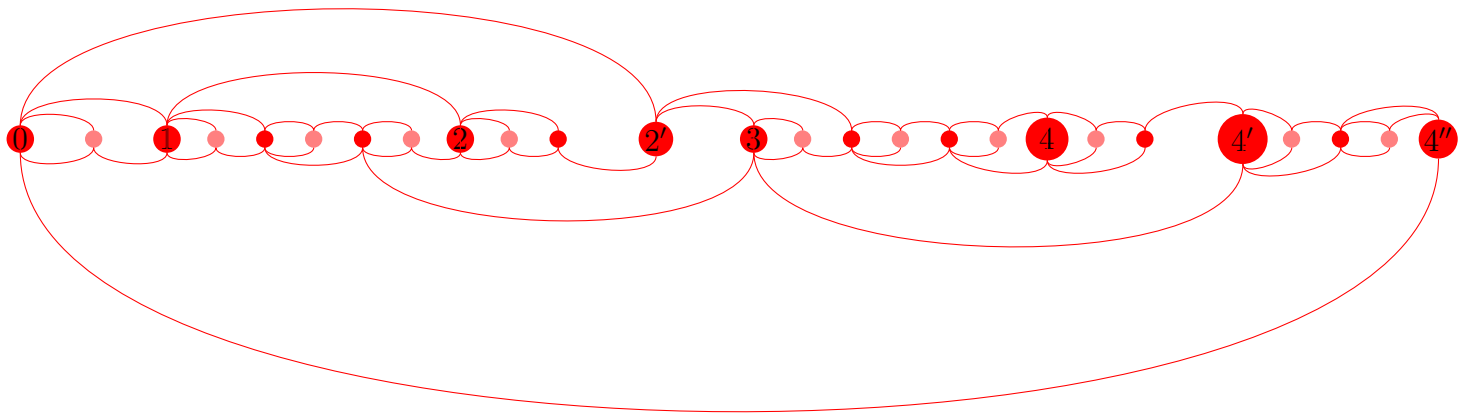
\begin{tikzpicture}[xscale=1.3, yscale=0.8, baseline=1cm]
  \node[dot, inner sep=0] (a0) at (0, 1) {$0$};  				  \node[dot, inner sep=0] (a1) at (1.5,1) {$1$};
  \node[dot] (b1) at (2.5,1) {};					  \node[dot] (b2) at (3.5,1) {};	
  \node[dot, inner sep=0] (a2) at (4.5,1) {$2$};		
    \node[dot] (d1) at (5.5,1) {};	
  \node[dot, inner sep =0] (d2) at (6.5,1) {$2'$};	
  \begin{scope}[xshift=1.5cm]
  \node[dot, inner sep=0] (a3) at (6,1) {$3$};
    \node[dot] (c1) at (7,1) {};					  \node[dot] (c2) at (8,1) {};	
    \node[dot, inner sep] (e1) at (9,1) {$4$};		\node[dot] (e2) at (10,1) {};
     \node[dot, inner sep] (e3) at (11,1) {$4'$};		\node[dot] (e4) at (12,1) {};
      \begin{scope}[xshift=4cm]
      \node[dot, inner sep=0] (a4) at (9,1) {$4''$};
      \end{scope}
  \end{scope}
  
    \node[pinkdot] (f1) at (.75,1) {};							\node[pinkdot] (f2) at (2,1) {};
    \node[pinkdot] (f3) at (3,1) {};							\node[pinkdot] (f4) at (4,1) {};
    \node[pinkdot] (f5) at (5,1) {};							\node[pinkdot] (f6) at (8,1) {};
    \node[pinkdot] (f7) at (9,1) {};							\node[pinkdot] (f8) at (10,1) {};
    \node[pinkdot] (f9) at (11,1) {};							\node[pinkdot] (f10) at (13,1) {};
    \node[pinkdot] (f11) at (14,1) {};

\begin{pgfonlayer}{foreground}    
  \draw[red, rounded corners] (a0) to [out=90, in=90] node [below, midway] {} (a1) ;  
  \draw[red, rounded corners] (a0) to [out=90, in=90]  node [midway, above] {} (d2);
  \draw[red, rounded corners] (a0)  to [out=-90, in=-90] node [midway, below] {}(a4) ;
  \draw[red, rounded corners] (a1)  to  [out=90, in=90] node [midway, below] {} (a2);
  \draw[red, rounded corners] (d2) to  [out=90, in=90] node [midway, below ] {} (a3);
  \draw[red, rounded corners] (a3) to  [out=-90, in=-90] (e3); 
  
  \draw [red] (a1.90) to [out=90, in=90] (b1.90) to (b1.-90) to [out=-90, in=-90] (b2.-90) to [out=-90, in=-90] (a3);
  \draw[red] (d2) to  [out=90, in=90] (c1.90) to (c1.-90) to  [out=-90, in=-90] (c2.-90) to  [out=-90, in=-90]   (e1); 
  \draw [red] (a2) to  [out=90, in=90] (d1.90) to (d1.-90) to [out=-90, in=-90] (d2);
    \draw[red] (e1) to  [out=-90, in=-90] (e2.-90) to (e2.90) to  [out=90, in=90]   (e3); 
  \draw[red] (e3) to  [out=-90, in=-90] (e4.-90) to (e4.90)  to  [out=90, in=90]   (a4); 

   \draw[red] (a1) to [out=-90, in=-90] (f1.-90) to [out=-90, in=-90] (a0.-90) to (a0.90) to[out=90, in=90] (f1.90);
   \draw[red] (b1) to [out=-90, in=-90] (f2.-90) to [out=-90, in=-90] (a1.-90) to (a1.90) to[out=90, in=90] (f2.90);
   \draw[red] (b2) to [out=90, in=90] (f3.90) to [out=90, in=90] (b1.90) to (b1.-90) to[out=-90, in=-90] (f3.-90);
   \draw[red] (a2) to [out=-90, in=-90] (f4.-90) to [out=-90, in=-90] (b2.-90) to (b2.90) to[out=90, in=90] (f4.90);
   \draw[red] (d1) to [out=-90, in=-90] (f5.-90) to [out=-90, in=-90] (a2.-90) to (a2.90) to[out=90, in=90] (f5.90);
   \draw[red] (c1) to [out=-90, in=-90] (f6.-90) to [out=-90, in=-90] (a3.-90) to (a3.90) to[out=90, in=90] (f6.90);
   \draw[red] (c2) to [out=90, in=90] (f7.90) to [out=90, in=90] (c1.90) to (c1.-90) to[out=-90, in=-90] (f7.-90);
   \draw[red] (e1) to [out=90, in=90] (f8.90) to [out=90, in=90] (c2.90) to (c2.-90) to[out=-90, in=-90] (f8.-90);
   \draw[red] (e2) to [out=90, in=90] (f9.90) to [out=90, in=90] (e1.90) to (e1.-90) to[out=-90, in=-90] (f9.-90);
   \draw[red] (e4) to [out=90, in=90] (f10.90) to [out=90, in=90] (e3.90) to (e3.-90) to[out=-90, in=-90] (f10.-90);
   \draw[red] (a4) to [out=90, in=90] (f11.90) to [out=90, in=90] (e4.90) to (e4.-90) to[out=-90, in=-90] (f11.-90);

\end{pgfonlayer}
\end{tikzpicture}

Thus we have created a signed planar graph of the sort that can arise from a Thompson group knot.  More than that, it \emph{does} arise from a Thompson group element.  The edges direct us how to place the crossings:  one per edge, with the edge passing through the unshaded regions and the sign of the crossing based on the sign of the edge.  Then the faces of the graph tell us how to connect the crossings: around the perimeter of each face (and the outer face).  The result is a knot which is a small isotopy away from having the Thompson group form.

\section{Multiplication of Tangles}

\begin{defn}\label{T}
The Thompson diagram $T(x)$ is the diagram from \ref{mchairs}:
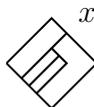
$$\begin{tikzpicture}[baseline=0, scale=0.2]
    \chair[thick](-2,2)(1:2);

	\draw (0,0)--(-3,3)--(-6,0);
	\draw (0,0)--(-3,-3)--(-6,0);

    \node at (-2,2) [above right] {$x$};
\end{tikzpicture}$$

\noindent The Thompson diagram $T(x_1,x_2,\cdots,x_{n-1},x_n, x_{n+1})$ for $x_i \geq 1$ is defined recursively: $T(x_1,x_2,\cdots,x_n,x_{n+1})$ is the rescaled concatenation of two Thompson diagrams, where the lower right edge of $T(x_{n+1})$ is joined to the upper left edge of $R(T(x_1,x_2,\cdots,x_n-1))$.  This process is illustrated below:

$$\begin{tikzpicture}[baseline=0, scale=0.2]

    \chair[thick](-6,6)(5:2);

    \filldraw[black!20!white] (0,0)--(-7,-7)--(-11,-3)--(-4,4)--(0,0);
    \node at (-6,-1) {$A$};
	 \draw (0,0)--(-7,7)--(-14,0);
	 \draw (0,0)--(-7,-7)--(-14,0);
	 \draw (-11,-3)--(-4,4);

    \node at (-6,6) [above right] {$x_n$};

    \node at (-7,-12) [below] {$T(x_1,x_2,\cdots,x_n)$};

\end{tikzpicture}
\rightarrow
\begin{tikzpicture}[baseline=0, scale=0.2,yscale=-1]

    \chair[thick](-6,6)(5:2);

    \filldraw[black!20!white] (0,0)--(-7,-7)--(-11,-3)--(-4,4)--(0,0);
    \node[yscale=-1] at (-6,-1) {$A$};
	 \draw (0,0)--(-7,7)--(-14,0);
	 \draw (0,0)--(-7,-7)--(-14,0);
	 \draw (-11,-3)--(-4,4);

    \node at (-6,6) [below right] {$x_n-1$};

    \node at (-7,12) [below] {$R(T(x_1,x_2,\cdots,x_n-1))$};

\end{tikzpicture}
\rightarrow
\begin{tikzpicture}[baseline=0, scale=0.2,yscale=-1]

    \chair[thick](-6,6)(5:2);
    \uchair[thick](-9,-9)(8:2);

    \filldraw[black!20!white] (0,0)--(-7,-7)--(-11,-3)--(-4,4)--(0,0);
    \node[yscale=-1] at (-6,-1) {$A$};
	 \draw (0,0)--(-10,10)--(-20,0);
	 \draw (0,0)--(-10,-10)--(-20,0);
	 \draw (-11,-3)--(-4,4);
	 \draw (-7,-7)--(-17,3);

    \node at (-6,6) [below right] {$x_n-1$};
    \node at (-9,-9) [above right] {$x_{n+1}$};

    \node at (-10,12) [below] {$R(T(x_1,x_2,\cdots,x_n, x_{n+1}))$};

\end{tikzpicture}
$$

\noindent where $A$ is the composition of former parts of the recursion relation.  Note that the diagram for $T(x_1,x_2,\cdots,x_{n-1},x_n)$ has one less chair in each position than the notation suggests (except for the final position), for reasons that will be clear soon.  

\end{defn}

\begin{eg}
Let us demonstrate the construction of T(3,4,2,5).

We begin with $T(3)$, from which we construct $T(3,4)$, and next $T(3,4,2)$:

\noindent
 
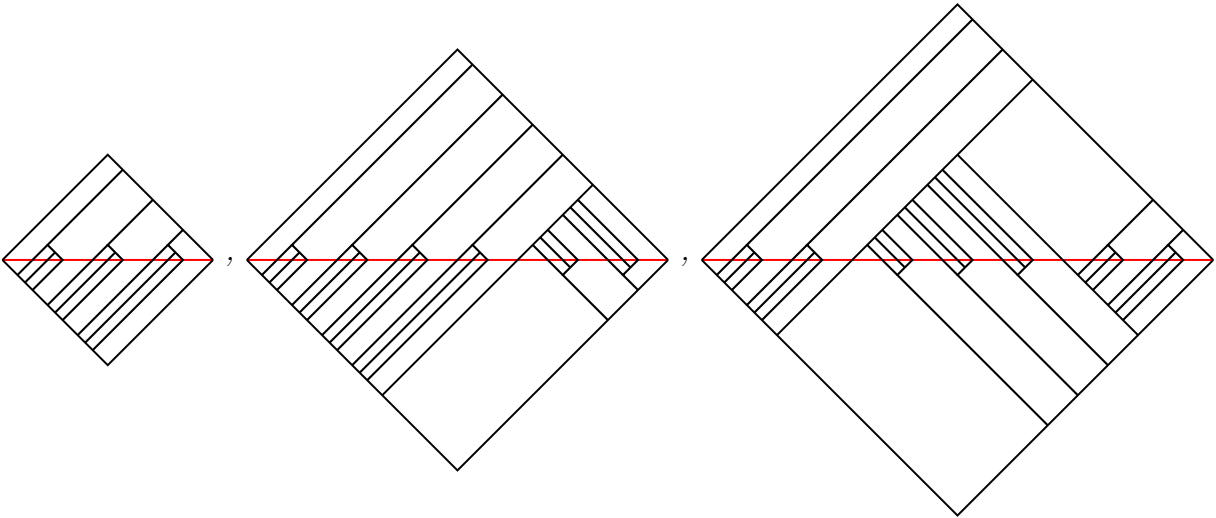
\begin{tikzpicture}[baseline=0, scale=0.2, yscale=-1]
    \draw[red](-14,0)--(0,0);
 
    \uchair[thick](-2,-2)(1:6);
	\uchair[thick](-4,-4)(3:4);
	\uchair[thick](-6,-6)(5:2);

	\draw (0,0)--(-7,7)--(-14,0);
	\draw (0,0)--(-7,-7)--(-14,0);
\end{tikzpicture} 
,
\begin{tikzpicture}[baseline=0, scale=0.2]
    \draw[red](-28,0)--(0,0);
 
    \uchair[thick](-2,-2)(1:4);
	\uchair[thick](-4,-4)(3:2);
	
	\chair[thick](-7,7)(6:8);
	\chair[thick](-9,9)(8:6);
	\chair[thick](-11,11)(10:4);
	\chair[thick](-13,13)(12:2);

	\draw (0,0)--(-14,14)--(-28,0);
	\draw (0,0)--(-14,-14)--(-28,0);
	\draw (-5,5)--(-19,-9);
\end{tikzpicture}
,
\begin{tikzpicture}[baseline=0, scale=0.2, yscale=-1]
    \draw[red](-34,0)--(0,0);
 
    \uchair[thick](-2,-2)(1:4);
	\uchair[thick](-4,-4)(3:2);
	
	\chair[thick](-7,7)(6:6);
	\chair[thick](-9,9)(8:4);
	\chair[thick](-11,11)(10:2);
	
	\uchair[thick](-14,-14)(13:4);
	\uchair[thick](-16,-16)(15:2);

	\draw (0,0)--(-17,17)--(-34,0);
	\draw (0,0)--(-17,-17)--(-34,0);
	\draw (-5,5)--(-17,-7);
	\draw (-12,-12) -- (-17,-7) -- (-29,5);
\end{tikzpicture}

\noindent
And finally we produce $T(3,4,2,5)$:

\noindent
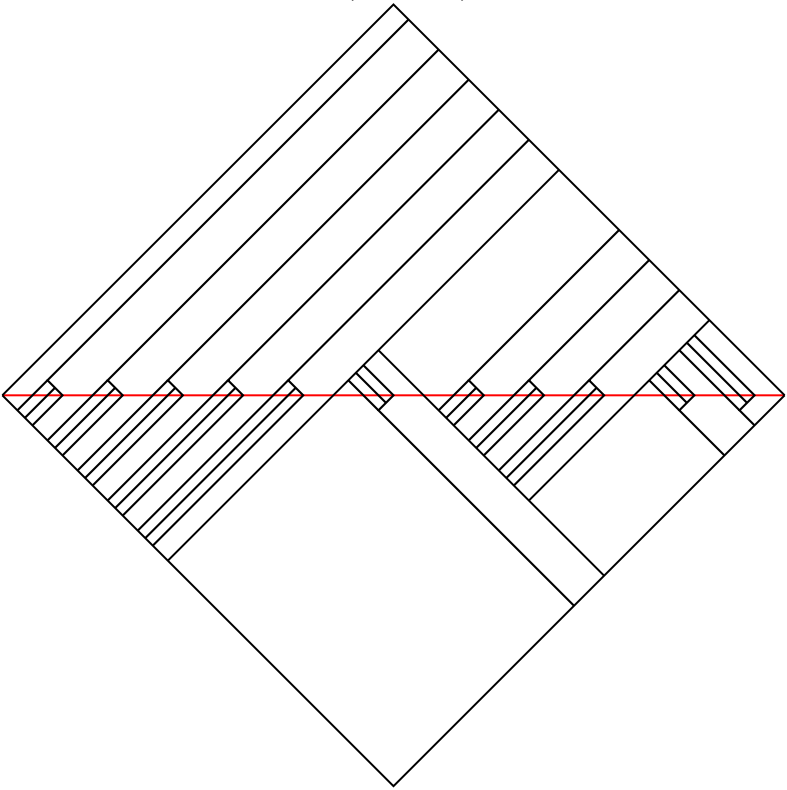
\begin{tikzpicture}[baseline=0, scale=0.2]
    \draw[red](-52,0)--(0,0);
 
    \uchair[thick](-2,-2)(1:4);
	\uchair[thick](-4,-4)(3:2);
	
	\chair[thick](-7,7)(6:6);
	\chair[thick](-9,9)(8:4);
	\chair[thick](-11,11)(10:2);
	
	\uchair[thick](-14,-14)(13:2);

	\chair[thick](-17,17)(16:10);
	\chair[thick](-19,19)(18:8);
	\chair[thick](-21,21)(20:6);
	\chair[thick](-23,23)(22:4);
	\chair[thick](-25,25)(24:2);

	\draw (0,0)--(-26,26)--(-52,0);
	\draw (0,0)--(-26,-26)--(-52,0);
	\draw (-5,5)--(-17,-7);
	\draw (-12,-12) -- (-17,-7) -- (-27,3);
	\draw (-15,15) -- (-27,3)--(-41,-11);
\end{tikzpicture}

\end{eg}

\begin{defn}\label{K1thrun}
When $n$ is odd, $K(x_1,x_2,\cdots,x_n)$ is the knot diagram  on the left, and when $n$ is even, the diagram on the right:

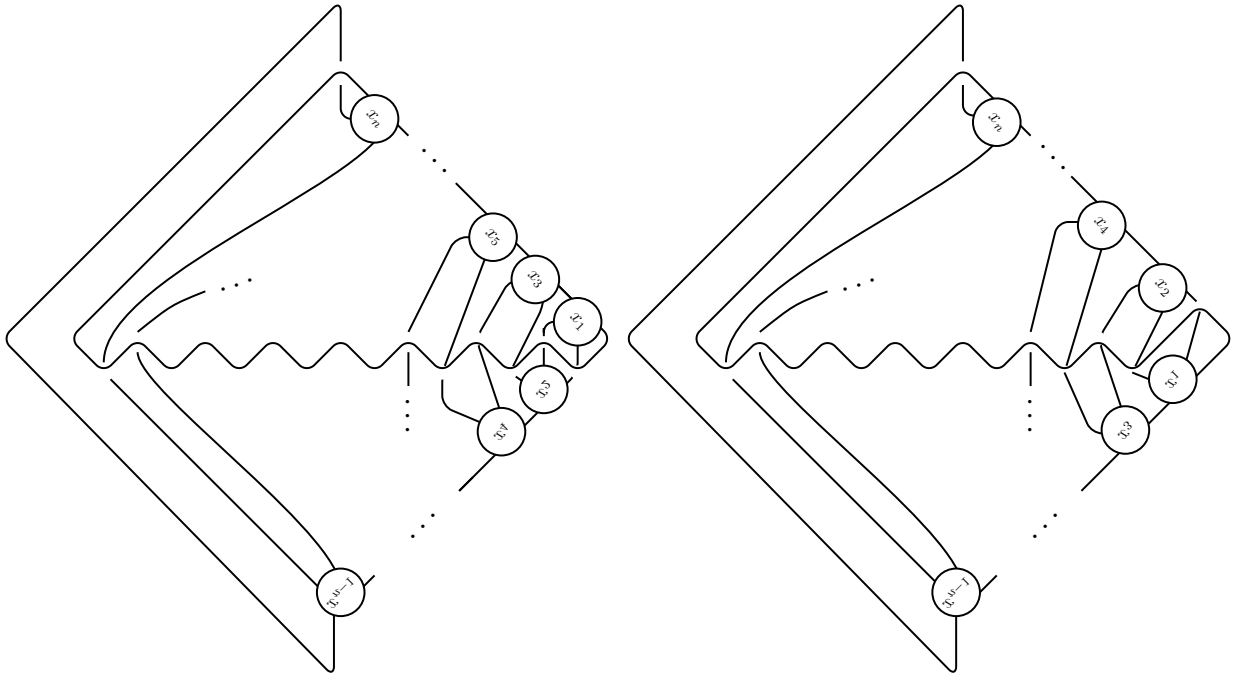
\begin{tikzpicture}[scale=0.9]
\draw(4,-4)--(0.6,-0.6);%

\draw (0.5,-0.35) .. controls (0.5,) and (5,2.5) .. (4.7,3.4);

\draw[rounded corners](7.5,0.25)--(7.5,-0.5)--(7,-1)--(5.75,-2.25);%
\draw[rounded corners](7.5,0.25)--(7,0.25)--(7,0)--(7,-0.85);

\draw[rounded corners](7.75,0.25)--(7.4,0.6);

\draw[rounded corners](7.45,0.55)--(7,1);

\draw[rounded corners](6.875,0.875)--(6.5,0.875)--(6,0)--(6.5,-1.5);

\draw(7.5,0.5)--(5.7,2.3);

\draw[rounded corners](7,0.8)--(7,0.5)--(6.5,-0.5);
\draw[rounded corners](6.5,-0.5)--(7,-0.85);
%

\draw[rounded corners](6.25,1.5)--(5.5,-0.5)--(5.5,-1)--(6.45,-1.4);

\draw[rounded corners](6.25,1.5)--(5.7,1.5)--(5,0.1);
\draw(5,-0.2)--(5,-0.7);
\node at (5,-1.1)[rotate=90]{$\cdots$};

\draw(5.95,-2.05)--(5.75,-2.25);

\node at (5.25,-2.75)[rotate=45]{$\cdots$};
\node at (5.4,2.6)[rotate=-45]{$\cdots$};

\fill[color=white] (7,-0.05) circle(2.5pt);
\fill[color=white] (6.5,-0.5) circle(3pt);
\fill[color=white] (6,0) circle(3pt);
\fill[color=white] (5.5,-0.5) circle(3pt);
\fill[color=white] (7.5,-0.5) circle(3pt);
\draw[rounded corners](7.5,0.5)--(8,0)--(7.5,-0.5)--(7,0)--(6.5,-0.5)--(6,0)--(5.5,-0.5)--(5,0)--(4.5,-0.5)--(4,0)--(3.5,-0.5)--(3,0)--(2.5,-0.5)--(2,0)--(1.5,-0.5)--(1,0)--(0.5,-0.5)--(0,0)--(4,4)--(5,3);

\draw[rounded corners](4,3.75)--(4,3.25)--(4.5,3.25);

\draw[rounded corners](4,4.1)--(4,5)--(-1,0)-- (3.9,-5)--(3.9,-4);

\draw(4,-3.8) .. controls (4.2,-3) and (1,-0.8) .. (1,-0.2);
\draw(1,0.1) .. controls (1.5,0.5) .. (2,0.7);

\draw(4,-4)--(4.5,-3.5);

\node at (2.5,0.8)[rotate=20]{$\cdots$};

\filldraw [fill=white](7.5,0.25) node[scale=0.6,rotate=-45]{$x_1$} circle(10pt);
\filldraw [fill=white](6.875,0.875) node[scale=0.6,rotate=-45]{$x_3$} circle(10pt);
\filldraw [fill=white](6.25,1.5) node[scale=0.6,rotate=-45]{$x_5$} circle(10pt);
\filldraw [fill=white](4.5,3.25) node[scale=0.6,rotate=-45]{$x_n$} circle(10pt);
\filldraw [fill=white](7,-0.75) node[scale=0.6,rotate=-135,xscale=-1,yscale=1]{${x_2}$} circle(10pt);
\filldraw [fill=white](6.375,-1.375) node[scale=0.6,rotate=-135,xscale=-1,yscale=1]{${x_4}$} circle(10pt);
\filldraw [fill=white](4,-3.75) node[scale=0.6,rotate=-135,xscale=-1,yscale=1]{${x_{n-1}}$} circle(10pt);
\end{tikzpicture}
\begin{tikzpicture}[scale=0.9]
\draw(4,-4)--(0.6,-0.6);
\draw (0.5,-0.35) .. controls (0.5,) and (5,2.5) .. (4.7,3.4);

\draw[rounded corners](7.25,-0.75)--(7.25,-0.5)--(7.5,0.4);

\draw(7.45,0.55)--(5.6,2.4);

\draw[rounded corners](7,0.8)--(6.5,0.8)--(6,0)--(6.5,-1.5);

\draw[rounded corners](7,0.8)--(7,0.5)--(6.5,-0.5);

\draw[rounded corners](6.1,1.725)--(6.1,1.5)--(5.5,-0.5);
\draw[rounded corners](6.1,1.725)--(5.4,1.725)--(5,0.1);
\draw(5,-0.2)--(5,-0.7);
\node at (5,-1.1)[rotate=90]{$\cdots$};

\node at (5.25,-2.75)[rotate=45]{$\cdots$};
\node at (5.35,2.65)[rotate=-45]{$\cdots$};

\fill[color=white] (6.5,-0.5) circle(3pt);
\fill[color=white] (6,0) circle(3pt);
\fill[color=white] (5.5,-0.5) circle(3pt);
\draw(6.5,-0.5)--(7,-0.7);
\draw[rounded corners](5.5,-0.5)--(5.9,-1.4)--(6.45,-1.4);
\draw[rounded corners](5.75,-2.25)--(6.5,-1.5)--(8,0)--(7.75,0.25)--(7.5,0.5)--(7,0)--(6.5,-0.5)--(6,0)--(5.5,-0.5)--(5,0)--(4.5,-0.5)--(4,0)--(3.5,-0.5)--(3,0)--(2.5,-0.5)--(2,0)--(1.5,-0.5)--(1,0)--(0.5,-0.5)--(0,0)--(4,4)--(5,3);

\draw(4,-3.8) .. controls (4.2,-3) and (1,-0.8) .. (1,-0.2);
\draw[rounded corners](4,3.75)--(4,3.3)--(4.3,3.3);
\draw[rounded corners](4,4.1)--(4,5)--(-1,0)-- (3.9,-5)--(3.9,-4);

\draw(4,-4)--(4.5,-3.5);
\draw(1,0.1) .. controls (1.5,0.5) .. (2,0.7);

\node at (2.5,0.8)[rotate=20]{$\cdots$};

\filldraw [fill=white](6.95,0.75) node[scale=0.6,rotate=-45]{$x_2$} circle(10pt);
\filldraw [fill=white](6.05,1.675) node[scale=0.6,rotate=-45]{$x_4$} circle(10pt);
\filldraw [fill=white](4.5,3.2) node[scale=0.6,rotate=-45]{$x_n$} circle(10pt);
\filldraw [fill=white] (7.1,-0.6) node[scale=0.6,rotate=-135,xscale=-1,yscale=1]{$x_1$} circle(10pt);
\filldraw [fill=white](6.4,-1.35)  node[scale=0.6,rotate=-135,xscale=-1,yscale=1]{$x_3$} circle(10pt);
\filldraw [fill=white](3.9,-3.75)  node[scale=0.6,rotate=-135,xscale=-1,yscale=1]{$x_{n-1}$} circle(10pt);
\end{tikzpicture}

\end{defn}

The $K$'s are related by a recursion relation similar to the one in Definition $\ref{T}$. This will be explained in more detail later.

\begin{eg}

The knot  $K(3,4,2,5)$:

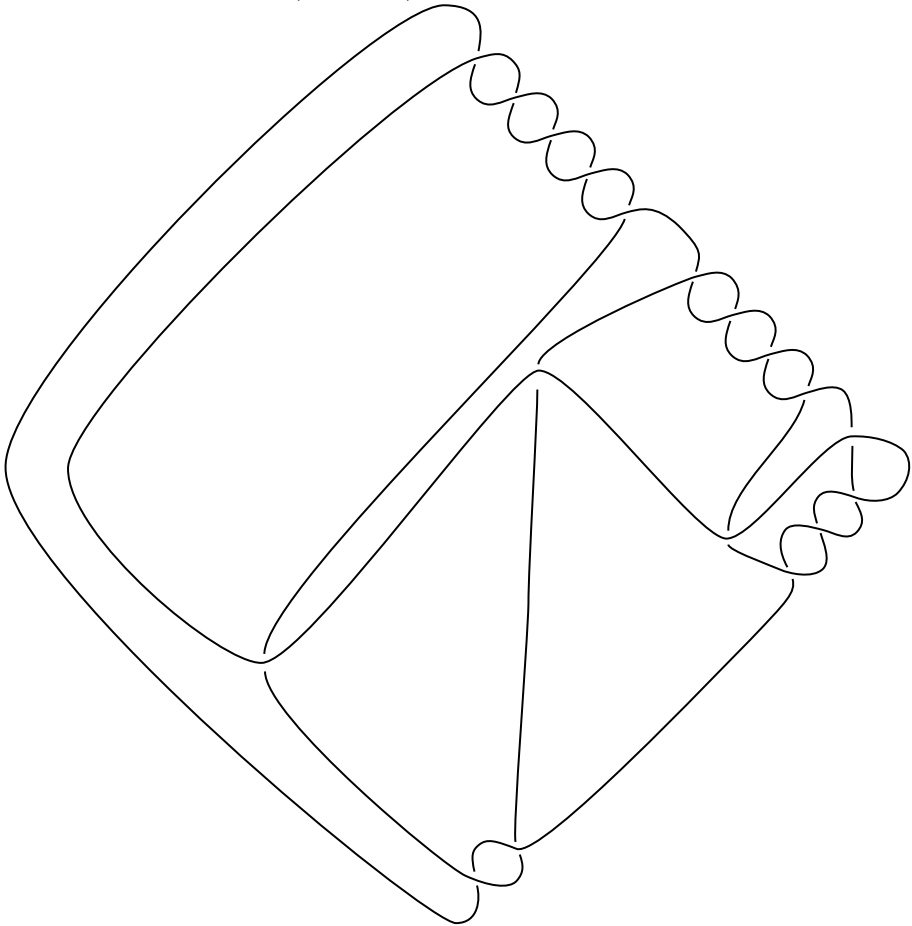
\tikzset{every picture/.style={line width=0.75pt}} 

\begin{tikzpicture}[x=0.75pt,y=0.75pt,yscale=-1,xscale=1,scale=.9]

\draw    (527,352.38) .. controls (529,359.94) and (523.86,364.96) .. (503.33,386.5) .. controls (484.33,405.5) and (387,508.67) .. (372.43,503.81) .. controls (360.71,499.24) and (355.75,497.79) .. (351.29,501.53) .. controls (347.5,504.54) and (346.71,507.38) .. (348.75,516.54) ;
\draw    (542.71,326.67) .. controls (544.14,331.81) and (548.71,341.81) .. (544.14,346.67) .. controls (538.14,352.67) and (524.3,348.95) .. (521,347.53) .. controls (514.39,344.68) and (491.11,336.44) .. (491,333.11) ;
\draw    (560.56,277.78) .. controls (560.56,286.22) and (559.51,295.84) .. (561.29,302.96) ;
\draw    (372,79.79) .. controls (374.73,70.52) and (375.36,67.03) .. (369.91,61.58) .. controls (364.64,56.49) and (359.73,57.58) .. (349.73,60.49) .. controls (301.36,76.67) and (121,255.92) .. (121,290.42) .. controls (120.5,328.42) and (206.4,399.97) .. (229.4,399.47) .. controls (254.9,399.47) and (368.17,235.32) .. (384.83,235.54) .. controls (404.83,235.32) and (473.22,329.33) .. (489.89,329.78) .. controls (505.89,329.78) and (544.11,272.89) .. (560.11,272.44) .. controls (571,272) and (583,274.89) .. (589.22,280.44) .. controls (594.78,286.67) and (592.61,297.27) .. (585.89,304.22) .. controls (581.02,309.14) and (570.81,309.78) .. (562.14,306.67) .. controls (553,303.53) and (545.57,301.53) .. (541.57,306.1) .. controls (537.29,311.24) and (539.29,315.81) .. (540.71,321.24) ;
\draw    (494.83,201) .. controls (496.17,196.56) and (499.28,190.56) .. (492.83,183.67) .. controls (487.28,178.33) and (479.5,181.22) .. (472.39,183) .. controls (464.83,185) and (385.08,218.54) .. (384.58,232.04) ;
\draw    (535.89,244.22) .. controls (537.22,239.78) and (541.51,233.9) .. (535.07,227.02) .. controls (529.51,221.68) and (521.73,224.57) .. (514.62,226.35) .. controls (507.07,228.35) and (499.8,233.73) .. (492.84,227.02) .. controls (487.93,222.15) and (489.07,216.35) .. (492.18,207.68) ;
\draw    (560.11,267.33) .. controls (560.11,263.11) and (560.37,250.75) .. (555.22,246.67) .. controls (550.07,242.59) and (539.51,246.79) .. (535.95,247.68) .. controls (528.4,249.68) and (521.13,255.07) .. (514.18,248.35) .. controls (509.26,243.48) and (510.4,237.68) .. (513.51,229.02) ;
\draw    (491,325.33) .. controls (490.67,301) and (527.67,275.33) .. (533.67,252) ;
\draw    (350.17,524.33) .. controls (352.17,531.89) and (351.17,545.34) .. (338.5,545.42) .. controls (323.17,546.34) and (86,352.17) .. (86,290.17) .. controls (84,236.17) and (292.43,30.1) .. (331.86,30.67) .. controls (355.29,30.96) and (353,45.38) .. (351,56.38) ;
\draw    (384,246.17) .. controls (384,264.17) and (379,343.92) .. (379,366.42) .. controls (378.5,389.42) and (370,492.5) .. (371.78,499.89) ;
\draw    (514.89,222.22) .. controls (516.22,217.78) and (520.51,211.9) .. (514.07,205.02) .. controls (508.51,199.68) and (500.73,202.57) .. (493.62,204.35) .. controls (486.07,206.35) and (478.8,211.73) .. (471.84,205.02) .. controls (466.93,200.15) and (468.07,194.35) .. (471.18,185.68) ;
\draw    (231.4,404.27) .. controls (231.07,429.6) and (330.9,513) .. (343.67,518.83) .. controls (356.43,524.67) and (368.16,526.52) .. (372.44,521.67) .. controls (376.73,516.81) and (376.14,512.67) .. (374.14,506.96) ;
\draw    (472.86,180.02) .. controls (475.57,170.31) and (476.86,168.31) .. (463.43,154.6) .. controls (450,141.74) and (441.33,144.97) .. (434.22,146.75) .. controls (426.67,148.75) and (419.4,154.13) .. (412.44,147.42) .. controls (407.53,142.55) and (408.67,136.75) .. (411.78,128.08) ;
\draw    (435.38,142.84) .. controls (436.71,138.4) and (441,132.53) .. (434.55,125.64) .. controls (429,120.3) and (421.22,123.19) .. (414.11,124.97) .. controls (406.55,126.97) and (399.29,132.36) .. (392.33,125.64) .. controls (387.41,120.77) and (388.55,114.97) .. (391.67,106.3) ;
\draw    (413.6,121.62) .. controls (414.93,117.18) and (419.22,111.3) .. (412.78,104.42) .. controls (407.22,99.08) and (399.44,101.97) .. (392.33,103.75) .. controls (384.78,105.75) and (377.51,111.13) .. (370.55,104.42) .. controls (365.64,99.55) and (367.64,94.21) .. (370.75,85.54) ;
\draw    (392.82,100.29) .. controls (394.16,95.84) and (398.44,89.97) .. (392,83.08) .. controls (386.44,77.75) and (378.67,80.64) .. (371.55,82.42) .. controls (364,84.42) and (356.73,89.8) .. (349.78,83.08) .. controls (344.86,78.21) and (346,72.42) .. (349.11,63.75) ;
\draw    (231,394.42) .. controls (231.44,357.31) and (423.67,178.83) .. (433,150.5) ;
\draw    (524.14,345.81) .. controls (519.29,336.96) and (519.19,329.81) .. (522.9,324.95) .. controls (526.62,320.1) and (537.69,323.27) .. (544.14,324.96) .. controls (550.59,326.64) and (557.86,331.53) .. (563,326.1) .. controls (568.14,320.67) and (565.76,315.81) .. (562.33,309.24) ;

\end{tikzpicture}

\end{eg}

\begin{defn}
When $n$ is odd, the tangle $\tilde{K}_l(x_1,x_2,\cdots, x_{n-1},x_n)$ is defined as the diagram on the left below, and when $n$ is even, it is defined as the diagram on the right:

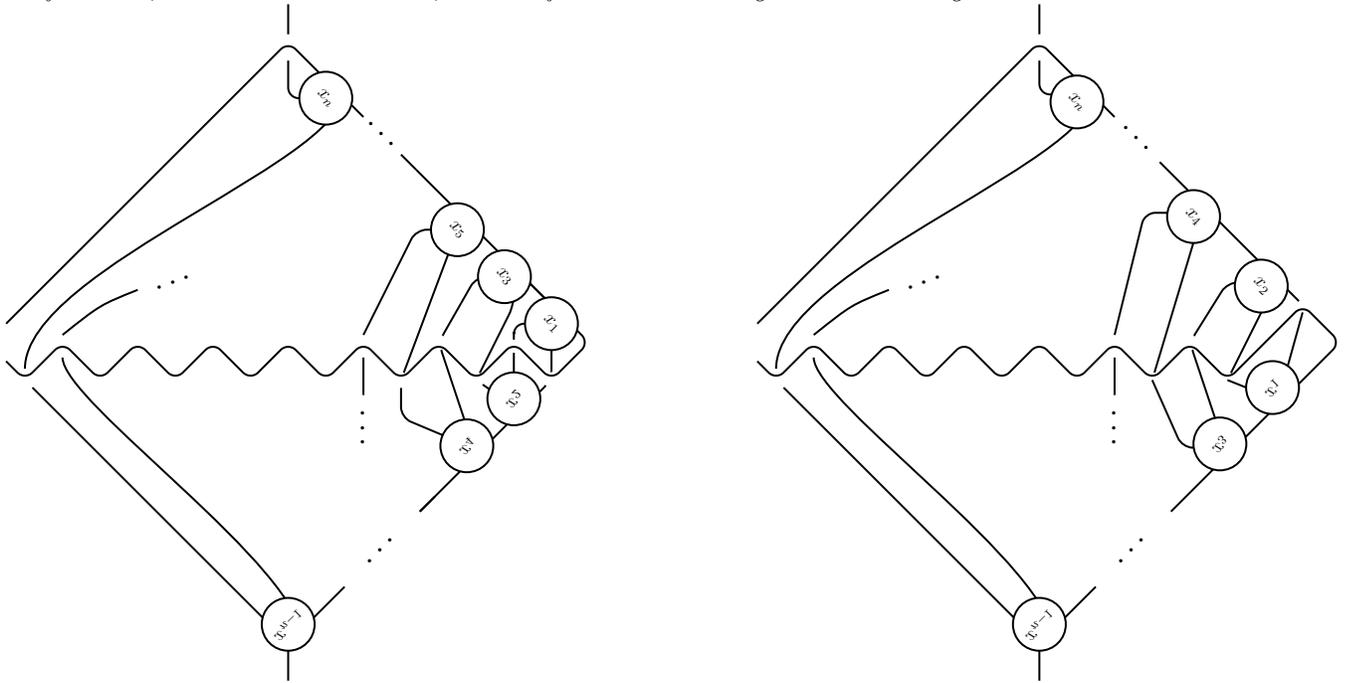
\begin{tikzpicture}
\draw[rounded corners](4.75,-3.25)--(4,-4)--(0.6,-0.6);%
\draw (0.5,-0.35) .. controls (0.5,) and (5,2.5) .. (4.7,3.4);

\draw[rounded corners](7.5,0.25)--(7.5,-0.5)--(7,-1)--(5.75,-2.25);%
\draw[rounded corners](7.5,0.25)--(7,0.25)--(7,0)--(7,-0.85);

\draw[rounded corners](7.75,0.25)--(7.4,0.6);

\draw[rounded corners](7.45,0.55)--(7,1);

\draw[rounded corners](6.875,0.875)--(6.5,0.875)--(6,0)--(6.5,-1.5);

\draw(7.5,0.5)--(5.7,2.3);

\draw[rounded corners](7,0.8)--(7,0.5)--(6.5,-0.5);
\draw[rounded corners](6.5,-0.5)--(7,-0.85);
%

\draw[rounded corners](6.25,1.5)--(5.5,-0.5)--(5.5,-1)--(6.45,-1.4);

\draw[rounded corners](6.25,1.5)--(5.7,1.5)--(5,0.1);

\draw(5.7,2.3)--(5.5,2.5);
\draw(5.95,-2.05)--(5.75,-2.25);

\node at (5.25,-2.75)[rotate=45]{$\cdots$};
\node at (5.25,2.75)[rotate=-45]{$\cdots$};

\fill[color=white] (7,-0.05) circle(2.5pt);
\fill[color=white] (6.5,-0.5) circle(3pt);
\fill[color=white] (6,0) circle(3pt);
\fill[color=white] (5.5,-0.5) circle(3pt);
\fill[color=white] (7.5,-0.5) circle(3pt);
\draw[rounded corners](7.5,0.5)--(8,0)--(7.5,-0.5)--(7,0)--(6.5,-0.5)--(6,0)--(5.5,-0.5)--(5,0)--(4.5,-0.5)--(4,0)--(3.5,-0.5)--(3,0)--(2.5,-0.5)--(2,0)--(1.5,-0.5)--(1,0)--(0.5,-0.5)--(0.25,-0.25);

\draw[rounded corners](4,3.75)--(4,3.25)--(4.5,3.25);

\draw[rounded corners](4,4.1)--(4,4.5);
\draw(4,-4.5)--(4,-4);

\draw(1,0.1) .. controls (1.5,0.5) .. (2,0.7);

\node at (2.5,0.8)[rotate=20]{$\cdots$};

\draw(5,-0.2)--(5,-0.7);
\node at (5,-1.1)[rotate=90]{$\cdots$};

\draw(4.2,-4) .. controls (4.2,-3) and (1,-0.8) .. (1,-0.2);

\draw[rounded corners](0.25,0.25)--(4,4)--(5,3);

\filldraw [fill=white](7.5,0.25) node[scale=0.6,rotate=-45]{$x_1$} circle(10pt);
\filldraw [fill=white](6.875,0.875) node[scale=0.6,rotate=-45]{$x_3$} circle(10pt);
\filldraw [fill=white](6.25,1.5) node[scale=0.6,rotate=-45]{$x_5$} circle(10pt);
\filldraw [fill=white](4.5,3.25) node[scale=0.6,rotate=-45]{$x_n$} circle(10pt);
\filldraw [fill=white](7,-0.75) node[scale=0.6,rotate=-135,xscale=-1,yscale=1]{${x_2}$} circle(10pt);
\filldraw [fill=white](6.375,-1.375) node[scale=0.6,rotate=-135,xscale=-1,yscale=1]{${x_4}$} circle(10pt);
\filldraw [fill=white](4,-3.75) node[scale=0.6,rotate=-135,xscale=-1,yscale=1]{${x_{n-1}}$} circle(10pt);
\end{tikzpicture}
\hspace{2cm}
\begin{tikzpicture}
\draw[rounded corners](4.75,-3.25)--(4,-4)--(0.6,-0.6);

\draw (0.5,-0.35) .. controls (0.5,) and (5,2.5) .. (4.7,3.4);

\draw[rounded corners](7.25,-0.75)--(7.25,-0.5)--(7.5,0.4);

\draw(7.45,0.55)--(5.6,2.4);

\draw[rounded corners](7,0.8)--(6.5,0.8)--(6,0)--(6.5,-1.5);

\draw[rounded corners](7,0.8)--(7,0.5)--(6.5,-0.5);

\draw[rounded corners](6.1,1.725)--(6.1,1.5)--(5.5,-0.5);
\draw[rounded corners](6.1,1.725)--(5.4,1.725)--(5,0.1);
\draw(5,-0.2)--(5,-0.7);
\node at (5,-1.1)[rotate=90]{$\cdots$};

\node at (5.25,-2.75)[rotate=45]{$\cdots$};
\node at (5.3,2.7)[rotate=-45]{$\cdots$};

\fill[color=white] (6.5,-0.5) circle(3pt);
\fill[color=white] (6,0) circle(3pt);
\fill[color=white] (5.5,-0.5) circle(3pt);
\draw(6.5,-0.5)--(7,-0.7);
\draw[rounded corners](5.5,-0.5)--(5.9,-1.4)--(6.45,-1.4);
\draw[rounded corners](5.75,-2.25)--(6.5,-1.5)--(8,0)--(7.75,0.25)--(7.5,0.5)--(7,0)--(6.5,-0.5)--(6,0)--(5.5,-0.5)--(5,0)--(4.5,-0.5)--(4,0)--(3.5,-0.5)--(3,0)--(2.5,-0.5)--(2,0)--(1.5,-0.5)--(1,0)--(0.5,-0.5)--(0.25,-0.25);

\draw[rounded corners](4,3.75)--(4,3.3)--(4.3,3.3);
\draw[rounded corners](4,4.1)--(4,4.5);

\draw(4,-4.5)--(4,-4);
\draw(4.2,-4) .. controls (4.2,-3) and (1,-0.8) .. (1,-0.2);

\draw[rounded corners](0.25,0.25)--(4,4)--(5,3);
\draw(1,0.1) .. controls (1.5,0.5) .. (2,0.7);

\node at (2.5,0.8)[rotate=20]{$\cdots$};

\filldraw [fill=white](6.95,0.75) node[scale=0.6,rotate=-45]{$x_2$} circle(10pt);
\filldraw [fill=white](6.05,1.675) node[scale=0.6,rotate=-45]{$x_4$} circle(10pt);
\filldraw [fill=white](4.5,3.2) node[scale=0.6,rotate=-45]{$x_n$} circle(10pt);
\filldraw [fill=white] (7.1,-0.6) node[scale=0.6,rotate=-135,xscale=-1,yscale=1]{$x_1$} circle(10pt);
\filldraw [fill=white](6.4,-1.35)  node[scale=0.6,rotate=-135,xscale=-1,yscale=1]{$x_3$} circle(10pt);
\filldraw [fill=white](4,-3.75)  node[scale=0.6,rotate=-135,xscale=-1,yscale=1]{$x_{n-1}$} circle(10pt);
\end{tikzpicture}

Note that $\tilde{K}_l$ is K but snipping the strand from top to bottom, and snipping the left-hand rounded corner.  In general, snipping $K$ is similar to snipping $T$, except when $K$ is snipped, we always cut the strand that goes from top to bottom in addition to the (rounded) corner indicated by the subscript.  

\end{defn}

This is almost the tangle with Conway notation $x_1x_2\cdots x_{n-1}x_n$.  Observe a zig zagging strand in the middle of the diagram that goes over all other strands it meets.  If we move this strand over so that it sits outside of $x_1$ and doesn't cross any other strands, it will be the Conway notation tangle mentioned above.
 
Note that tangle multiplication is not associative, so $x_1,x_2,x_3$ means $(x_1x_2)x_3$.

\begin{thm}
$\psi'(\tilde{T}_l(x_1,x_2,\cdots,x_{n-1},x_n))=\tilde{K}_l(x_1,x_2,\cdots,x_{n-1},x_n)$

\end{thm}
\begin{proof}
We will illustrate the case where $n$ is even; the case for odd $n$ works analogously.

Base case: 
$$\tilde{T}_l(x_1) = \hspace{12pt}
\begin{tikzpicture}[baseline=0, scale=0.4]
    \chair[thick](-2,2)(1:2);

	\draw (0,0)--(-3,3)--(-5.5,0.5);
	\draw (0,0)--(-3,-3)--(-5.5,-0.5);

    \node at (-2,2) [above right] {$x_1$};
\end{tikzpicture}
\hspace{6pt}; \hspace{24pt}
\psi'(\tilde{T}_l(x_1)) \hspace{12pt} = \hspace{12pt}
\begin{tikzpicture}[baseline,scale=0.5]
\draw[rounded corners](3,1)--(2,2)--(0.25,0.25);
\draw[rounded corners](0.25,-0.5)--(2,-2)--(4,0)--(3,1);
\draw[rounded corners](3,1)--(2,1)--(2,1.7);
\draw(3,1) .. controls (3,-0.5) .. (2,-1.7);
\filldraw[fill=white] (2.8,0.8) circle(20pt);
\node at (2.8,0.8)[rotate=-45,scale=1.2]{$x_1$};
\draw[rounded corners](2,2.2)--(2,3);
\draw(2,-3)--(2,-2.2);

\end{tikzpicture}
= \tilde{K}(x_1).
$$

Inductive hypothesis: Suppose we've shown the statement true for some $k\geq 1$, so we've shown:

$$
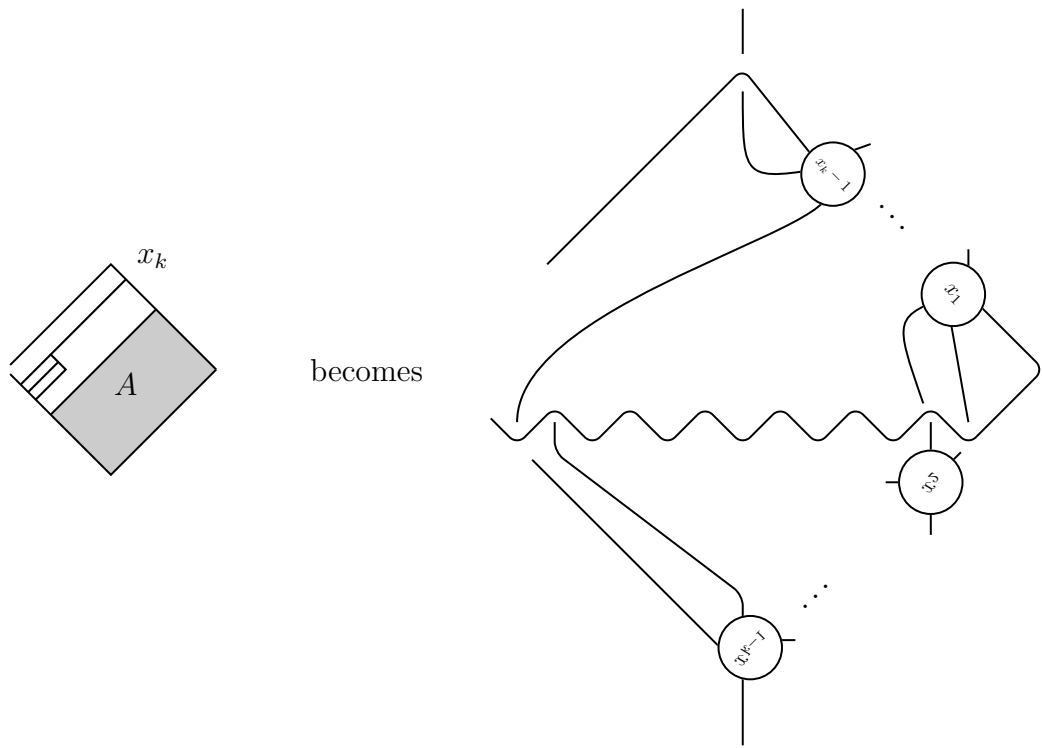
\begin{tikzpicture}[baseline=0, scale=0.2]

    \chair[thick](-6,6)(5:2);

    \filldraw[black!20!white] (0,0)--(-7,-7)--(-11,-3)--(-4,4)--(0,0);
    \node at (-6,-1) {$A$};
	 \draw (0,0)--(-7,7)--(-13.7,0.3);
	 \draw (0,0)--(-7,-7)--(-13.7,-0.3);
	 \draw (-11,-3)--(-4,4);

    \node at (-6,6) [above right] {$x_k$};

\node at (10,0) {becomes};
\begin{scope}[scale=5,shift={(3,0)}]
\draw(4,4.2)--(4,4.8);
\draw[rounded corners](4,-5)--(4,-3)--(1.5,-1.1)--(1.5,-0.7);
\draw(5.2,2.7) .. controls (4,2.5) .. (4,3.7);
\draw(6.7,1)--(7,-0.7);
\draw(1,-0.7) .. controls (1,1)  and (6,2) .. (5,2.5);
\draw[color=white](6.5,1.5)--(5.7,2.3);
\draw[rounded corners](5.7,3)--(5,2.75)--(4,4)--(1.4,1.4);
\draw[rounded corners](0.65,-0.65)--(1,-1)--(1.5,-0.5)--(2,-1)--(2.5,-0.5)--(3,-1)--(3.5,-0.5)--(4,-1)--(4.5,-0.5)--(5,-1)--(5.5,-0.5)--(6,-1)--(6.5,-0.5)--(7,-1)--(8,0)--(7,1)--(7,1.6);
\draw[rounded corners](1.2,-1.2)--(4,-4);
\draw(6.9,-1.1)--(6.5,-1.5);
\draw[rounded corners](6.5,-0.7)--(6.5,-1.5);
\draw(5.9,-1.5)--(6.5,-1.5)--(6.5,-2.2);
\begin{scope}[shift={(-0.7,0.5)}]
\draw(7.5,0.5) .. controls (6.7,0.2) .. (7.1,-0.95);
\filldraw [color=white] (7.5,0.5) circle (12pt);
\draw (7.5,0.5) circle (12pt);
\node at (7.5,0.5)[scale=0.7,rotate=-45]{$x_1$};
\end{scope}
\begin{scope}[shift={(0.7,-0.9)}]
\filldraw [color=white] (4.5,3.5) circle (12pt);
\draw (4.5,3.5) circle (12pt);
\node at (4.5,3.5)[scale=0.5,rotate=-45]{$x_k-1$};
\end{scope}
\draw(4.7,-3.6)--(4,-3.6);
\filldraw[color=white] (4.1,-3.7) circle (12pt);
\draw (4.1,-3.7) circle (12pt);
\node at (4.1,-3.7)[scale=0.7,rotate=45,xscale=1,yscale=-1]{$x_{k-1}$};
\filldraw[color=white] (6.5,-1.5) circle (12pt);
\draw (6.5,-1.5) circle (12pt);
\node at (6.5,-1.5)[rotate=45,scale=0.7,xscale=1,yscale=-1]{$x_2$};
\node at (6,2)[rotate=-45]{$\cdots$};
\node at (5,-3)[rotate=45]{$\cdots$};

\end{scope}
\end{tikzpicture}
$$
.

Reflecting the previous picture, and bringing in the $k+1$st chair set, we have 

$$
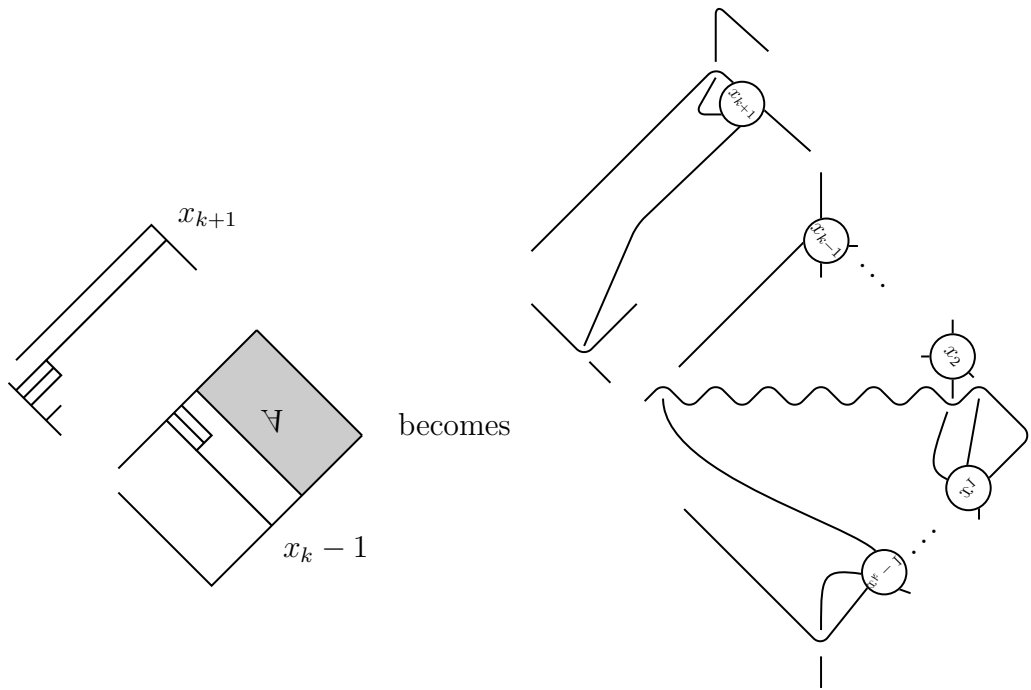
\begin{tikzpicture}[baseline=0, scale=0.2,yscale=-1]

   \chair[thick](-6,6)(5:2);

    \filldraw[black!20!white] (0,0)--(-7,-7)--(-11,-3)--(-4,4)--(0,0);
    \node[yscale=-1] at (-6,-1) {$A$};
	 \draw (0,0)--(-10,10)--(-16.2,3.8);
	\draw(0,0)--(-7,-7);

	 \draw (-11,-3)--(-4,4);
	 \draw (-7,-7)--(-16.2,2.2);
	 	 
	 \begin{scope}[shift={(-4,-4)}]
	 \draw (-7,-7)--(-10,-10)--(-19,-1);
	 \draw(-19.5,0.5)--(-16,4);
	 \draw(-17,3)--(-16,2);
	 
	 	    \node at (-9,-9) [above right] {$x_{k+1}$};
	 	    \uchair[thick](-9,-9)(8:2);
	 \end{scope}

    \node at (-6,6) [below right] {$x_k-1$};

\end{tikzpicture}
\;
\text{ becomes }
\;
\begin{tikzpicture}[xscale=1,yscale=-1,scale=0.7, baseline=0]
\draw(4,4.2)--(4,4.8);
\draw(4,-5)--(4,-3.6);
\draw(5.2,2.7) .. controls (4,2.5) .. (4,3.7);
\draw(6.7,1)--(7,-0.7);
\draw(1,-0.7) .. controls (1,1)  and (6,2) .. (5,2.5);
\draw[color=white](6.5,1.5)--(5.7,2.3);
\draw[rounded corners](5.7,3)--(5,2.75)--(4,4)--(1.4,1.4);
\draw[rounded corners](0.65,-0.65)--(1,-1)--(1.5,-0.5)--(2,-1)--(2.5,-0.5)--(3,-1)--(3.5,-0.5)--(4,-1)--(4.5,-0.5)--(5,-1)--(5.5,-0.5)--(6,-1)--(6.5,-0.5)--(7,-1)--(8,0)--(7,1)--(7,1.6);
\draw[rounded corners](1.3,-1.3)--(4,-4);
\draw(6.9,-1.1)--(6.5,-1.5);
\draw[rounded corners](6.5,-0.7)--(6.5,-1.5);
\draw(5.9,-1.5)--(6.5,-1.5)--(6.5,-2.2);
\begin{scope}[shift={(-0.7,0.5)}]
\draw(7.5,0.5) .. controls (6.7,0.2) .. (7.1,-0.95);
\filldraw [color=white] (7.5,0.5) circle (12pt);
\draw (7.5,0.5) circle (12pt);
\node at (7.5,0.5)[scale=0.7,rotate=45,xscale=1,yscale=-1]{$x_1$};
\end{scope}
\begin{scope}[shift={(0.7,-0.9)}]
\filldraw [color=white] (4.5,3.5) circle (12pt);
\draw (4.5,3.5) circle (12pt);
\node at (4.5,3.5)[scale=0.5,rotate=45,xscale=1,yscale=-1]{$x_k-1$};
\end{scope}
\draw(4.7,-3.6)--(4,-3.6)--(4,-3);
\filldraw[color=white] (4.1,-3.7) circle (12pt);
\draw (4.1,-3.7) circle (12pt);
\node at (4.1,-3.7)[scale=0.7,rotate=-45]{$x_{k-1}$};
\filldraw[color=white] (6.5,-1.5) circle (12pt);
\draw (6.5,-1.5) circle (12pt);
\node at (6.5,-1.5)[rotate=-45,scale=0.7]{$x_2$};
\node at (6,2)[rotate=45]{$\cdots$};
\node at (5,-3)[rotate=-45]{$\cdots$};

\begin{scope}[shift={(0,-1)}]
\begin{scope}[xscale=1,yscale=-1]

\draw[rounded corners](3.8,4.4)--(2,6)--(-1.5,2.5);
\draw[rounded corners](-1.5,1.5)--(-0.5,0.5)--(0.5,1.5);
\draw(-0.4,0.4)--(0,0);

\draw[rounded corners](2.8,5)--(2.5,5)--(2.5,5.5);

\draw[rounded corners](2,6.1)--(2,7.2)--(3,6.3);
\end{scope}

\draw[rounded corners](-0.5,-0.7)--(0.5,-3)--(2.7,-5.1);
\draw[rounded corners](2.4,-5.1)--(1.6,-5.1)--(2,-5.8);
\filldraw [color=white] (2.5,-5.3) circle (12pt);
\draw (2.5,-5.3) circle (12pt);
\node at (2.5,-5.3)[scale=0.6,rotate=-45]{$x_{k+1}$};
\end{scope}
\end{tikzpicture}
$$

When connected, these become

$$
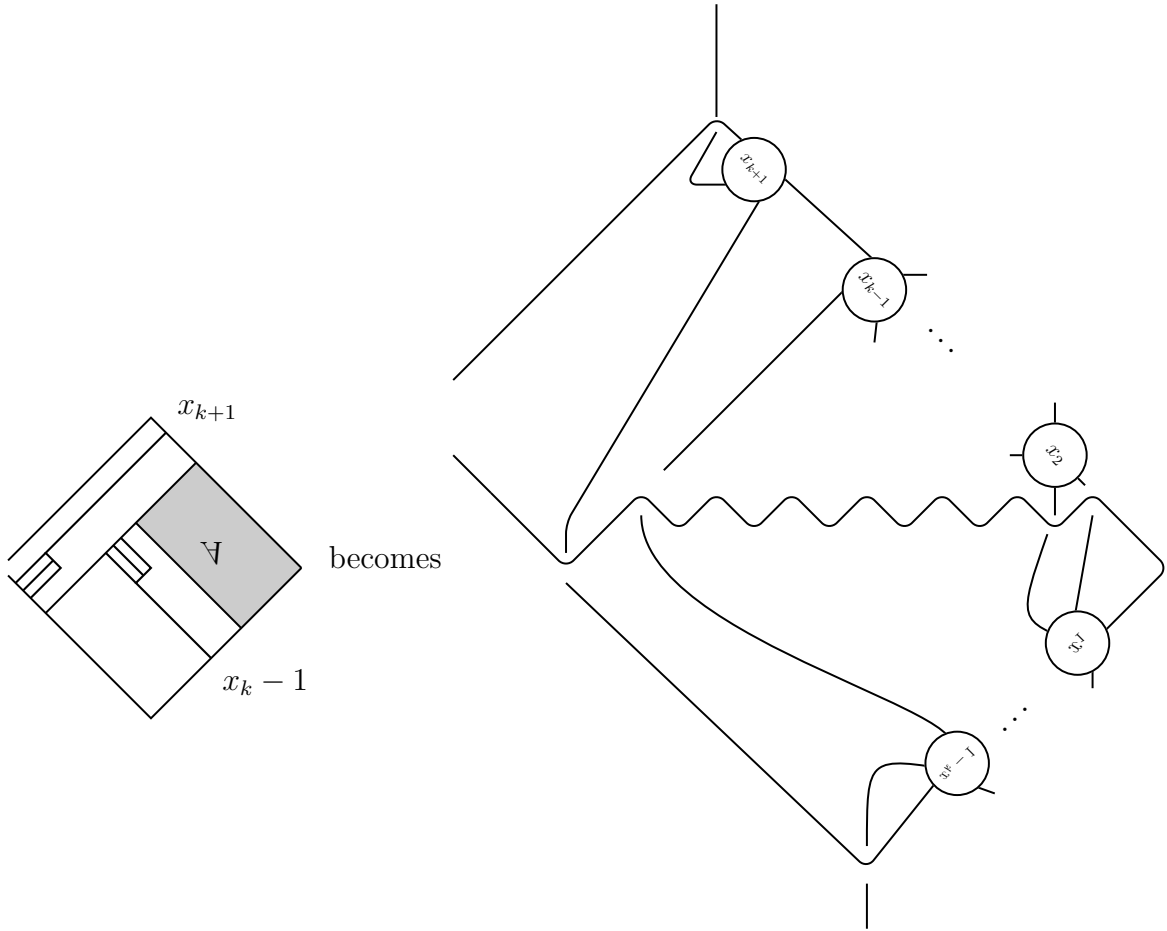
\begin{tikzpicture}[baseline=0, scale=0.2,yscale=-1]

    \chair[thick](-6,6)(5:2);
    \uchair[thick](-9,-9)(8:2);

    \filldraw[black!20!white] (0,0)--(-7,-7)--(-11,-3)--(-4,4)--(0,0);
    \node[yscale=-1] at (-6,-1) {$A$};
	 \draw (0,0)--(-10,10)--(-19.5,0.5);
	 \draw (0,0)--(-10,-10)--(-19.5,-0.5);
	 \draw (-11,-3)--(-4,4);
	 \draw (-7,-7)--(-17,3);

    \node at (-6,6) [below right] {$x_k-1$};
    \node at (-9,-9) [above right] {$x_{k+1}$};

\end{tikzpicture}
\text{ becomes }
\begin{tikzpicture}[yscale=-1,baseline=0]

\draw[rounded corners](4,4.2)--(4,4.8);
\draw(2,-7.5)--(2,-6);
\draw(5.2,2.7) .. controls (4,2.5) .. (4,3.7);
\draw(6.7,1)--(7,-0.7);

\draw(1,-0.7) .. controls (1,1)  and (6,2) .. (5,2.5);
\draw[rounded corners](5.7,3)--(5,2.75)--(4,4)--(0,0.2);

\draw[rounded corners](-1,-1)--(0,0)--(1,-1)--(1.5,-0.5)--(2,-1)--(2.5,-0.5)--(3,-1)--(3.5,-0.5)--(4,-1)--(4.5,-0.5)--(5,-1)--(5.5,-0.5)--(6,-1)--(6.5,-0.5)--(7,-1)--(8,0)--(7,1)--(7,1.6);
\draw[rounded corners](1.3,-1.3)--(4,-4);
\draw(6.9,-1.1)--(6.5,-1.5);
\draw[rounded corners](6.5,-0.7)--(6.5,-1.5);
\draw(5.9,-1.5)--(6.5,-1.5)--(6.5,-2.2);
\begin{scope}[shift={(-0.7,0.5)}]
\draw(7.5,0.5) .. controls (6.7,0.2) .. (7.1,-0.95);
\filldraw [color=white] (7.5,0.5) circle (12pt);
\draw (7.5,0.5) circle (12pt);
\node at (7.5,0.5)[scale=0.7,rotate=45,xscale=1,yscale=-1]{$x_1$};
\end{scope}
\begin{scope}[shift={(0.7,-0.9)}]
\filldraw [color=white] (4.5,3.5) circle (12pt);
\draw (4.5,3.5) circle (12pt);
\node at (4.5,3.5)[scale=0.5,rotate=45,xscale=1,yscale=-1]{$x_k-1$};
\end{scope}

\filldraw[color=white] (6.5,-1.5) circle (12pt);
\draw (6.5,-1.5) circle (12pt);
\node at (6.5,-1.5)[rotate=-45,scale=0.7]{$x_2$};
\node at (6,2)[rotate=45]{$\cdots$};
\node at (5,-3)[rotate=-45]{$\cdots$};

\begin{scope}[xscale=1,yscale=-1]

\draw[rounded corners](4.2,4)--(2,6)--(-1.5,2.5);
\draw[rounded corners](-1.5,1.5)--(-1,1);
\draw[rounded corners](2.5,5)--(2.5,5.6);

\end{scope}

\draw(4.1,-3)--(4.2,-3.9)--(4.8,-3.9);
\filldraw[color=white] (4.1,-3.7) circle (12pt);
\draw (4.1,-3.7) circle (12pt);
\node at (4.1,-3.7)[scale=0.7,rotate=-45]{$x_{k-1}$};

\draw[rounded corners](0,-0.2)--(0,-0.6)--(2.7,-5.1);
\draw[rounded corners](2.4,-5.1)--(1.6,-5.1)--(2,-5.8);
\filldraw [color=white] (2.5,-5.3) circle (12pt);
\draw (2.5,-5.3) circle (12pt);
\node at (2.5,-5.3)[scale=0.6,rotate=-45]{$x_{k+1}$};

\end{tikzpicture}
$$

This is exactly $\tilde{K}_l(x_1,x_2,\cdots, x_{k-1},x_k,x_{k+1})$.

\end{proof}

\begin{cor}
It is obvious that the closure of $\tilde{K}_l$ is $K$, so \\ $\psi'(T(x_1,x_2,\cdots,x_{n-1},x_n,x_{n+1})) = K(x_1,x_2,\cdots,x_{n-1},x_n,x_{n+1})$.

\end{cor}

By performing an ambient isotopy on the strand in $K(x_1,x_2,\cdots, x_{n-1},x_n,x_{n+1})$ that goes over the middle of the diagram from one end to the other, we see that we have the knot which is the closure of the product of tangles $\prod_{i=1}^{n}{x_i}$.

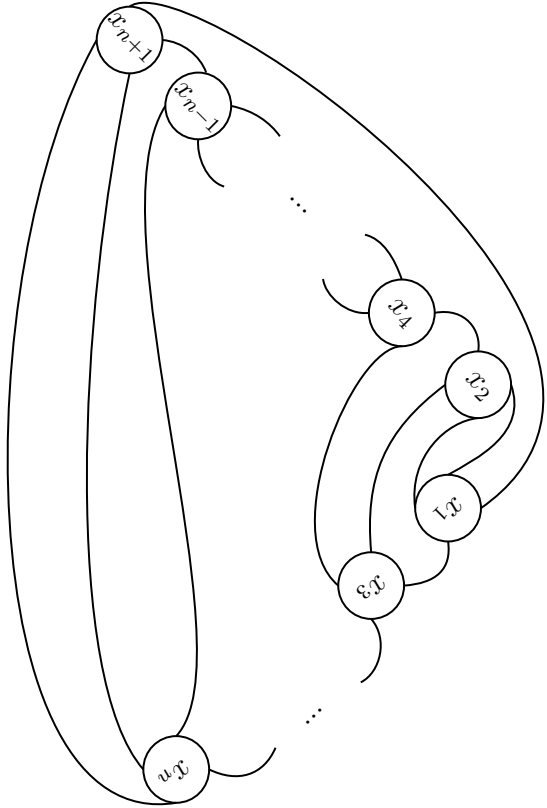
\tikzset{every picture/.style={line width=0.75pt}} 

\begin{tikzpicture}[x=0.75pt,y=0.75pt,yscale=-1,xscale=1]

\draw   (299.22,27.79) .. controls (299.22,18.5) and (306.65,10.96) .. (315.81,10.96) .. controls (324.97,10.96) and (332.4,18.5) .. (332.4,27.79) .. controls (332.4,37.08) and (324.97,44.61) .. (315.81,44.61) .. controls (306.65,44.61) and (299.22,37.08) .. (299.22,27.79) -- cycle ;
\draw    (508.2,201.79) .. controls (516.77,226.44) and (491.03,238.99) .. (476.43,247.24) ;
\draw    (459.84,264.07) .. controls (456.33,238.1) and (475.84,220.1) .. (491.61,218.62) ;
\draw    (437.64,286.34) .. controls (434.4,249.01) and (443.87,225.16) .. (475.02,201.79) ;
\draw    (453.09,182.34) .. controls (428.4,182.61) and (389.54,276.47) .. (421.05,303.16) ;
\draw    (469.68,165.51) .. controls (485.4,162.61) and (493.4,174.61) .. (491.61,184.97) ;
\draw    (333.84,61.07) .. controls (296.11,120.57) and (377.39,337.98) .. (339.33,379.14) ;
\draw    (332.4,27.79) .. controls (342.4,28.61) and (351.4,35.61) .. (354.43,44.24) ;
\draw    (322.74,395.96) .. controls (282.72,350.17) and (289.14,168.36) .. (315.81,44.61) ;
\draw    (339.33,412.79) .. controls (234.4,430.01) and (233.41,144.18) .. (299.22,27.79) ;
\draw    (454.23,303.16) .. controls (470.4,302.01) and (478.4,295.01) .. (476.43,280.89) ;
\draw    (315.81,10.96) .. controls (346.96,-12.4) and (614.51,179.06) .. (493.02,264.07) ;
\draw    (350.43,77.89) .. controls (349.5,86.15) and (354.5,99.15) .. (363.5,101.75) ;
\draw    (367.02,61.07) .. controls (378.61,63.41) and (386.61,69.41) .. (391.61,76.41) ;
\draw    (413.4,148.41) .. controls (414.4,156.81) and (425.4,166.81) .. (436.5,165.51) ;
\draw    (434.5,126.07) .. controls (443.5,128.07) and (449.5,138.07) .. (453.09,148.69) ;
\draw    (355.92,395.96) .. controls (370.4,403.01) and (383.4,399.01) .. (389.4,385.01) ;
\draw    (437.64,319.99) .. controls (446.4,330.01) and (442.4,347.01) .. (432.4,352.01) ;
\draw   (333.84,61.07) .. controls (333.84,51.78) and (341.27,44.24) .. (350.43,44.24) .. controls (359.59,44.24) and (367.02,51.78) .. (367.02,61.07) .. controls (367.02,70.36) and (359.59,77.89) .. (350.43,77.89) .. controls (341.27,77.89) and (333.84,70.36) .. (333.84,61.07) -- cycle ;
\draw   (436.5,165.51) .. controls (436.5,156.22) and (443.93,148.69) .. (453.09,148.69) .. controls (462.26,148.69) and (469.68,156.22) .. (469.68,165.51) .. controls (469.68,174.81) and (462.26,182.34) .. (453.09,182.34) .. controls (443.93,182.34) and (436.5,174.81) .. (436.5,165.51) -- cycle ;
\draw   (475.02,201.79) .. controls (475.02,192.5) and (482.45,184.97) .. (491.61,184.97) .. controls (500.77,184.97) and (508.2,192.5) .. (508.2,201.79) .. controls (508.2,211.09) and (500.77,218.62) .. (491.61,218.62) .. controls (482.45,218.62) and (475.02,211.09) .. (475.02,201.79) -- cycle ;
\draw   (459.84,264.07) .. controls (459.84,254.78) and (467.27,247.24) .. (476.43,247.24) .. controls (485.59,247.24) and (493.02,254.78) .. (493.02,264.07) .. controls (493.02,273.36) and (485.59,280.89) .. (476.43,280.89) .. controls (467.27,280.89) and (459.84,273.36) .. (459.84,264.07) -- cycle ;
\draw   (421.05,303.16) .. controls (421.05,293.87) and (428.47,286.34) .. (437.64,286.34) .. controls (446.8,286.34) and (454.23,293.87) .. (454.23,303.16) .. controls (454.23,312.45) and (446.8,319.99) .. (437.64,319.99) .. controls (428.47,319.99) and (421.05,312.45) .. (421.05,303.16) -- cycle ;
\draw   (322.74,395.96) .. controls (322.74,386.67) and (330.17,379.14) .. (339.33,379.14) .. controls (348.49,379.14) and (355.92,386.67) .. (355.92,395.96) .. controls (355.92,405.25) and (348.49,412.79) .. (339.33,412.79) .. controls (330.17,412.79) and (322.74,405.25) .. (322.74,395.96) -- cycle ;

\draw (476.2,265.08) node [inner sep=0.75pt]  [rotate=-134.97]  {$x_{1}$};
\draw (489.62,191.79) node [anchor=north west][inner sep=0.75pt]  [rotate=-45]  {$x_{2}$};
\draw (448.41,302.17) node [anchor=north west][inner sep=0.75pt]  [rotate=-135]  {$x_{3}$};
\draw (450.5,155) node [anchor=north west][inner sep=0.75pt]  [rotate=-45]  {$x_{4}$};
\draw (349.39,395.38) node [anchor=north west][inner sep=0.75pt]  [rotate=-135]  {$x_{n}$};
\draw (342.44,46.06) node [anchor=north west][inner sep=0.75pt]  [font=\small,rotate=-45]  {$x_{n-1}$};
\draw (308.82,10.78) node [anchor=north west][inner sep=0.75pt]  [font=\small,rotate=-45]  {$x_{n+1}$};
\draw (401,111) node [inner sep=0.75pt]  [rotate=-45]  {$...$};
\draw (416,365) node [anchor=north west][inner sep=0.75pt]  [rotate=-135]  {$...$};

\end{tikzpicture}

\section{Concatenation of Tangles}

\begin{defn}[Flipped Reflected Chair]
This is a reflected chair rotated by 180 degrees.  It is left-oriented, so $n$ of them will produce an opposite $n$-crossing-tangle when acted upon by $\psi'$.
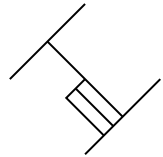
$$
\begin{tikzpicture}

\draw(0,2)--(1,3);
\draw(1,1)--(2,2);
\draw(1.5,1.5)--(0.5,2.5);
\draw(1,2)--(0.75,1.75)--(1.25,1.25);
\draw(0.875,1.875)--(1.375,1.375);

\end{tikzpicture}
$$
\end{defn}

The building blocks of our diagrams relating concatenations of tangles are called `stories':

\begin{defn}\label{story}
A \emph{story} of a concatenation Thompson diagram has `furniture' consisting of a single pillar to the left of $x_i-2$ flipped reflected chairs.  The top left edge of the diagram is the ceiling, and the bottom right edge is the floor.

$$
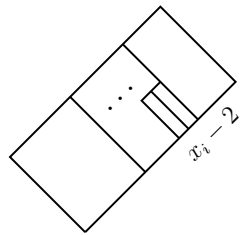
\begin{tikzpicture}
\draw(0,0)--(-1,1)--(1,3)--(2,2)--(0,0);
\begin{scope}[shift={(0.5,0.5)}]
\draw(0.3,0.3)--(-0.7,1.3);
\end{scope}
\draw(1.5,1.5)--(0.5,2.5);
\draw(1,2)--(0.75,1.75)--(1.25,1.25);
\draw(0.875,1.875)--(1.375,1.375);
\node at (1.7,1.3)[rotate=45,scale=0.8]{$x_i-2$};
\node at (0.5,1.8)[rotate=45]{$\cdots$};
\end{tikzpicture}
$$

\end{defn}

\begin{defn}
The Thompson diagram $U(x_1,x_2,\cdots,x_{k-1},x_k)$ is an apartment building\footnote{or perhaps the leaning tower of Pisa} made out of stories as defined above.  We define $U(x_1,x_2,\cdots,x_{k-1},x_k)$ recursively:
 
$$
U(x_1)=
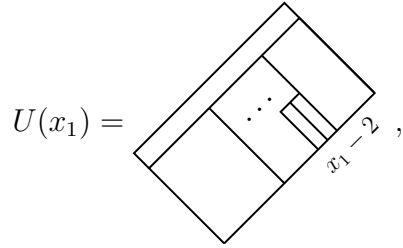
\begin{tikzpicture}[baseline=1.5cm]
\draw(-1,1)--(-1.2,1.2)--(0.8,3.2)--(2,2);
\draw(0,0)--(-1,1)--(1,3)--(2,2)--(0,0);
\begin{scope}[shift={(0.5,0.5)}]
\draw(0.3,0.3)--(-0.7,1.3);
\end{scope}
\draw(1.5,1.5)--(0.5,2.5);
\draw(1,2)--(0.75,1.75)--(1.25,1.25);
\draw(0.875,1.875)--(1.375,1.375);
\node at (1.7,1.3)[rotate=45,scale=0.8]{$x_1-2$};
\node at (0.5,1.8)[rotate=45]{$\cdots$};
\end{tikzpicture},
$$

and $U(x_1,x_2,\cdots,x_{n-1},x_n,x_{n+1})$ is created from $U(x_1,x_2,\cdots,x_{n-1},x_n)$ by appending a bottom story with $x_{n+1}-2$ chairs, so that the pillar of the bottom story is to right of all the chairs of the story above it (as illustrated by the red dot and line):

$U(x_1,x_2,\cdots,x_{n-1},x_n,x_{n+1})=
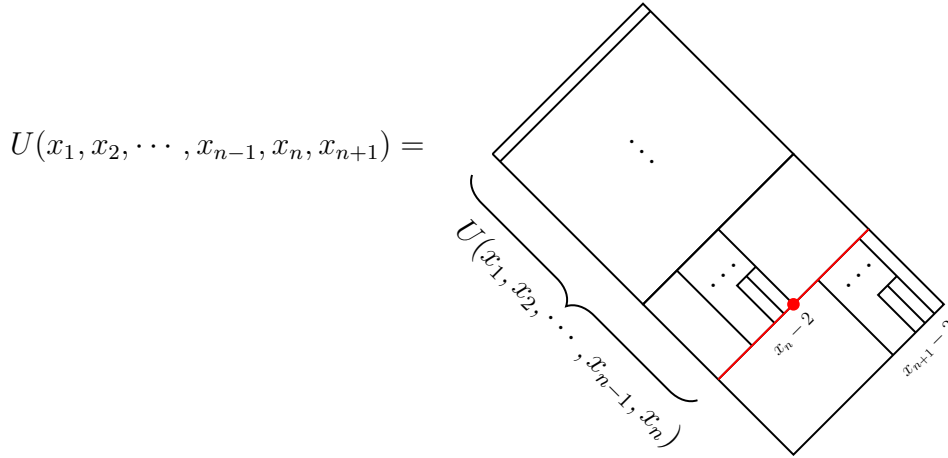
\begin{tikzpicture}[baseline=0cm]
\draw(0,0)--(2,2)--(4,0)--(2,-2)--(0,0);
\draw (0.1,-0.1)--(2.1,1.9);
\draw[decorate,decoration= {brace,amplitude=10pt},xshift=-5pt,yshift=-5pt]  (2.9,-3.1)--(-0.1,-0.1);
\node at (1,-2.4)[ rotate=-45] {$U(x_1,x_2,\cdots,x_{n-1},x_n) $};

\draw(3,-3)--(4,-4)--(6,-2)--(5,-1);

\begin{scope}[shift={(3,-3)}]
\draw(0,0)--(-1,1)--(1,3)--(2,2)--(0,0);

\draw(0.45,0.45)--(-0.55,1.45);

\begin{scope}[shift={(-0.5,-0.5)}]
\draw(1.5,1.5)--(0.5,2.5);
\draw(1,2)--(0.75,1.75)--(1.25,1.25);
\draw(0.875,1.875)--(1.375,1.375);
\end{scope}

\begin{scope}[shift={(1.625,-0.375)}]

\begin{scope}[shift={(0.25,0.25)}]

\draw(0.45,0.45)--(-0.55,1.45);

\end{scope}
\begin{scope}[shift={(-0.25,-0.25)}]
\draw(1.5,1.5)--(0.5,2.5);
\draw(1,2)--(0.75,1.75)--(1.25,1.25);
\draw(0.875,1.875)--(1.375,1.375);
\end{scope}

\end{scope}

\end{scope}

\node at (4,-2.4)[rotate=45,scale=0.6]{$x_{n}-2$};
\node at (2,0)[rotate=-45]{$\cdots$};
\draw[color=red](3,-3)--(5,-1);
\filldraw[color=red] (4,-2) circle (2pt);
\node at (3.1,-1.6)[rotate=45]{$\cdots$};
\node at (4.9,-1.6)[rotate=45]{$\cdots$};
\node at (5.8,-2.6)[rotate=45,scale=0.6]{$x_{n+1}-2$};
\end{tikzpicture}
$

\end{defn}

\begin{defn}  
The knot diagram $J(x_1,x_2,\cdots,x_{n-1},x_n)$ is defined as the diagram below left, and the 6-tangle diagram $\tilde{J}_{rb}(x_1,x_2,\cdots,x_{n-1},x_n)$ is defined as the diagram below right:

$$
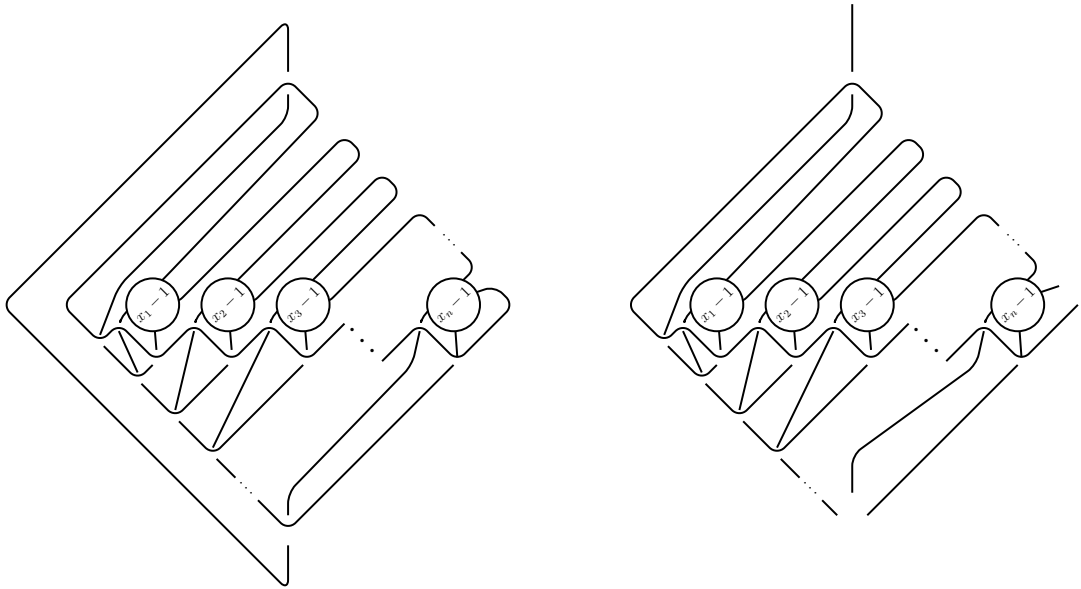
\begin{tikzpicture}[xscale=-1,yscale=1]

\draw[rounded corners](5.4,-2.6)--(5,-3)--(2.8,-0.8);

\draw[rounded corners](7,0)--(6.65,0.35)--(4.55,2.55)--(5,3)--(8,0)--(7.5,-0.5)--(7.25,-0.25)--(6.75,-0.75)--(6.25,-0.25)--(5.75,-0.75)--(5.25,-0.25)--(4.75,-0.75)--(4.25,-0.25);
\draw[rounded corners](3.75,-0.75)--(3.25,-0.25)--(2.75,-0.75)--(2,0)--(2.25,0.25)--(3,0); 

\node at (4,-0.5)[rotate=-45]{$\cdots$};

\draw[rounded corners](7,0)--(6.5,0)--(4.25,2.25)--(4,2)--(5.65,0.35)--(6,0);

\draw[rounded corners](6,0)--(5.5,0)--(3.75,1.75)--(3.5,1.5)--(4.65,0.35)--(5,0);

\draw[rounded corners](7.45,-0.55)--(7,-1)--(6.8,-0.8);

\draw[rounded corners](6.95,-1.05)--(6.5,-1.5)--(5.8,-0.8);
\draw[rounded corners](6.45,-1.55)--(6,-2)--(4.8,-0.8);
\draw(5.95,-2.05)--(5.7,-2.3);

\draw[rounded corners](7.5,-0.4)--(7.25,0.25)--(5,2.5)--(5,2.8);

\draw[rounded corners](3.05,1.05)--(3.25,1.25)--(4.5,0)--(5,0);

\draw[rounded corners](3.25,-0.35)--(3.35,-0.8)--(5,-2.5)--(5,-2.8);
\draw[rounded corners](5,3.1)--(5,3.8)--(8.8,0)--(5,-3.8)--(5,-3.2);

\begin{scope}[shift={(2,0)}]

\draw[rounded corners](1,0)--(0.5,0.5)--(0.7,0.7);
\draw[rounded corners](0.8,0)--(0.75,-0.7);
\draw[rounded corners](1,0)--(1.2,-0.10)--(1.225,-0.2);

\end{scope}

\filldraw[fill=white](2.8,0) circle (10pt);
\node at (2.8,0)[scale=0.5,rotate=45]{$x_n-1$};

\begin{scope}[shift={(4,0)}]

\draw(0.8,0)--(0.75,-0.6);
\draw[rounded corners](1,0)--(1.2,-0.10)--(1.225,-0.2);
\draw(1.25,-0.35)--(2,-1.9);
\filldraw[fill=white](0.8,0) circle (10pt);
\node at (0.8,0)[scale=0.5,rotate=45]{$x_3-1$};

\end{scope}

\begin{scope}[shift={(5,0)}]

\draw(0.8,0)--(0.75,-0.6);
\draw[rounded corners](1,0)--(1.2,-0.10)--(1.225,-0.2);
\draw(1.25,-0.35)--(1.5,-1.4);
\filldraw[fill=white](0.8,0) circle (10pt);
\node at (0.8,0)[scale=0.5,rotate=45]{$x_2-1$};
\end{scope}

\begin{scope}[shift={(6,0)}]

\draw(0.8,0)--(0.75,-0.6);
\draw[rounded corners](1,0)--(1.2,-0.10)--(1.225,-0.2);
\draw(1.25,-0.35)--(1,-0.9);
\filldraw[fill=white](0.8,0) circle (10pt);
\node at (0.8,0)[scale=0.5,rotate=45]{$x_1-1$};
\end{scope}
\node at (5.55,-2.45)[scale=0.5,rotate=-45]{$\cdots$};

\node at (2.85,0.85)[scale=0.5,rotate=-45]{$\cdots$};

\begin{scope}[shift={(-7.5,0)}]
\draw(5.4,-2.6)--(5.2,-2.8);
\draw(4.8,-2.8)--(2.8,-0.8);

\draw[rounded corners](7,0)--(6.65,0.35)--(4.55,2.55)--(5,3)--(8,0)--(7.5,-0.5)--(7.25,-0.25)--(6.75,-0.75)--(6.25,-0.25)--(5.75,-0.75)--(5.25,-0.25)--(4.75,-0.75)--(4.25,-0.25);
\draw[rounded corners](3.75,-0.75)--(3.25,-0.25)--(2.75,-0.75)--(2,0);

\node at (4,-0.5)[rotate=-45]{$\cdots$};
\draw(2.25,0.25)--(3,0); 

\draw[rounded corners](7,0)--(6.5,0)--(4.25,2.25)--(4,2)--(5.65,0.35)--(6,0);

\draw[rounded corners](6,0)--(5.5,0)--(3.75,1.75)--(3.5,1.5)--(4.65,0.35)--(5,0);

\draw[rounded corners](7.45,-0.55)--(7,-1)--(6.8,-0.8);

\draw[rounded corners](6.95,-1.05)--(6.5,-1.5)--(5.8,-0.8);
\draw[rounded corners](6.45,-1.55)--(6,-2)--(4.8,-0.8);
\draw(5.95,-2.05)--(5.7,-2.3);

\draw[rounded corners](7.5,-0.4)--(7.25,0.25)--(5,2.5)--(5,2.8);

\draw[rounded corners](3.05,1.05)--(3.25,1.25)--(4.5,0)--(5,0);

\draw[rounded corners](3.25,-0.35)--(3.35,-0.8)--(5,-2)--(5,-2.5);
\draw(5,3.1)--(5,4);

\begin{scope}[shift={(2,0)}]

\draw[rounded corners](1,0)--(0.5,0.5)--(0.7,0.7);
\draw[rounded corners](0.8,0)--(0.75,-0.7);
\draw[rounded corners](1,0)--(1.2,-0.10)--(1.225,-0.2);

\end{scope}

\filldraw[fill=white](2.8,0) circle (10pt);
\node at (2.8,0)[scale=0.5,rotate=45]{$x_n-1$};

\begin{scope}[shift={(4,0)}]

\draw(0.8,0)--(0.75,-0.6);
\draw[rounded corners](1,0)--(1.2,-0.10)--(1.225,-0.2);
\draw(1.25,-0.35)--(2,-1.9);
\filldraw[fill=white](0.8,0) circle (10pt);
\node at (0.8,0)[scale=0.5,rotate=45]{$x_3-1$};

\end{scope}

\begin{scope}[shift={(5,0)}]

\draw(0.8,0)--(0.75,-0.6);
\draw[rounded corners](1,0)--(1.2,-0.10)--(1.225,-0.2);
\draw(1.25,-0.35)--(1.5,-1.4);
\filldraw[fill=white](0.8,0) circle (10pt);
\node at (0.8,0)[scale=0.5,rotate=45]{$x_2-1$};
\end{scope}

\begin{scope}[shift={(6,0)}]

\draw(0.8,0)--(0.75,-0.6);
\draw[rounded corners](1,0)--(1.2,-0.10)--(1.225,-0.2);
\draw(1.25,-0.35)--(1,-0.9);
\filldraw[fill=white](0.8,0) circle (10pt);
\node at (0.8,0)[scale=0.5,rotate=45]{$x_1-1$};
\end{scope}
\node at (5.55,-2.45)[scale=0.5,rotate=-45]{$\cdots$};

\node at (2.85,0.85)[scale=0.5,rotate=-45]{$\cdots$};
\end{scope}
\end{tikzpicture}
$$
\end{defn}

\begin{thm}

$\psi'(\tilde{U}_{br}(x_1,x_2,\cdots,x_{k-1},x_k)) = \tilde{J}_{br}(x_1,x_2,\cdots,x_{n-1},x_n)$

\end{thm}

\begin{proof}

Base case: When $n=1$, we have

$$
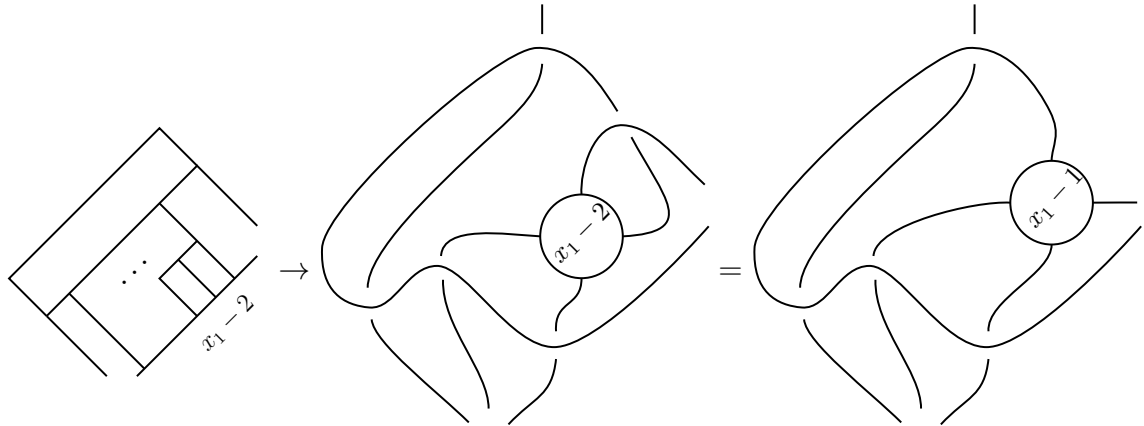
\begin{tikzpicture}[baseline=1.5cm]

\draw(0.2,0.2)--(1.8,1.8);
\draw(-1,1)--(-1.5,1.5)--(0.5,3.5)--(1,3);

\draw(-0.2,0.2)--(-1,1);
\draw(-1,1)--(1,3)--(1.8,2.2);
\draw(0.3,0.3)--(-0.7,1.3);
\draw(1.5,1.5)--(0.5,2.5);
\draw(1,2)--(0.5,1.5)--(1,1);
\draw(0.75,1.75)--(1.25,1.25);

\node at (1.4,0.9)[rotate=45,scale=0.8]{$x_1-2$};
\node at (0.2,1.6)[rotate=45]{$\cdots$};

\end{tikzpicture}
\rightarrow
\begin{tikzpicture}[x=0.75pt,y=0.75pt,yscale=-1,xscale=1, baseline=-250]

\draw    (147,322) .. controls (148,306) and (179.4,311.04) .. (197.02,312.08) ;
\draw   (197.02,312.08) .. controls (197.02,300.42) and (206.34,290.97) .. (217.84,290.97) .. controls (229.34,290.97) and (238.66,300.42) .. (238.66,312.08) .. controls (238.66,323.75) and (229.34,333.2) .. (217.84,333.2) .. controls (206.34,333.2) and (197.02,323.75) .. (197.02,312.08) -- cycle ;
\draw    (198,225) .. controls (197,242) and (180,254) .. (162,271) .. controls (144,288) and (110,320) .. (110,338) ;
\draw    (243,262) .. controls (259,292) and (278,310) .. (238.66,312.08) ;
\draw    (148,334) .. controls (147.4,358.04) and (171,380) .. (171,399) ;
\draw    (205,360) .. controls (204,344) and (218,348) .. (217.84,333.2) ;
\draw    (198,195) -- (198,210) ;
\draw    (181,407) .. controls (194,393) and (203,393) .. (205,374) ;
\draw    (282,307) .. controls (272.75,317.65) and (226,368) .. (204,368) .. controls (182,368) and (160,327) .. (145,327) .. controls (130,327) and (124,349) .. (111,348) .. controls (98,347) and (86,340) .. (87,318) .. controls (88,296) and (174,218) .. (195,217) .. controls (216,216) and (229,238) .. (236,249) ;
\draw    (217.84,290.97) .. controls (217,279) and (224,255) .. (239,256) .. controls (254,257) and (265,272) .. (280,286) ;
\draw    (112,354) .. controls (113,366) and (145.33,390.47) .. (161,406) ;

\draw (199,320) node [anchor=north west][inner sep=0.75pt]  [font=\footnotesize,rotate=-315]  {$x_{1} -2$};

\end{tikzpicture}
=
\begin{tikzpicture}[x=0.75pt,y=0.75pt,yscale=-1,xscale=1, baseline=-250]

\draw    (147,322) .. controls (148,306) and (198.4,294.04) .. (216.02,295.08) ;
\draw   (216.02,295.08) .. controls (216.02,283.42) and (225.34,273.97) .. (236.84,273.97) .. controls (248.34,273.97) and (257.66,283.42) .. (257.66,295.08) .. controls (257.66,306.75) and (248.34,316.2) .. (236.84,316.2) .. controls (225.34,316.2) and (216.02,306.75) .. (216.02,295.08) -- cycle ;
\draw    (198,225) .. controls (197,242) and (180,254) .. (162,271) .. controls (144,288) and (110,320) .. (110,338) ;
\draw    (280,295) .. controls (267,295) and (276,295) .. (257.66,295.08) ;
\draw    (148,334) .. controls (147.4,358.04) and (171,380) .. (171,399) ;
\draw    (205,360) .. controls (204,344) and (237,331) .. (236.84,316.2) ;
\draw    (198,195) -- (198,210) ;
\draw    (181,407) .. controls (194,393) and (203,393) .. (205,374) ;
\draw    (282,307) .. controls (272.75,317.65) and (226,368) .. (204,368) .. controls (182,368) and (160,327) .. (145,327) .. controls (130,327) and (124,349) .. (111,348) .. controls (98,347) and (86,340) .. (87,318) .. controls (88,296) and (174,218) .. (195,217) .. controls (216,216) and (230,237) .. (236,249) .. controls (242,261) and (237,266) .. (236.84,273.97) ;
\draw    (112,354) .. controls (113,366) and (145.33,390.47) .. (161,406) ;

\draw (220.22,302.86) node [anchor=north west][inner sep=0.75pt]  [font=\footnotesize,rotate=-312.33]  {$x_{1} -1$};

\end{tikzpicture}
$$

Now supposed we've proven the statement true for some natural number $k \geq 1$, so we've shown:

$$
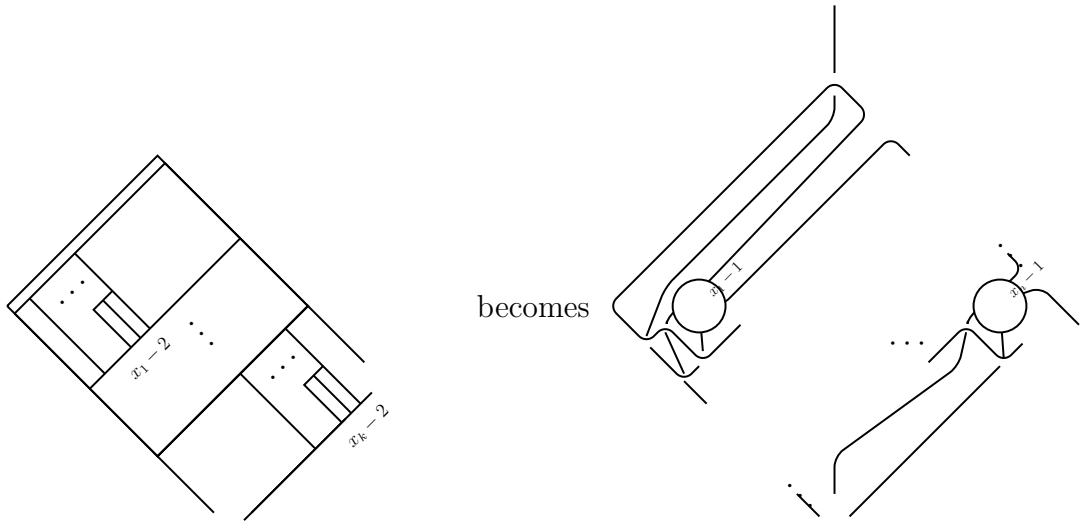
\begin{tikzpicture}
\draw(0,0)--(2,2)--(4,0)--(2,-2)--(0,0);

\begin{scope}[shift={(0.4,-1.8)}]
\draw(1.5,1.5)--(0.5,2.5);

\draw(1,2)--(0.75,1.75)--(1.25,1.25);
\draw(0.875,1.875)--(1.375,1.375);

\draw(0.7,0.7)--(2.7,2.7);
\draw(0.9,0.9)--(-0.1,1.9);
\node at (0.5,2)[rotate=45] {$\cdots$};
\end{scope}

\begin{scope}[shift={(3,-3)}]
\draw(-0.25,0.25)--(-1,1)--(1,3)--(1.75,2.25);
\draw(1.85,1.85)--(0.25,0.25);

\begin{scope}[shift={(0.8,0.8)}]
\draw(0.3,0.3)--(-0.7,1.3);

\end{scope}
\begin{scope}[shift={(0.2,0.2)}]
\draw(1.5,1.5)--(0.5,2.5);

\draw(1,2)--(0.75,1.75)--(1.25,1.25);
\draw(0.875,1.875)--(1.375,1.375);
\end{scope}

\end{scope}
\node at (1.9,-0.7)[rotate=45,scale=0.6]{$x_1-2$};
\node at (4.8,-1.6)[rotate=45,scale=0.6]{$x_{k}-2$};
\draw(0.1,-0.1)--(2.1,1.9)--(4,0)--(2,-2)--(0.1,-0.1);
\node at (2.6,-0.4)[rotate=-45]{$\cdots$};
\draw(3.15,-2.85)--(4.5,-1.5);

\node at (3.7,-0.8)[rotate=45]{$\cdots$};

\hspace{1cm};
\node at (6,0){becomes};

\begin{scope}[xscale=-1,yscale=1,shift={(-16,0)}]
\draw(7,-1)--(6.7,-1.3);
\draw(4.8,-2.8)--(2.8,-0.8);

\draw[rounded corners](7,0)--(6.65,0.35)--(4.55,2.55)--(5,3)--(8,0)--(7.5,-0.5)--(7.25,-0.25)--(6.75,-0.75)--(6.25,-0.25);

\draw[rounded corners](3.75,-0.75)--(3.25,-0.25)--(2.75,-0.75)--(2.5,-0.5);

\node at (5,-0.5){$\cdots$};
\draw[rounded corners](1.75,-0.25)--(2.25,0.25)--(3,0); 

\draw[rounded corners](7,0)--(6.5,0)--(4.25,2.25)--(4,2);

\draw[rounded corners](7.45,-0.55)--(7,-1)--(6.8,-0.8);

\draw(5.2,-2.8)--(5.5,-2.5);

\draw[rounded corners](7.5,-0.4)--(7.25,0.25)--(5,2.5)--(5,2.8);

\draw[rounded corners](3.25,-0.35)--(3.35,-0.8)--(5,-2)--(5,-2.5);
\draw(5,3.1)--(5,4);

\begin{scope}[shift={(2,0)}]

\draw[rounded corners](1,0)--(0.5,0.5)--(0.7,0.7);
\draw[rounded corners](0.8,0)--(0.75,-0.7);
\draw[rounded corners](1,0)--(1.2,-0.10)--(1.225,-0.2);

\end{scope}

\filldraw[fill=white](2.8,0) circle (10pt);
\node at (2.8,0)[scale=0.5,rotate=45]{$x_n-1$};

\begin{scope}[shift={(6,0)}]

\draw(0.8,0)--(0.75,-0.6);
\draw[rounded corners](1,0)--(1.2,-0.10)--(1.225,-0.2);
\draw(1.25,-0.35)--(1,-0.9);
\filldraw[fill=white](0.8,0) circle (10pt);
\node at (0.8,0)[scale=0.5,rotate=45]{$x_1-1$};
\end{scope}
\node at (6.15,-1.85)[rotate=-45]{$\cdots$};

\node at (3.35,1.35)[rotate=-45]{$\cdots$};
\end{scope}
\end{tikzpicture}
$$

(The above right is an abbreviated version of $\tilde{J}_{br}(x_1,x_2,\cdots,x_{k-1},x_k$.))

Bringing in the $k+1$st chair set gives:

$$
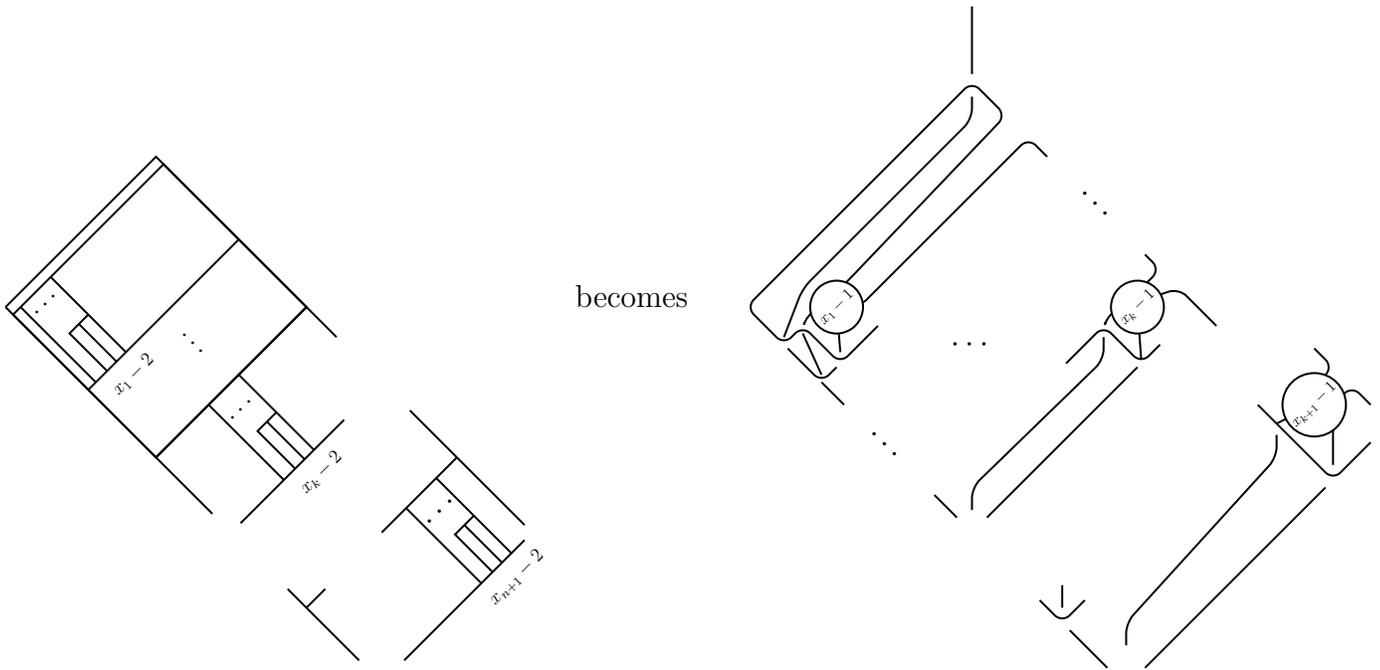
\begin{tikzpicture}[baseline]
\draw(0,0)--(2,2)--(4,0)--(2,-2)--(0,0);
\begin{scope}[shift={(1.1,-1.1)}]
\draw(0,0)--(2,2);

\begin{scope}[shift={(-0.2,-0.2)}]
\draw(0.3,0.3)--(-0.7,1.3);

\end{scope}
\begin{scope}[shift={(-1,-1)}]
\draw(1.5,1.5)--(0.5,2.5);

\draw(1,2)--(0.75,1.75)--(1.25,1.25);
\draw(0.875,1.875)--(1.375,1.375);
\end{scope}
\end{scope}
\node at (1.7,-0.9)[rotate=45,scale=0.6]{$x_1-2$};
\node at (0.55,0.05)[scale=0.8,rotate=45]{$\cdots$};

\begin{scope}[shift={(3,-3)}]
\draw(-0.25,0.25)--(-1,1)--(1,3)--(1.4,2.6);

\begin{scope}[shift={(0.25,0.25)}]
\draw(0.45,0.45)--(-0.55,1.45);
\end{scope}

\begin{scope}[shift={(-0.4,-0.4)}]
\draw(1.5,1.5)--(0.5,2.5);
\draw(1,2)--(0.75,1.75)--(1.25,1.25);
\draw(0.875,1.875)--(1.375,1.375);
\end{scope}
\end{scope}
\node at (4.2,-2.2)[rotate=45,scale=0.6]{$x_{k}-2$};
\draw(0.1,-0.1)--(2.1,1.9)--(4,0)--(2,-2)--(0.1,-0.1);
\node at (2.5,-0.5)[scale=0.8,rotate=-45]{$\cdots$};
\draw(3.125,-2.875)--(4.5,-1.5);

\node at (3.15,-1.35)[scale=0.8,rotate=45]{$\cdots$};

\begin{scope}[xshift=1cm,yshift=-1cm]
\draw(3,-3)--(3.7,-3.7);
\draw(4.3,-3.7)--(5.9,-2.1);
\draw(5.9,-1.9)--(5,-1);
\begin{scope}[shift={(3,-3)}]

\begin{scope}[shift={(1.625,-0.375)}]

\begin{scope}[shift={(0.4,0.4)}]
\draw(0.3,0.3)--(-0.7,1.3);
\end{scope}
\begin{scope}[shift={(-0.4,-0.4)}]
\draw(1.5,1.5)--(0.5,2.5);
\draw(1,2)--(0.75,1.75)--(1.25,1.25);
\draw(0.875,1.875)--(1.375,1.375);
\end{scope}

\end{scope}

\end{scope}

\draw(3,-3)--(3.25,-2.75);
\draw(4,-2)--(5,-1);
\draw(3,-3)--(2.75,-2.75);
\draw(5,-1)--(4.375,-0.375);

\node at (4.8,-1.7)[rotate=45]{$\cdots$};
\node at (5.8,-2.6)[rotate=45,scale=0.6]{$x_{n+1}-2$};
\end{scope}
\end{tikzpicture}
\;
\text{ becomes }
\hspace{0.5cm}
\;
\begin{tikzpicture}[baseline,xscale=-1,yscale=1]

\draw(7,-1)--(6.7,-1.3);
\draw(4.8,-2.8)--(2.8,-0.8);

\draw[rounded corners](7,0)--(6.65,0.35)--(4.55,2.55)--(5,3)--(8,0)--(7.5,-0.5)--(7.25,-0.25)--(6.75,-0.75)--(6.25,-0.25);

\draw[rounded corners](3.75,-0.75)--(3.25,-0.25)--(2.75,-0.75)--(2.5,-0.5);

\node at (5,-0.5){$\cdots$};
\draw[rounded corners](1.75,-0.25)--(2.25,0.25)--(3,0); 

\draw[rounded corners](7,0)--(6.5,0)--(4.25,2.25)--(4,2);

\draw[rounded corners](7.45,-0.55)--(7,-1)--(6.8,-0.8);

\draw(5.2,-2.8)--(5.5,-2.5);

\draw[rounded corners](7.5,-0.4)--(7.25,0.25)--(5,2.5)--(5,2.8);

\draw[rounded corners](3.25,-0.4)--(3.25,-0.7)--(5,-2.4)--(5,-2.7);
\draw(5,3.1)--(5,4);

\begin{scope}[shift={(2,0)}]

\draw[rounded corners](1,0)--(0.5,0.5)--(0.7,0.7);
\draw[rounded corners](0.8,0)--(0.75,-0.7);
\draw[rounded corners](1,0)--(1.2,-0.10)--(1.225,-0.2);

\end{scope}

\filldraw[fill=white](2.8,0) circle (10pt);
\node at (2.8,0)[scale=0.5,rotate=45]{$x_k-1$};

\begin{scope}[shift={(6,0)}]

\draw(0.8,0)--(0.75,-0.6);
\draw[rounded corners](1,0)--(1.2,-0.10)--(1.225,-0.2);
\draw(1.25,-0.35)--(1,-0.9);
\filldraw[fill=white](0.8,0) circle (10pt);
\node at (0.8,0)[scale=0.5,rotate=45]{$x_1-1$};
\end{scope}
\node at (6.15,-1.85)[rotate=-45]{$\cdots$};

\node at (3.35,1.35)[rotate=-45]{$\cdots$};


\begin{scope}[shift={(-1.3,-1.3)}]

\draw[rounded corners](1.75,-0.25)--(1.25,0.25)--(1,0);
\draw(1.6,-1.1)--(4,-3.5);
\draw(4.5,-3.5)--(5,-3);
\draw[rounded corners](4.25,-3.2)--(4.25,-2.9)--(2.25,-0.7)--(2.25,-0.4);
\draw[rounded corners](1.5,0)--(1.75,0)--(2.25,-0.25);
\draw(1.5,0)--(1.5,-0.8);
\draw[rounded corners](5.4,-2.6)--(5.1,-2.9)--(4.8,-2.6);
\draw(5.1,-2.7)--(5.1,-2.4);
\draw[rounded corners](2.5,-0)--(1.5,-1)--(1,-0.5);
\draw[rounded corners](2,0)--(1.5,0.5)--(1.75,0.75);

\filldraw[fill=white](1.75,0) circle (12pt);
\node at (1.75,0)[scale=0.5,rotate=45]{$x_{k+1}-1$};
\end{scope}

\end{tikzpicture}
$$

When connected, these become
 
$$
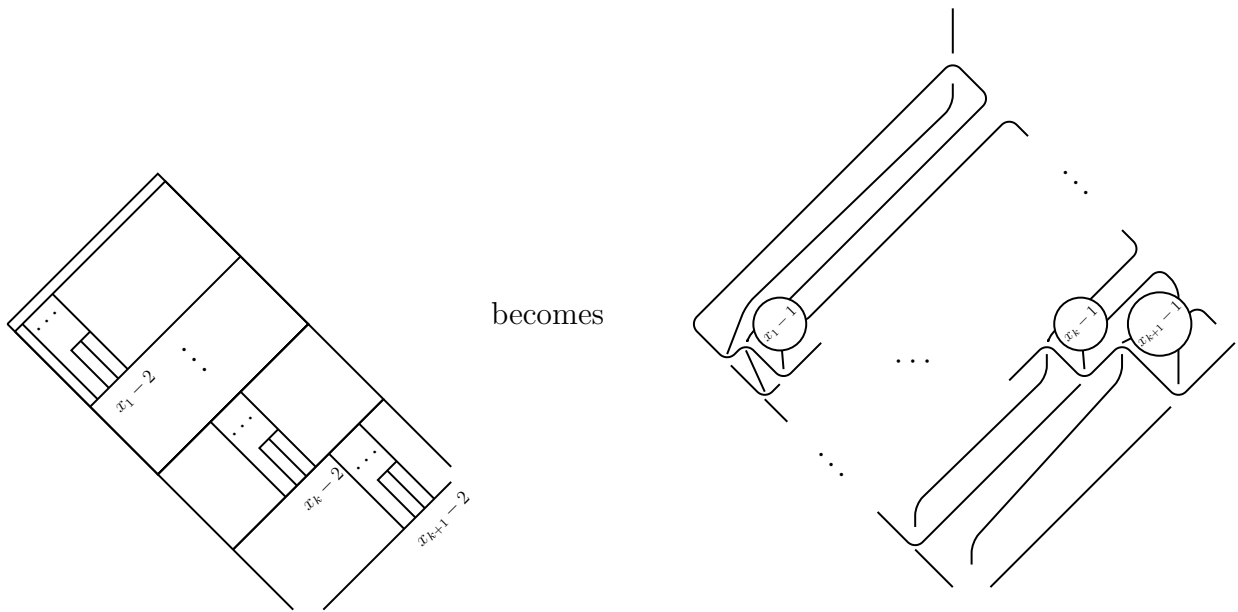
\begin{tikzpicture}[baseline]
\begin{scope}[shift={(9,0)}]

\draw(0,0)--(2,2)--(4,0)--(2,-2)--(0,0);
\draw(3,-3)--(3.8,-3.8);
\draw(4.2,-3.8)--(5.9,-2.1);
\draw(5.9,-1.9)--(5,-1);

\begin{scope}[shift={(1.1,-1.1)}]
\draw(0,0)--(2,2);

\begin{scope}[shift={(-0.2,-0.2)}]
\draw(0.3,0.3)--(-0.7,1.3);

\end{scope}
\begin{scope}[shift={(-1,-1)}]
\draw(1.5,1.5)--(0.5,2.5);

\draw(1,2)--(0.75,1.75)--(1.25,1.25);
\draw(0.875,1.875)--(1.375,1.375);
\end{scope}
\end{scope}

\node at (1.7,-0.9)[rotate=45,scale=0.6]{$x_1-2$};
\node at (0.55,0.05)[scale=0.8,rotate=45]{$\cdots$};

\begin{scope}[shift={(3,-3)}]
\draw(0,0)--(-1,1)--(1,3)--(2,2)--(0,0);
\draw(0.7,0.7)--(-0.3,1.7);

\begin{scope}[shift={(-0.4,-0.4)}]
\draw(1.5,1.5)--(0.5,2.5);
\draw(1,2)--(0.75,1.75)--(1.25,1.25);
\draw(0.875,1.875)--(1.375,1.375);
\end{scope}

\begin{scope}[shift={(1.625,-0.375)}]

\draw(0.65,0.65)--(-0.35,1.65);

\begin{scope}[shift={(-0.45,-0.45)}]
\draw(1.5,1.5)--(0.5,2.5);
\draw(1,2)--(0.75,1.75)--(1.25,1.25);
\draw(0.875,1.875)--(1.375,1.375);
\end{scope}

\end{scope}

\end{scope}
\node at (4.2,-2.2)[rotate=45,scale=0.6]{$x_{k}-2$};
\draw(0.1,-0.1)--(2.1,1.9)--(4,0)--(2,-2)--(0.1,-0.1);
\node at (2.5,-0.5)[rotate=-45]{$\cdots$};
\draw(3,-3)--(5,-1);

\node at (3.15,-1.35)[scale=0.8,rotate=45]{$\cdots$};
\node at (4.8,-1.8)[scale=0.8,rotate=45]{$\cdots$};
\node at (5.8,-2.6)[rotate=45,scale=0.6]{$x_{k+1}-2$};
\end{scope}
\end{tikzpicture}
\;
\text{becomes}
\hspace{1cm}
\;
\begin{tikzpicture}[baseline,xscale=-1,yscale=1]
\draw(7,-1)--(6.7,-1.3);

\draw[rounded corners](7,0)--(6.65,0.35)--(4,3)--(4.5,3.5)--(8,0)--(7.5,-0.5)--(7.25,-0.25)--(6.75,-0.75)--(6.25,-0.25);

\draw[rounded corners](3.75,-0.75)--(3.25,-0.25)--(2.75,-0.75)--(2.25,-0.25)--(1.5,-1)--(0.75,-0.25);

\node at (5,-0.5){$\cdots$};
\draw[rounded corners](1.5,0)--(1.5,0.5)--(1.75,0.75)--(2.75,-0.25); 

\draw[rounded corners](7,0)--(6.5,0)--(3.75,2.75)--(3.5,2.5);

\draw[rounded corners](7.45,-0.55)--(7,-1)--(6.8,-0.8);

\draw[rounded corners](2.8,-0.8)--(5,-3)--(5.5,-2.5);

\draw[rounded corners](7.5,-0.4)--(7.25,0.25)--(4.5,2.9)--(4.5,3.2);

\draw[rounded corners](3.25,-0.4)--(3.25,-0.7)--(5,-2.4)--(5,-2.7);
\draw(4.5,3.6)--(4.5,4.2);

\begin{scope}[shift={(2,0)}]

\draw[rounded corners](1,0)--(0,1)--(0.25,1.25);
\draw[rounded corners](0.8,0)--(0.75,-0.6);
\draw[rounded corners](1,0)--(1.2,-0.10)--(1.225,-0.2);

\end{scope}

\filldraw[fill=white](2.8,0) circle (10pt);
\node at (2.8,0)[scale=0.5,rotate=45]{$x_k-1$};

\begin{scope}[shift={(6,0)}]

\draw(0.8,0)--(0.75,-0.6);
\draw[rounded corners](1,0)--(1.2,-0.10)--(1.225,-0.2);
\draw(1.25,-0.35)--(1,-0.9);
\filldraw[fill=white](0.8,0) circle (10pt);
\node at (0.8,0)[scale=0.5,rotate=45]{$x_1-1$};
\end{scope}
\node at (6.1,-1.9)[rotate=-45]{$\cdots$};

\node at (2.85,1.85)[rotate=-45]{$\cdots$};

\draw[rounded corners](1.75,-0.25)--(1.25,0.25)--(1,0);
\draw(1.6,-1.1)--(4,-3.5);
\draw(4.5,-3.5)--(5,-3);
\draw[rounded corners](4.25,-3.2)--(4.25,-2.9)--(2.25,-0.7)--(2.25,-0.4);
\draw[rounded corners](1.5,0)--(1.75,0)--(2.25,-0.25);
\draw(1.5,0)--(1.5,-0.8);

\filldraw[fill=white](1.75,0) circle (12pt);
\node at (1.75,0)[scale=0.5,rotate=45]{$x_{k+1}-1$};

\end{tikzpicture}
$$

This is exactly $\tilde{J}_{br}(x_1,x_2,\cdots,x_{k-1},x_k,x_{k+1})$.
\end{proof}

It is obvious that the closure of $\tilde{J}_{br}(x_1,x_2,\cdots,x_{k-1},x_k,x_{k+1})$ is $J(x_1,x_2,\cdots,x_{k-1},x_k,x_{k+1})$.  Thus, $\psi'(U(x_1,x_2,\cdots,x_{n-1},x_n))=J(x_1,x_2,\cdots,x_{k-1},x_n,x_{n+1})$.

Notice a single strand zig-zagging over all other strands in the middle of the diagram.  Moving this strand over and merging adjacent crossings with the $x_{i-1}$-tanlges, we have the knot that is the closure of the tangle $x_1,x_2,...,x_n$ in Conway notation.

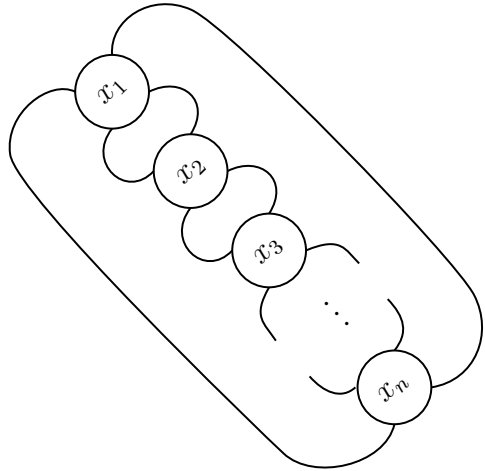
\tikzset{every picture/.style={line width=0.75pt}} 
$$
\begin{tikzpicture}[x=0.75pt,y=0.75pt,yscale=-1,xscale=1]

\draw   (172.2,136.55) .. controls (172.2,126.28) and (180.53,117.95) .. (190.8,117.95) .. controls (201.08,117.95) and (209.4,126.28) .. (209.4,136.55) .. controls (209.4,146.83) and (201.08,155.15) .. (190.8,155.15) .. controls (180.53,155.15) and (172.2,146.83) .. (172.2,136.55) -- cycle ;
\draw   (211.8,176.6) .. controls (211.8,166.33) and (220.13,158) .. (230.4,158) .. controls (240.67,158) and (249,166.33) .. (249,176.6) .. controls (249,186.87) and (240.67,195.2) .. (230.4,195.2) .. controls (220.13,195.2) and (211.8,186.87) .. (211.8,176.6) -- cycle ;
\draw   (251.4,216.65) .. controls (251.4,206.37) and (259.72,198.05) .. (270,198.05) .. controls (280.27,198.05) and (288.6,206.37) .. (288.6,216.65) .. controls (288.6,226.92) and (280.27,235.25) .. (270,235.25) .. controls (259.72,235.25) and (251.4,226.92) .. (251.4,216.65) -- cycle ;
\draw   (314.6,285.75) .. controls (314.6,275.48) and (322.93,267.15) .. (333.2,267.15) .. controls (343.48,267.15) and (351.8,275.48) .. (351.8,285.75) .. controls (351.8,296.03) and (343.48,304.35) .. (333.2,304.35) .. controls (322.93,304.35) and (314.6,296.03) .. (314.6,285.75) -- cycle ;
\draw    (190.8,117.95) .. controls (189.4,100.2) and (210.4,86.2) .. (230.61,95.11) .. controls (250.81,104.02) and (363.55,220.05) .. (373.4,239.2) .. controls (383.25,258.35) and (374.4,282.2) .. (351.8,285.75) ;
\draw    (172.2,136.55) .. controls (159.4,130.2) and (136.4,147.2) .. (139.4,168.2) .. controls (142.4,189.2) and (262.4,305.2) .. (277.4,319.2) .. controls (292.4,333.2) and (328.4,325.2) .. (333.2,304.35) ;
\draw    (211.8,176.6) .. controls (197.8,191.6) and (177.4,172.2) .. (190.8,155.15) ;
\draw    (251.4,216.65) .. controls (237.4,231.65) and (217,212.25) .. (230.4,195.2) ;
\draw    (273.4,262.2) .. controls (267.4,253.2) and (261.4,250.2) .. (270,235.25) ;
\draw    (230.4,158) .. controls (241.21,141.11) and (226.21,128.11) .. (209.4,136.55) ;
\draw    (270,198.05) .. controls (280.8,181.15) and (265.8,168.15) .. (249,176.6) ;
\draw    (315.4,223.2) .. controls (308.4,216.2) and (306.4,209.2) .. (288.6,216.65) ;
\draw    (333.2,267.15) .. controls (340.4,258.2) and (338.4,250.2) .. (329.8,241.2) ;
\draw    (313.5,285.75) .. controls (306.4,292.2) and (297.4,288.2) .. (290.4,280.2) ;

\draw (180.02,138.54) node [anchor=north west][inner sep=0.75pt]  [rotate=-315]  {$x_{1}$};
\draw (220,178.59) node [anchor=north west][inner sep=0.75pt]  [rotate=-315]  {$x_{2}$};
\draw (259,219) node [anchor=north west][inner sep=0.75pt]  [rotate=-315]  {$x_{3}$};
\draw (321.5,288.74) node [anchor=north west][inner sep=0.75pt]  [rotate=-315]  {$x_{n}$};
\draw (299.62,237.37) node [anchor=north west][inner sep=0.75pt]  [rotate=-45]  {$\cdots $};

\end{tikzpicture}
$$

\section{Future Directions}\label{future}
In this article we describe how a specific class of tangle-closure knots can be built from Thompson diagrams.  These are knots which are the closure of a product of positive $n$-crossing tangles, or the closure of a concatenation of $n$-crossing tangles.

The immediate next steps are to duplicate these procedures for products or concatenations consisting of a combination of negative and positive crossings; and tangle addition ought also to be considered.  After that, it would be desirable to understand the Thompson diagram analogs of tangle operations of multiplication, concatenation and addition in a truly general way -- for inputs that are arbitrary and not restricted to $n$-tangles.  

A more general question is whether there is a Markov-type theorem for links coming from Thompson group elements.  In the theory of braid groups, Markov's theorem explains that two braids have closures that are equivalent links if and only if they are related by a series of allowed moves.  These moves include some Reidemeister-type moves on the braids, as well as conjugation and stabilization.  It would be interesting to find a Markov-type relation on Thompson group elements, explaining when two distinct elements of $F$ generate the same link.  One may generate quite of lot of data about distinct elements that generate the same link by using the algorithms of Jones \cite{MR3589908}, of this article, and of Aiello-Baader \cite{MR4403999}, to find three different presentations of any knot coming from tangle multiplications or concatenations of positive $n$-crossing tangles.


\begin{bibdiv}
\begin{biblist}

\bib{MR4403999}{article}{
   author={Aiello, Valeriano},
   author={Baader, Sebastian},
   title={Arborescence of positive Thompson links},
   journal={Pacific J. Math.},
   volume={316},
   date={2022},
   number={2},
   pages={237--248},
   issn={0030-8730},
   review={\MR{4403999}},
   doi={10.2140/pjm.2022.316.237},
}

\bib{MR4294116}{article}{
   author={Aiello, Valeriano},
   author={Brothier, Arnaud},
   author={Conti, Roberto},
   title={Jones representations of Thompson's group $F$ arising from
   Temperley-Lieb-Jones algebras},
   journal={Int. Math. Res. Not. IMRN},
   date={2021},
   number={15},
   pages={11209--11245},
   issn={1073-7928},
   review={\MR{4294116}},
   doi={10.1093/imrn/rnz240},
}

\bib{MR4156130}{article}{
   author={Aiello, Valeriano},
   author={Jones, Vaughan F. R.},
   title={On spectral measures for certain unitary representations of R.
   Thompson's group $F$},
   journal={J. Funct. Anal.},
   volume={280},
   date={2021},
   number={1},
   pages={Paper No. 108777, 27},
   issn={0022-1236},
   review={\MR{4156130}},
   doi={10.1016/j.jfa.2020.108777},
}

\bib{MR4000616}{article}{
   author={Brothier, Arnaud},
   author={Jones, Vaughan F. R.},
   title={On the Haagerup and Kazhdan properties of R. Thompson's groups},
   journal={J. Group Theory},
   volume={22},
   date={2019},
   number={5},
   pages={795--807},
   issn={1433-5883},
   review={\MR{4000616}},
   doi={10.1515/jgth-2018-0114},
}

\bib{MR3989149}{article}{
   author={Brothier, Arnaud},
   author={Jones, Vaughan F. R.},
   title={Pythagorean representations of Thompson's groups},
   journal={J. Funct. Anal.},
   volume={277},
   date={2019},
   number={7},
   pages={2442--2469},
   issn={0022-1236},
   review={\MR{3989149}},
   doi={10.1016/j.jfa.2019.02.009},
}



\bib{MR1426438}{article}{
   author={Cannon, J. W.},
   author={Floyd, W. J.},
   author={Parry, W. R.},
   title={Introductory notes on Richard Thompson's groups},
   journal={Enseign. Math. (2)},
   volume={42},
   date={1996},
   number={3-4},
   pages={215--256},
   issn={0013-8584},
   review={\MR{1426438}},
}

\bib{MR2856142}{article}{
   author={Cannon, J. W.},
   author={Floyd, W. J.},
   title={What is $\ldots$ Thompson's group?},
   journal={Notices Amer. Math. Soc.},
   volume={58},
   date={2011},
   number={8},
   pages={1112--1113},
   issn={0002-9920},
   review={\MR{2856142}},
}
\bib{MR0258014}{article}{
   author={Conway, J. H.},
   title={An enumeration of knots and links, and some of their algebraic
   properties},
   conference={
      title={Computational Problems in Abstract Algebra},
      address={Proc. Conf., Oxford},
      date={1967},
   },
   book={
      publisher={Pergamon, Oxford},
   },
   date={1970},
   pages={329--358},
   review={\MR{0258014}},
}


\bib{MR3589908}{article}{
   author={Jones, Vaughan},
   title={Some unitary representations of Thompson's groups $F$ and $T$},
   journal={J. Comb. Algebra},
   volume={1},
   date={2017},
   number={1},
   pages={1--44},
   issn={2415-6302},
   review={\MR{3589908}},
   doi={10.4171/JCA/1-1-1},
}

\bib{MR3986040}{article}{
   author={Jones, Vaughan F. R.},
   title={On the construction of knots and links from Thompson's groups},
   conference={
      title={Knots, low-dimensional topology and applications},
   },
   book={
      series={Springer Proc. Math. Stat.},
      volume={284},
      publisher={Springer, Cham},
   },
   date={2019},
   pages={43--66},
   review={\MR{3986040}},
}


	

\bib{MR3764571}{article}{
   author={Jones, Vaughan F. R.},
   title={A no-go theorem for the continuum limit of a periodic quantum spin
   chain},
   journal={Comm. Math. Phys.},
   volume={357},
   date={2018},
   number={1},
   pages={295--317},
   issn={0010-3616},
   review={\MR{3764571}},
   doi={10.1007/s00220-017-2945-3},
}









\end{biblist}
\end{bibdiv}
\end{document}